\documentclass[12pt]{article}
\pdfoutput=1

\usepackage[utf8]{inputenc}

\usepackage{geometry,comment,color}
\usepackage{amsmath,amsthm,amssymb,mathrsfs,bbm,nicefrac,stmaryrd}

\usepackage[shortlabels,inline]{enumitem}
\usepackage[nosort,capitalise]{cleveref}

\crefname{enumi}{}{}
\crefname{equation}{}{}
\crefname{subsection}{Subsection}{Subsections}

\usepackage{xparse}
\NewDocumentEnvironment{cproof}{m}
{\begin{proof}[Proof of \cref{#1}]}%
{\noindent The proof of \cref{#1} is thus complete.
\end{proof}}

\NewDocumentEnvironment{cproof2}{m}
{\begin{proof}[Proof of \cref{#1}]}%
{the proof of \cref{#1}.
\end{proof}}

\newcommand{\N}{\ensuremath{\mathbb N}}

\newcommand{\R}{\ensuremath{\mathbb R}}

\newcommand{\bbB}{\ensuremath{\mathbb B}}

\renewcommand{\P}{\ensuremath{\mathbb P}}
\newcommand{\E}{\ensuremath{\mathbb E}}
\newcommand{\ud}{\,\mathrm{d}}
\newcommand{\dd}{\mathrm{d}}
\DeclareMathOperator*{\argmin}{arg\,min}

\newcommand{\cA}{\mathcal{A}}
\newcommand{\cB}{\mathcal{B}}
\newcommand{\cC}{\mathcal{C}}

\newcommand{\cE}{\mathcal{E}}
\newcommand{\cF}{\mathcal{F}}
\newcommand{\cG}{\mathcal{G}}

\newcommand{\cL}{\mathcal{L}}
\newcommand{\cM}{\mathcal{M}}
\newcommand{\cN}{\mathcal{N}}

\newcommand{\cR}{\mathcal{R}}
\newcommand{\cS}{\mathcal{S}}

\newcommand{\cX}{\mathcal{X}}
\newcommand{\cY}{\mathcal{Y}}

\newcommand{\fC}{\mathfrak{C}}

\newcommand{\fK}{\mathfrak{K}}
\newcommand{\fL}{\mathfrak{L}}
\newcommand{\fM}{\mathfrak{M}}
\newcommand{\fN}{\mathfrak{N}}

\newcommand{\fR}{\mathfrak{R}}

\newcommand{\fc}{\mathfrak{c}}
\newcommand{\fd}{\mathfrak{d}}

\newcommand{\fn}{\mathfrak{n}}

\newcommand{\fr}{\mathfrak{r}}

\newcommand{\fx}{\mathfrak{x}}

\newcommand{\bfa}{\mathbf{a}}

\newcommand{\bfd}{\mathbf{d}}

\newcommand{\bfk}{\mathbf{k}}
\newcommand{\bfl}{\mathbf{l}}
\newcommand{\bfm}{\mathbf{m}}

\newcommand{\bfJ}{\mathbf{J}}

\newcommand{\bfL}{\mathbf{L}}

\newcommand{\bfN}{\mathbf{N}}

\newcommand{\bfR}{\mathbf{R}}

\newcommand{\scrD}{\mathscr{D}}
\newcommand{\scrE}{\mathscr{E}}
\newcommand{\scrF}{\mathscr{F}}

\newcommand{\scrN}{\mathscr{N}}

\newcommand{\scrR}{\mathscr{R}}

\newcommand{\clippedNN}[4]{\mathscr{N}^{#1,#2}_{#3,#4}}

\newcommand{\ltriplevert}{\mathopen{\lvert\kern-0.25ex\lvert\kern-0.25ex\lvert}}
\newcommand{\rtriplevert}{\rvert\kern-0.25ex\rvert\kern-0.25ex\rvert}
\newcommand{\biggltriplevert}{\mathopen{\biggl\lvert\kern-0.25ex\biggl\lvert\kern-0.25ex\biggl\lvert}}
\newcommand{\biggrtriplevert}{\biggr\rvert\kern-0.25ex\biggr\rvert\kern-0.25ex\biggr\rvert}

\DeclareFontFamily{U}{mathb}{\hyphenchar\font45}
\DeclareFontShape{U}{mathb}{m}{n}{
      <5> <6> <7> <8> <9> <10> gen * mathb
      <10.95> mathb10 <12> <14.4> <17.28> <20.74> <24.88> mathb12
      }{}
\DeclareSymbolFont{mathb}{U}{mathb}{m}{n}
\DeclareMathSymbol{\llfloor}{\mathopen}{mathb}{'172}
\DeclareMathSymbol{\rrfloor}{\mathclose}{mathb}{'173}

\DeclareMathOperator*{\smallsum}{\textstyle\sum}
\DeclareMathOperator*{\smallprod}{\textstyle\prod}
\DeclareMathOperator*{\smallbigcup}{\textstyle\bigcup}

\makeatletter
\newcommand{\pushright}[1]{\ifmeasuring@#1\else\omit\hfill$\displaystyle#1$\fi\ignorespaces}
\makeatother

\usepackage{scalerel}
\newcommand{\medint}[1]{{\stretchrel*{\scalerel*[2ex]{\int}{\int}}{\sum}}_{ \mkern-10mu \smash{#1} }}

\newtheorem{theorem}{Theorem}[section]
\newtheorem{definition}[theorem]{Definition}
\newtheorem{proposition}[theorem]{Proposition}
\newtheorem{corollary}[theorem]{Corollary}
\newtheorem{lemma}[theorem]{Lemma}

\hypersetup{colorlinks=true}
\newcommand\yesnumber{\refstepcounter{equation}\tag{\theequation}}

\usepackage[backend=biber,style=numeric-comp,giveninits=true,sorting=nyt,maxbibnames=99,useprefix=true]{biblatex}
\DeclareNameAlias{default}{last-first}
\renewbibmacro{in:}{}
\AtBeginBibliography{%

}
\DefineBibliographyStrings{english}{%
	andothers = {\addcomma\addspace\textsc{et\addabbrvspace al}\adddot},
	and = {\textsc{and}}
}

\DeclareFieldFormat
[article,inbook,incollection,inproceedings,patent,thesis,unpublished]
{title}{#1}
\renewbibmacro{in:}{%
	\ifentrytype{article}{%
	}{%
	\printtext{\bibstring{in}\intitlepunct}%
}%
}
\renewbibmacro*{volume+number+eid}{%
	\printfield{volume}%
	\setunit*{\addcomma\space}%
	\printfield{number}%
	\setunit{\addcomma\space}%
	\printfield{eid}}
\DeclareFieldFormat{pages}{#1}
\renewbibmacro*{publisher+location+date}{%
	\printlist{publisher}%
	\setunit*{\addcomma\space}%
	\printlist{location}%
	\setunit*{\addcomma\space}%
	\usebibmacro{date}%
	\newunit}

\renewbibmacro*{doi+eprint+url}{%
  \iftoggle{bbx:doi}
    {\iffieldundef{url}{\printfield{doi}}{}}
    {}%
  \newunit\newblock
  \iftoggle{bbx:eprint}
    {\usebibmacro{eprint}}
    {}%
  \newunit\newblock
  \iftoggle{bbx:url}
    {\usebibmacro{url+urldate}}
    {}}

\title{Overall error analysis for the training
of\\ deep neural networks via
stochastic\\ gradient descent with random initialisation}

\author{Arnulf Jentzen$^{1}$ and Timo Welti$^{2}$\bigskip\\
\small{$^1$ Faculty of Mathematics and Computer Science, University of M\"unster,}\\
\small{M\"unster, Germany; e-mail: \texttt{ajentzen}\textcircled{\texttt{a}}\texttt{uni-muenster.de}}\\
\small{$^2$ SAM, Department of Mathematics, ETH Z\"urich,}\\
\small{Z\"urich, Switzerland; e-mail: \texttt{twelti}\textcircled{\texttt{a}}\texttt{twelti.org}}
}

\begin{document}

\maketitle

\begin{abstract}
In spite of
the accomplishments of deep learning based algorithms in numerous applications
and very broad corresponding research interest,
at the moment there is still no rigorous understanding of the reasons why
such algorithms produce useful results in certain situations.
A thorough mathematical analysis of deep learning based algorithms
seems to be crucial in order
to improve our understanding
and to make their implementation more effective and efficient. 
In this article we provide
a mathematically rigorous full error analysis
of deep learning based empirical risk minimisation
with quadratic loss function
in the probabilistically strong sense,
where the underlying deep neural networks
are trained using stochastic gradient descent
with random initialisation.
The convergence speed
we obtain
is presumably far from optimal
and
suffers under the curse of dimensionality.
To the best of our knowledge,
we establish, however,
the first full error analysis in the scientific literature
for a deep learning based algorithm
in the probabilistically strong sense
and, moreover,
the first full error analysis in the scientific literature
for a deep learning based algorithm
where stochastic gradient descent
with random initialisation
is the employed optimisation method.
\end{abstract}

\begin{center}
full error analysis,
approximation, generalisation, optimisation,
strong\\ convergence,
stochastic gradient descent,
random initialisation
\end{center}

\pagebreak

\tableofcontents

\pagebreak

\section{Introduction}

Deep learning based algorithms have been applied
extremely successfully to overcome fundamental challenges
in many different areas, such as
image recognition, natural language processing, game intelligence, autonomous driving,
and computational advertising,
just to name a few.
In line with this,
researchers from a wide range of different fields,
including, for example, computer science, mathematics, chemistry, medicine, and finance,
are investing significant efforts into studying such algorithms
and employing them to tackle challenges arising in their fields.
In spite of this broad research interest
and the accomplishments of deep learning based algorithms in numerous applications,
at the moment there is still no rigorous understanding of the reasons why such algorithms
produce useful results in certain situations.
Consequently,
there is no rigorous way to predict,
before actually implementing a deep learning based algorithm,
in which situations it
might perform reliably and in which situations it might fail.
This necessitates in many cases a trial-and-error approach
in order to move forward,
which can cost a lot of time and resources.
A thorough mathematical analysis of deep learning based algorithms
(in scenarios where it is possible to formulate such an analysis)
seems to be crucial in order to make progress on these issues.
Moreover, such an analysis may lead to new insights that enable the design of more effective and efficient algorithms.

The aim of this article
is to provide a mathematically rigorous full error analysis
of deep learning based empirical risk minimisation
with quadratic loss function
in the probabilistically strong sense,
where the underlying deep neural networks (DNNs)
are trained using stochastic gradient descent (SGD)
with random initialisation (cf.\ \cref{thm:intro} below).
For a brief illustration of deep learning based empirical risk minimisation
with quadratic loss function,
consider
natural numbers $ d, \bfd \in \N $,
a probability space
$ ( \Omega, \cF, \P ) $,
random variables
$ X \colon \Omega \to [ 0, 1 ]^d $
and
$ Y \colon \Omega \to [ 0, 1 ] $,
and a measurable function
$ \cE \colon [ 0, 1 ]^d \to [ 0, 1 ] $
satisfying
$ \P $-a.s.\ that
$ \cE( X ) = \E[ Y \vert X ] $.
The goal is
to find a DNN with appropriate architecture
and appropriate parameter vector $ \theta \in \R^\bfd $
(collecting its weights and biases)
such that its realisation
$ \scrN_\theta \colon \R^d \to \R $
approximates the target function $ \cE $ well
in the sense that
the error
$ \E[ \lvert \scrN_\theta( X ) - \cE( X ) \rvert^p ]
=
\int_{ \smash{ [ 0, 1 ]^d } }
    \lvert \scrN_\theta( x ) - \cE( x ) \rvert^p
\, \P_{ X }( \dd x )
\in [ 0, \infty ) $
for some $ p \in [ 1, \infty ) $
is as small as possible.
In other words,
given $ X $
we want $ \scrN_\theta( X ) $ to predict $ Y $
as reliably as possible.
Due to the well-known
bias--variance decomposition
(cf., e.g.,
Beck, Jentzen, \& Kuckuck~\cite[Lemma~4.1]{BeckJentzenKuckuck2019arXiv}),
for the case $ p = 2 $
minimising the error function
$ \R^\bfd \ni \theta \mapsto
\E[ \lvert \scrN_\theta( X ) - \cE( X ) \rvert^2 ]
\in [ 0, \infty ) $
is equivalent to minimising the risk function
$ \R^\bfd \ni \theta \mapsto
\E[ \lvert \scrN_\theta( X ) - Y \rvert^2 ]
\in [ 0, \infty ) $
(corresponding to a quadratic loss function).
Since in practice the joint distribution of $ X $ and $ Y $ is typically not known,
the risk function is replaced
by an empirical risk function based on i.i.d.\ training samples of $ ( X, Y ) $.
This empirical risk
is then approximatively minimised using
an optimisation method such as SGD.
As is often the case for deep learning based algorithms,
the overall error arising from this procedure consists of the following three different parts
(cf.~\cite[Lemma~4.3]{BeckJentzenKuckuck2019arXiv}
and \cref{prop:error_decomposition} below):
\begin{enumerate*}[(i)]
\item
the \emph{approximation error}
(cf., e.g.,~\cite{Cybenko1989,
Funahashi1989,
HornikStinchcombeWhite1989,
HornikStinchcombeWhite1990,
HartmanKeelerKowalski1990,
BlumLi1991,
ParkSandberg1991,
Hornik1991,
Hornik1993,
LeshnoLinPinkusSchocken1993,
Barron1993,
Barron1994,
ChuiMhaskar1994,
Ellacott1994}
and the references in the introductory paragraph in \cref{sec:approximation_error}),
which arises from approximating the target function $ \cE $
by the considered class of DNNs,
\item
the \emph{generalisation error}
(cf., e.g.,~\cite{VanDeGeer2000,
CuckerSmale2002,
GyorfiKohlerKrzyzakWalk2002,
BartlettBousquetMendelson2005,
Massart2007,
Shalev_ShwartzBen_David2014,
BernerGrohsJentzen2018arXiv,
BeckJentzenKuckuck2019arXiv,
EMaWu2019,
EMaWang2019arXiv,
EMaWu2020online}),
which arises from replacing the true risk by the empirical risk,
and
\item
the \emph{optimisation error}
(cf., e.g.,~\cite{AroraDuHuLiWang2019,
BachMoulines2013,
BeckBeckerGrohsJaafariJentzen2018arXiv,
BeckJentzenKuckuck2019arXiv,
BercuFort2013,
ChauMoulinesRasonyiSabanisZhang2019arXiv,
DereichKassing2019arXiv,
DereichMuller_Gronbach2019,
DuLeeLiWangZhai2019,
DuZhaiPoczosSingh2018arXiv,
FehrmanGessJentzen2019arXiv,
JentzenKuckuckNeufeldVonWurstemberger2018arXiv,
JentzenvonWurstemberger2020,
KarimiMiasojedowMoulinesWai2019arXiv,
LeiHuLiTang2019arXiv,
Shamir2019,
ZhangMartensGrosse2019,
ZouCaoZhouGu2019}),
which arises from computing only an approximate minimiser
using the selected optimisation method.
\end{enumerate*}

In this work we 
derive strong convergence rates for
the approximation error,
the generalisation error,
and the optimisation error
separately
and combine these findings
to establish strong convergence results
for the overall error
(cf.\ \cref{sec:overall_error_strong_rate,sec:SGD}),
as illustrated in \cref{thm:intro} below.
The convergence speed
we obtain (cf.~\cref{eq:thm:intro} in \cref{thm:intro})
is presumably far from optimal,
suffers under the curse of dimensionality
(cf., e.g.,
Bellman~\cite{Bellman1957}
and
\begingroup\renewcommand{\mkbibbrackets}[1]{[#1}%
Novak \& Wo\'{z}niakowski~\cite[Chapter~1]{NovakWozniakowski2008};
\renewcommand{\mkbibbrackets}[1]{#1]}%
\cite[Chapter~9]{NovakWozniakowski2010}\endgroup),
and is, as a consequence, very slow.
To the best of our knowledge,
\cref{thm:intro}
is, however,
the first full error result in the scientific literature
for a deep learning based algorithm
in the probabilistically strong sense
and, moreover,
the first full error result in the scientific literature
for a deep learning based algorithm
where SGD with random initialisation
is the employed optimisation method.
We now present \cref{thm:intro},
the statement of which is entirely self-contained,
before we add further
explanations and intuitions
for the mathematical objects that are introduced.

\begin{theorem}
\label{thm:intro}
Let
$ d, \bfd, \bfL, \bfJ, M, K, N \in \N $,
$ \gamma, L \in \R $,
$ c \in [ \max\{ 2, L \}, \infty ) $,
$ \bfl = ( \bfl_0, \ldots, \bfl_\bfL )
\allowbreak
\in \N^{ \bfL + 1 } $,
$ \bfN \subseteq \{ 0, \ldots, N \} $,
assume
$ 0 \in \bfN $,
$ \bfl_0 = d $,
$ \bfl_\bfL = 1 $,
and
$ \bfd \geq \sum_{i=1}^{\bfL} \bfl_i( \bfl_{ i - 1 } + 1 ) $,
for every
$ m, n \in \N $,
$ s \in \N_0 $,
$ \theta = ( \theta_1, \ldots, \theta_\bfd ) \in \R^\bfd $
with
$ \bfd \geq s + m n + m $
let
$ \cA_{ m, n }^{ \theta, s }
\colon \R^n \to \R^m $
satisfy for all
$ x = ( x_1, \ldots, x_n ) \in \R^n $
that
\begin{equation}
\label{eq:affine_linear}
\cA_{ m, n }^{ \theta, s }( x )
= 
\begin{pmatrix}
\theta_{ s + 1 }
& \theta_{ s + 2 }
& \cdots
& \theta_{ s + n }
\\
\theta_{ s + n + 1 }
& \theta_{ s + n + 2 }
& \cdots
& \theta_{ s + 2 n }
\\
\vdots
& \vdots
& \ddots
& \vdots
\\
\theta_{ s + ( m - 1 ) n + 1 }
& \theta_{ s + ( m - 1 ) n + 2 }
& \cdots
& \theta_{ s + m n }
\end{pmatrix}
\begin{pmatrix}
x_1
\\
x_2
\\
\vdots
\\
x_n
\end{pmatrix}
+
\begin{pmatrix}
\theta_{ s + m n + 1 }
\\
\theta_{ s + m n + 2 }
\\
\vdots 
\\
\theta_{ s + m n + m }
\end{pmatrix}
,
\end{equation}
let
$ \bfa_i \colon \R^{ \bfl_i } \to \R^{ \bfl_i } $,
$ i \in \{ 1, \ldots, \bfL \} $,
satisfy for all
$ i \in \N \cap [ 0, \bfL ) $,
$ x = ( x_1, \ldots, x_{ \bfl_i } ) \in \R^{ \bfl_i } $
that
$ \bfa_i( x ) =
( \max\{ x_1, 0 \},
\ldots,
\max\{ x_{ \bfl_i }, 0 \} ) $,
assume for all
$ x \in \R $
that
$ \bfa_\bfL( x )
= \max\{ \min\{ x, 1 \}, 0 \} $,
for every
$ \theta \in \R^\bfd $
let
$ \scrN_\theta \colon \R^d \to \R $
satisfy
$ \scrN_\theta
=
\bfa_\bfL \circ \cA_{ \bfl_\bfL, \bfl_{ \bfL - 1 } }^{ \theta, \smash{ \sum_{i = 1}^{\bfL-1} \bfl_i ( \bfl_{ i - 1 } + 1 ) } }
\circ \bfa_{ \bfL - 1 } \circ \smash{ \cA_{ \bfl_{ \bfL - 1 },  \bfl_{ \bfL - 2 } }^{ \theta, \smash{ \sum_{i = 1}^{\bfL-2} \bfl_i ( \bfl_{ i - 1 } + 1 ) } } }
\circ \ldots
\circ \bfa_1 \circ \cA_{  \bfl_1, \bfl_0 }^{ \theta, 0 } $,
let
$ ( \Omega, \cF, \P ) $
be a probability space,
let
$ X^{ k, n }_{ \smash{j} }
\colon \Omega \to [ 0, 1 ]^d $,
$ k, n, j \in \N_0 $,
and
$ Y^{ k, n }_{ \smash{j} }
\colon \Omega \to [ 0, 1 ] $,
$ k, n, j \in \N_0 $,
be functions,
assume that
$ ( X^{ 0, 0 }_{ \smash{j} }, Y^{ 0, 0 }_{ \smash{j} } ) $,
$ j \in \N $,
are i.i.d.\ random variables,
let
$ \cE \colon [ 0, 1 ]^d \to [ 0, 1 ] $
satisfy
$ \P $-a.s.\
that
$ \cE( X_{ \smash{1} }^{ 0, 0 } )
= \E[ Y_{ \smash{1} }^{ 0, 0 } \vert X_{ \smash{1} }^{ 0, 0 } ] $,
assume for all
$ x, y \in [ 0, 1 ]^d $
that
$ \lvert \cE( x ) - \cE( y ) \rvert \leq L \lVert x - y \rVert_{ 1 } $,
let
$ \Theta_{ k, n } \colon \Omega \to \R^{\bfd} $,
$ k, n \in \N_0 $,
and
$ \bfk \colon \Omega \to ( \N_0 )^2 $
be random variables,
assume
$ \bigl( \bigcup_{ k = 1 }^{ \infty }
\Theta_{ k, 0 }( \Omega ) \bigr)
\subseteq [ -c, c ]^\bfd $,
assume that
$ \Theta_{ k, 0 } $,
$ k \in \N $,
are i.i.d.,
assume that
$ \Theta_{ 1, 0 } $ is continuous uniformly distributed on $ [ -c, c ]^\bfd $,
let
$ \cR^{ k, n }_{ \smash{J} } \colon \R^{\bfd} \times \Omega \to [ 0, \infty ) $,
$ k, n, J \in \N_0 $,
and
$ \cG^{ k, n } \colon \R^{\bfd} \times \Omega \to \R^{\bfd} $,
$ k, n \in \N $,
satisfy for all
$ k, n \in \N $,
$ \omega \in \Omega $,
$ \theta \in
\{ \vartheta \in \R^{\bfd} \colon
( \cR^{ k, n }_{ \smash{ \bfJ } } ( \cdot, \omega ) \colon \R^{\bfd} \to [ 0, \infty )
\text{ is differentiable at } \vartheta ) \} $
that
$ \cG^{ k, n }( \theta, \omega )
=
( \nabla_\theta \cR^{ k, n }_{ \smash{ \bfJ } } )( \theta, \omega ) $,
assume for all
$ k, n \in \N $
that
$ \Theta_{ k, n } = \Theta_{ k, n -1 } - \gamma \cG^{ k, n }( \Theta_{ k, n -1 } ) $,
and
assume for all
$ k, n \in \N_0 $,
$ J \in \N $,
$ \theta \in \R^{\bfd} $,
$ \omega \in \Omega $
that
\begin{align}
\label{eq:empirical_risk}
\cR^{ k, n }_{ \smash{J} }( \theta, \omega )
=
\frac{1}{J}
\biggl[
\smallsum_{j=1}^{J}
    \lvert \scrN_\theta( X^{ k, n }_{ \smash{j} }( \omega ) ) - Y^{ k, n }_{ \smash{j} }( \omega ) \rvert^2
\biggr]
\qquad
\text{and}
\\
\label{eq:choice_bfk}
\bfk( \omega ) \in
\argmin\nolimits_{ ( l, m ) \in \{ 1, \ldots, K \} \times \bfN, \, \lVert \Theta_{ l, m }( \omega ) \rVert_{ \infty } \leq c }
\cR^{ 0, 0 }_{ \smash{ M } }( \Theta_{ l, m }( \omega ), \omega )
.
\end{align}
Then
\begin{equation}
\label{eq:thm:intro}
\begin{split}
&
\E\Bigl[
\medint{[ 0, 1 ]^d}
        \lvert \scrN_{ \Theta_{\bfk} }( x ) - \cE( x ) \rvert
    \, \P_{ X^{ 0, 0 }_{ \smash{1}\vphantom{x} } }( \dd x )
\Bigr]
\\ &
\leq
\frac{
d c^3
}{
[ \min\{ \bfL, \bfl_1, \ldots, \bfl_{ \bfL - 1 } \} ]^{ \nicefrac{1}{d} }
}
+
\frac{
    c^3
    \bfL ( \lVert \bfl \rVert_{ \infty } + 1 )
    \ln( e M )
}{ M^{ \nicefrac{1}{4} } }
+
\frac{
\bfL
( \lVert \bfl \rVert_{ \infty } + 1 )^\bfL
c^{ \bfL + 1 }
}{
K^{ [ ( 2 \bfL )^{-1} ( \lVert \bfl \rVert_{ \infty } + 1 )^{-2} ] }
}
.
\end{split}
\end{equation}
\end{theorem}

\noindent
Recall that
we denote
for every $ p \in [ 1, \infty ] $
by
$ \lVert \cdot \rVert_p \colon \bigl( \bigcup_{n=1}^\infty \R^n \bigr) \to [ 0, \infty ) $
the
$ p $-norm
of vectors in $ \bigcup_{n=1}^\infty \R^n $
(cf.~\cref{def:p-norm}).
In addition,
note that
the function
$ \Omega \times [ 0, 1 ]^d \ni ( \omega, x )
\mapsto
\lvert \scrN_{ \smash{ \Theta_{ \smash{ \bfk( \omega ) } }( \omega ) } }( x ) - \cE( x ) \rvert
\in [ 0, \infty ) $
is measurable
(cf.\ \cref{lem:measurability})
and that the expression on the left hand side of~\cref{eq:thm:intro} above is thus well-defined.
\cref{thm:intro} follows directly from \cref{cor:SGD_simplfied} in \cref{sec:SGD},
which, in turn, is a consequence of the main result of this article,
\cref{thm:main} in \cref{sec:overall_error_strong_rate}.

In the following we provide additional explanations and intuitions
for \cref{thm:intro}.
For every
$ \theta \in \R^\bfd $
the functions
$ \scrN_\theta \colon \R^d \to \R $
are realisations
of fully connected feedforward artificial neural networks
with $ \bfL + 1 $ layers
consisting of an input layer of dimension $ \bfl_0 = d $,
of $ \bfL - 1 $ hidden layers
of dimensions $ \bfl_1, \ldots, \bfl_{ \bfL - 1 } $, respectively,
and
of an output layer of dimension $ \bfl_{ \bfL } = 1 $
(cf.\ \cref{def:clipped_NN}).
The weights and biases stored in the DNN parameter vector $ \theta \in \R^\bfd $
determine the corresponding $ \bfL $ affine linear transformations
(cf.~\cref{eq:affine_linear} above).
As activation functions
we employ the multidimensional versions
$ \bfa_1, \ldots, \bfa_{ \bfL - 1 } $
(cf.\ \cref{def:multidimensional_version})
of the \emph{rectifier function}
$ \R \ni x \mapsto \max\{ x, 0 \} \in \R $
(cf.\ \cref{def:ReLu})
just in front of each of the hidden layers
and
the \emph{clipping function}
$ \bfa_\bfL $
(cf.\ \cref{def:clipping_function})
just in front of the output layer.
Furthermore,
observe that we assume the target function
$ \cE \colon [ 0, 1 ]^d \to [ 0, 1 ] $,
the values of which we intend to approximately predict with the trained DNN,
to be Lipschitz continuous with Lipschitz constant $ L $.
Moreover,
for every
$ k, n \in \N_0 $,
$ J \in \N $
the function
$ \cR^{ k, n }_{ \smash{J} } \colon \R^{\bfd} \times \Omega \to [ 0, \infty ) $
is the empirical risk
based on the $ J $ training samples
$ ( X^{ k, n }_{ \smash{j} }, Y^{ k, n }_{ \smash{j} } ) $,
$ j \in \{ 1, \ldots, J \} $
(cf.~\cref{eq:empirical_risk} above).
Derived from the empirical risk,
for every
$ k, n \in \N $
the function
$ \cG^{ k, n } \colon \R^{\bfd} \times \Omega \to \R^{\bfd} $
is a (generalised) gradient
of the empirical risk $ \cR^{ k, n }_{ \smash{\bfJ} } $
with respect to its first argument,
that is, with respect to the DNN parameter vector $ \theta \in \R^\bfd $.
These gradients are required
in order to formulate the training dynamics
of the (random) DNN parameter vectors
$ \Theta_{ k, n } \in \R^\bfd $,
$ k \in \N $,
$ n \in \N_0 $,
given by the SGD optimisation method
with learning rate $ \gamma $.
Note that the subscript $ n \in \N_0 $
of these SGD iterates (i.e., DNN parameter vectors)
is the current training step number,
whereas the subscript $ k \in \N $
counts the number of times the SGD iteration has been started from scratch so far.
Such a new start entails the corresponding initial DNN parameter vector
$ \Theta_{ k, 0 } \in \R^\bfd $
to be drawn
continuous uniformly from the hypercube $ [ -c, c ]^\bfd $,
in accordance with Xavier initialisation
(cf.\ Glorot \& Bengio~\cite{GlorotBengio2010}).
The (random) double index $ \bfk \in \N \times \N_0 $
represents the final choice made for the DNN parameter vector
$ \Theta_\bfk \in \R^\bfd $
(cf.~\cref{eq:thm:intro} above),
concluding the training procedure,
and is selected as follows.
During training
the empirical risk
$ \cR^{ 0, 0 }_{ \smash{ M } } $
has been calculated for the subset
of the SGD iterates
indexed by
$ \bfN \subseteq \{ 0, \ldots, N \} $
provided that they have not left the hypercube $ [ -c, c ]^\bfd $
(cf.~\cref{eq:choice_bfk} above).
After the SGD iteration has been started and finished $ K $ times (with maximally $ N $ training steps in each case)
the final choice for the DNN parameter vector $ \Theta_\bfk \in \R^\bfd $
is made among those SGD iterates
for which the calculated empirical risk is minimal
(cf.~\cref{eq:choice_bfk} above).
Observe that
we use mini-batches of size $ \bfJ $
consisting,
during SGD iteration number $ k \in \{ 1, \ldots, K \} $
for training step number $ n \in \{ 1, \ldots, N \} $,
of the training samples
$ ( X^{ k, n }_{ \smash{j} }, Y^{ k, n }_{ \smash{j} } ) $,
$ j \in \{ 1, \ldots, \bfJ \} $,
and that
we reserve
the $ M $ training samples
$ ( X^{ 0, 0 }_{ \smash{j} }, Y^{ 0, 0 }_{ \smash{j} } ) $,
$ j \in \{ 1, \ldots, M \} $,
for checking the value of
the empirical risk
$ \cR^{ 0, 0 }_{ \smash{ M } } $.
Regarding the conclusion of \cref{thm:intro},
note that
the left hand side of~\cref{eq:thm:intro}
is the expectation of the overall $ L^1 $-error,
that is,
the expected $ L^1 $-distance
between the trained DNN $ \scrN_{ \Theta_{ \smash{ \bfk } } } $ and the target function $ \cE $.
It is bounded from above by the
right hand side of~\cref{eq:thm:intro},
which consists of following three summands:
\begin{enumerate*}[(i)]
\item
the first summand corresponds to the
\emph{approximation error}
and converges to zero as
the number of hidden layers $ \bfL - 1 $
as well as the hidden layer dimensions
$ \bfl_1, \ldots, \bfl_{ \bfL - 1 } $
increase to infinity,
\item
the second summand corresponds to the
\emph{generalisation error}
and converges to zero as
number of training samples $ M $
used for calculating the empirical risk
increases to infinity,
and
\item
the third summand corresponds to the
\emph{optimisation error}
and converges to zero as
total number of times $ K $
the SGD iteration has been started from scratch
increases to infinity.
\end{enumerate*}
We would like to point out that the
the second summand (corresponding to the generalisation error)
does not suffer under the curse of dimensionality with respect to any of the variables involved.

The main result of this article,
\cref{thm:main} in \cref{sec:overall_error_strong_rate},
covers,
in comparison with \cref{thm:intro},
the more general cases
where
$ L^p $-norms of the overall $ L^2 $-error
instead of the expectation of the overall $ L^1 $-error
are considered
(cf.~\cref{eq:thm:main} in \cref{thm:main}),
where the training samples are not restricted to unit hypercubes,
and
where a general stochastic approximation algorithm
(cf., e.g., Robbins \& Monro~\cite{RobbinsMonro1951})
with random initialisation
is used for optimisation.
Our convergence proof for the optimisation error
relies, in fact, on the convergence of the Minimum Monte Carlo method
(cf.\ \cref{prop:minimum_MC_rate} in \cref{sec:optimisation_error})
and thus only exploits random initialisation
but not the specific dynamics of the employed optimisation method
(cf.~\cref{eq:reduce_to_minimum_MC} in the proof of \cref{prop:main}).
In this regard,
note that \cref{thm:intro} above also includes the application of
deterministic gradient descent instead of SGD
for optimisation
since we do not assume the samples used for gradient iterations to be i.i.d.
Parts of our derivation of
\cref{thm:intro} and \cref{thm:main}, respectively,
are inspired by
Beck, Jentzen, \& Kuckuck~\cite{BeckJentzenKuckuck2019arXiv},
Berner, Grohs, \& Jentzen~\cite{BernerGrohsJentzen2018arXiv},
and
Cucker \& Smale~\cite{CuckerSmale2002}.

This article is structured in the following way.
\cref{sec:DNNs}
recalls some basic definitions related to DNNs
and thereby introduces the corresponding notation
we use in the subsequent parts of this article.
In \cref{sec:approximation_error}
we examine the approximation error
and, in particular,
establish a convergence result for the approximation of
Lipschitz continuous functions by DNNs.
The following section,
\cref{sec:generalisation_error},
contains our strong convergence analysis
of the generalisation error.
In \cref{sec:optimisation_error},
in turn,
we address the optimisation error
and derive in connection with this strong convergence rates for the Minimum Monte Carlo method.
Finally,
we combine in
\cref{sec:DL_risk_minimisation}
a decomposition of the overall error
(cf.\ \cref{sec:overall_error})
with our results for the different error sources
from \cref{sec:approximation_error,sec:generalisation_error,sec:optimisation_error}
to prove
strong convergence results for the overall error.
The employed optimisation method
is initially allowed to be a general stochastic approximation algorithm with random initialisation
(cf.\ \cref{sec:overall_error_strong_rate})
and is afterwards specialised to the setting of SGD with random initialisation
(cf.\ \cref{sec:SGD}).

\section{Basics on deep neural networks (DNNs)}
\label{sec:DNNs}

In this section we present the mathematical description of DNNs
which we use throughout the remainder of this article.
It is a vectorised description in the sense that
all the weights and biases associated to the DNN under consideration
are collected in a single parameter vector
$ \theta \in \R^\bfd $
with $ \bfd \in \N $ sufficiently large
(cf.\ \cref{def:NN,def:clipped_NN}).
The content of this section is taken from
Beck, Jentzen, \& Kuckuck~\cite[Section~2.1]{BeckJentzenKuckuck2019arXiv}
and
is based on well-known material from the scientific literature,
see, e.g.,
Beck et al.~\cite{BeckBeckerGrohsJaafariJentzen2018arXiv},
Beck, E, \& Jentzen~\cite{BeckEJentzen2019},
Berner, Grohs, \& Jentzen~\cite{BernerGrohsJentzen2018arXiv},
E, Han, \& Jentzen~\cite{EHanJentzen2017},
Goodfellow, Bengio, \& Courville~\cite{GoodfellowBengioCourville2016},
and
Grohs et al.~\cite{GrohsHornungJentzenZimmermann2019arXiv}. 
In particular,
\cref{def:affine_function}
is \cite[Definition~2.1]{BeckJentzenKuckuck2019arXiv}
(cf., e.g., (25) in \cite{BeckEJentzen2019}),
\cref{def:NN}
is \cite[Definition~2.2]{BeckJentzenKuckuck2019arXiv}
(cf., e.g., (26) in \cite{BeckEJentzen2019}),
\cref{def:multidimensional_version}
is \cite[Definition~2.3]{BeckJentzenKuckuck2019arXiv}
(cf., e.g., \cite[Definition~2.2]{GrohsHornungJentzenZimmermann2019arXiv}),
and
\cref{def:ReLu,def:ReLu_multidim,def:clipping_function,def:clipping_multidim,def:clipped_NN}
are \cite[Definitions~2.4, 2.5, 2.6, 2.7, and 2.8]{BeckJentzenKuckuck2019arXiv}
(cf., e.g., \cite[Setting~2.5]{BernerGrohsJentzen2018arXiv}
and
\cite[Section~6.3]{GoodfellowBengioCourville2016}).

\subsection{Vectorised description of DNNs}
\label{sec:vectorised_description}

\begin{definition}[Affine function]
\label{def:affine_function}
Let
$ \bfd, m, n \in \N $,
$ s \in \N_0 $,
$ \theta = ( \theta_1, \theta_2, \ldots, \theta_\bfd ) \in \R^\bfd $
satisfy
$ \bfd \geq s + m n + m $.
Then
we denote by
$ \cA_{ m, n }^{ \theta, s }
\colon \R^n \to \R^m $
the function which satisfies for all
$ x = ( x_1, x_2, \ldots, x_n ) \in \R^n $
that
\begin{align*}
\yesnumber
& \cA_{ m, n }^{ \theta, s }( x )
= 
\begin{pmatrix}
\theta_{ s + 1 }
& \theta_{ s + 2 }
& \cdots
& \theta_{ s + n }
\\
\theta_{ s + n + 1 }
& \theta_{ s + n + 2 }
& \cdots
& \theta_{ s + 2 n }
\\
\theta_{ s + 2 n + 1 }
& \theta_{ s + 2 n + 2 }
& \cdots
& \theta_{ s + 3 n }
\\
\vdots
& \vdots
& \ddots
& \vdots
\\
\theta_{ s + ( m - 1 ) n + 1 }
& \theta_{ s + ( m - 1 ) n + 2 }
& \cdots
& \theta_{ s + m n }
\end{pmatrix}
\begin{pmatrix}
x_1
\\
x_2
\\
x_3
\\
\vdots
\\
x_n
\end{pmatrix}
+
\begin{pmatrix}
\theta_{ s + m n + 1 }
\\
\theta_{ s + m n + 2 }
\\
\theta_{ s + m n + 3 }
\\
\vdots 
\\
\theta_{ s + m n + m }
\end{pmatrix}
\\ &
=
\biggl( \!
\biggl[
    \smallsum_{ i = 1 }^n \theta_{ s + i } x_i
\biggr]
+ \theta_{ s + m n + 1 },
\biggl[
    \smallsum_{ i = 1 }^n \theta_{ s + n + i } x_i
\biggr]
+ \theta_{ s + m n + 2 },
\ldots,
\biggl[
    \smallsum_{i = 1}^n \theta_{ s + ( m - 1 ) n + i } x_i
\biggr]
+ \theta_{ s + m n + m }
\! \biggr)
.
\end{align*}
\end{definition}

\begin{definition}[Fully connected feedforward artificial neural network]
\label{def:NN}
Let
$ \bfd, \bfL, \bfl_0, \bfl_1, \ldots,
\allowbreak
\bfl_\bfL \in \N $,
$ s \in \N_0 $,
$ \theta \in \R^\bfd $
satisfy
$ \bfd \geq s + \sum_{i=1}^{\bfL} \bfl_i( \bfl_{ i - 1 } + 1 ) $
and
let
$ \bfa_i \colon \R^{\bfl_i} \to \R^{\bfl_i} $,
$ i \in \{ 1, 2, \ldots, \bfL \} $,
be functions.
Then we denote by 
$ \cN^{ \theta, s, \bfl_0 }_{ \bfa_1, \bfa_2, \ldots, \bfa_\bfL } \colon
\R^{\bfl_0} \to \R^{\bfl_\bfL} $
the function which satisfies
for all
$ x \in \R^{\bfl_0} $
that
\begin{equation}
\begin{split}
\bigl(
    \cN^{ \theta, s, \bfl_0 }_{ \bfa_1, \bfa_2, \ldots, \bfa_\bfL }
\bigr)(x)
=
\bigl(
    \bfa_\bfL
    \circ \cA_{ \bfl_\bfL, \bfl_{ \bfL - 1 } }^{ \theta, s + \sum_{i = 1}^{\bfL-1} \bfl_i ( \bfl_{ i - 1 } + 1) }
    \circ \bfa_{ \bfL - 1 }
    \circ \cA_{ \bfl_{ \bfL - 1 }, \bfl_{ \bfL - 2 } }^{ \theta, s + \sum_{i = 1}^{\bfL-2} \bfl_i ( \bfl_{ i - 1 } + 1) }
    \circ 
    \ldots
\qquad\qquad\\
\pushright{
    \ldots
    \circ \bfa_2
    \circ  \cA_{ \bfl_2, \bfl_1 }^{ \theta, s + \bfl_1 ( \bfl_0 + 1) }
    \circ \bfa_1
    \circ \cA_{ \bfl_1, \bfl_0 }^{ \theta, s }
\bigr)(x)
}
\end{split}
\end{equation}
(cf.~\cref{def:affine_function}).
\end{definition}

\subsection{Activation functions}

\begin{definition}[Multidimensional version]
\label{def:multidimensional_version}
Let
$ d \in \N $
and
let
$ \bfa \colon \R \to \R $
be a function.
Then
we denote by
$ \fM_{ \bfa, d } \colon \R^d \to \R^d $
the function which satisfies for all 
$ x = ( x_1, x_2, \ldots, x_d ) \in \R^d $
that
\begin{equation}
\fM_{ \bfa, d }( x )
=
( \bfa( x_1 ),
\bfa( x_2 ),
\ldots,
\bfa( x_d ) )
.
\end{equation}
\end{definition}

\begin{definition}[Rectifier function]
\label{def:ReLu}
We denote by
$ \fr \colon \R \to \R $
the function which satisfies for all
$ x \in \R $
that 
\begin{equation}
\fr( x )
=
\max\{ x, 0 \}
.
\end{equation}
\end{definition}

\begin{definition}[Multidimensional rectifier function]
\label{def:ReLu_multidim}
Let
$ d \in \N $.
Then we denote by
$ \fR_d \colon \R^d \to \R^d $
the function given by
\begin{equation}
\fR_d
=
\fM_{ \fr, d }
\end{equation}
(cf.~\cref{def:multidimensional_version,def:ReLu}).
\end{definition}

\begin{definition}[Clipping function]
\label{def:clipping_function}
Let
$ u \in [ -\infty, \infty ) $,
$ v \in ( u, \infty ] $.
Then
we denote by
$ \fc_{ u, v } \colon \R \to \R $
the function which satisfies for all
$ x \in \R $
that 
\begin{equation}
\fc_{ u, v }( x )
=
\max\{ u, \min\{ x, v \} \}
.
\end{equation}
\end{definition}

\begin{definition}[Multidimensional clipping function]
\label{def:clipping_multidim}
Let
$ d \in \N $,
$ u \in [ -\infty, \infty ) $,
$ v \in ( u, \infty ] $.
Then
we denote by
$ \fC_{ u, v, d } \colon \R^d \to \R^d $
the function given by
\begin{equation}
\fC_{ u, v, d }
=
\fM_{ \fc_{ u, v }, d }
\end{equation}
(cf.~\cref{def:multidimensional_version,def:clipping_function}).
\end{definition}

\subsection{Rectified DNNs}

\begin{definition}[Rectified clipped DNN]
\label{def:clipped_NN}
Let 
$ \bfd, \bfL \in \N $,
$ u \in [ -\infty, \infty ) $,
$ v \in ( u, \infty ] $,
$ \bfl = ( \bfl_0, \bfl_1, \ldots, \bfl_\bfL ) \in \N^{ \bfL + 1 } $,
$ \theta \in \R^\bfd $
satisfy
$ \bfd \geq \sum_{i=1}^{\bfL} \bfl_i( \bfl_{ i - 1 } + 1 ) $.
Then
we denote by
$ \clippedNN{\theta}{\bfl}{u}{v} \colon
\R^{ \bfl_0 } \to \R^{ \bfl_\bfL } $
the function which satisfies for all
$ x \in \R^{ \bfl_0 } $
that
\begin{equation}
\clippedNN{\theta}{\bfl}{u}{v}( x )
=
\begin{cases}
\bigl(
    \cN^{ \theta, 0, \bfl_0 }_{ \fC_{ u, v, \bfl_\bfL } }
\bigr)( x )
&
\colon
\bfL = 1
\\
\bigl(
    \cN^{ \theta, 0, \bfl_0 }_{ \fR_{ \bfl_1 }, \fR_{ \bfl_2 }, \ldots, \fR_{ \bfl_{ \bfL - 1 } }, \fC_{ u, v, \bfl_\bfL } }
\bigr)( x )
&
\colon
\bfL > 1
\end{cases}
\end{equation}
(cf.~\cref{def:NN,def:ReLu_multidim,def:clipping_multidim}).
\end{definition}

\section{Analysis of the approximation error}
\label{sec:approximation_error}

This section is devoted to establishing a convergence result for the approximation
of Lipschitz continuous functions by DNNs
(cf. \cref{prop:approximation_error}).
More precisely,
\cref{prop:approximation_error}
establishes that
a Lipschitz continuous function
defined on a $ d $-dimensional hypercube
for $ d \in \N $
can be approximated by DNNs
with convergence rate $ \nicefrac{1}{d} $
with respect to
a parameter $ A \in ( 0, \infty ) $
that bounds the architecture size (that is, depth and width) of the approximating DNN from below.
Key ingredients of the proof of
\cref{prop:approximation_error}
are
Beck, Jentzen, \& Kuckuck~\cite[Corollary~3.8]{BeckJentzenKuckuck2019arXiv}
as well as
the elementary 
covering number estimate
in \cref{lem:covering_number_cube}.
In order to improve the accessibility of \cref{lem:covering_number_cube},
we recall
the definition of covering numbers associated to a metric space
in \cref{def:covering_number},
which is \cite[Definition~3.11]{BeckJentzenKuckuck2019arXiv}.
\cref{lem:covering_number_cube}
provides upper bounds for the covering numbers of hypercubes
equipped with the metric induced by the $ p $-norm
(cf.\ \cref{def:p-norm})
for $ p \in [ 1, \infty ] $
and
is a generalisation of
Berner, Grohs, \& Jentzen~\cite[Lemma~2.7]{BernerGrohsJentzen2018arXiv}
(cf.\ Cucker \& Smale~\cite[Proposition~5]{CuckerSmale2002}
and \cite[Proposition~3.12]{BeckJentzenKuckuck2019arXiv}).
Furthermore,
we present in \cref{lem:approximation_error}
an elementary upper bound
for the error
arising when Lipschitz continuous functions
defined on a hypercube
are approximated by certain DNNs.
Additional DNN approximation results
can be found, e.g., in~\cite{SirignanoSpiliopoulos2018,
HanLong2018arXiv,
CheriditoJentzenRossmannek2019arXiv,
GrohsJentzenSalimova2019arXiv,
GuliyevIsmailov2018a,
GuliyevIsmailov2018b,
Bach2017,
Barron1993,
Barron1994,
BlumLi1991,
BolcskeiGrohsKutyniokPetersen2019,
BurgerNeubauer2001,
Candes1998,
ChenChen1995,
ChuiMhaskar1994,
Cybenko1989,
DeVoreOskolkovPetrushev1997,
EWang2018,
ElbraechterGrohsJentzenSchwab2018arXiv,
EldanShamir2016,
Ellacott1994,
Funahashi1989,
GribonvalKutyniokNielsenVoigtlaender2019,
GrohsHornungJentzenVonWurstemberger2018arXiv,
GrohsHornungJentzenZimmermann2019arXiv,
GrohsPerekrestenkoElbrachterBolcskei2019arXiv,
GuhringKutyniokPetersen2019arXiv,
HartmanKeelerKowalski1990,
Hornik1991,
Hornik1993,
HornikStinchcombeWhite1989,
HornikStinchcombeWhite1990,
HutzenthalerJentzenKruseNguyen2019arXiv,
JentzenSalimovaWelti2018arXiv,
KutyniokPetersenRaslanSchneider2019arXiv,
LeshnoLinPinkusSchocken1993,
Mhaskar1996,
MhaskarMicchelli1995,
MhaskarPoggio2016,
MontanelliDu2019,
Nguyen_ThienTran_Cong1999,
ParkSandberg1991,
PerekrestenkoPhilippElbrachterBolcskei2018arXiv,
PetersenRaslanVoigtlaender2018arXiv,
PetersenVoigtlaender2018arXiv,
PetersenVoigtlaender2018,
Pinkus1999,
ReisingerZhang2019arXiv,
Schmitt2000,
SchwabZech2019,	
ShahamCloningerCoifman2018,
ShenYangZhang2019arXiv,
ShenYangZhang2019,
VoigtlaenderPetersen2019arXiv,
Yarotsky2017,
Yarotsky2018arXiv}
and the references therein.

\subsection{A covering number estimate}

\begin{definition}[$p$-norm]
\label{def:p-norm}
We denote by
$ \lVert \cdot \rVert_p \colon \bigl( \bigcup_{d=1}^\infty \R^d \bigr) \to [ 0, \infty ) $, $ p \in [ 1, \infty ] $,
the functions which satisfy for all
$ p \in [ 1, \infty ) $,
$ d \in \N $,
$ \theta = ( \theta_1, \theta_2, \ldots, \theta_d ) \in \R^d $
that
\begin{equation}
\lVert \theta \rVert_p
=
\biggl(
    \smallsum_{i=1}^d
        \lvert \theta_i \rvert^p
\biggr)^{ \!\! \nicefrac{1}{p} }
\qquad
\text{and}
\qquad
\lVert \theta \rVert_\infty
=
\max_{ i \in \{ 1, 2, \ldots, d \} }
    \lvert \theta_i \rvert
.
\end{equation}
\end{definition}

\begin{definition}[Covering number]
\label{def:covering_number}
Let
$ ( E, \delta ) $
be a metric space
and
let
$ r \in [ 0, \infty ] $.
Then we denote by
$ \cC_{ ( E, \delta ), r } \in \N_0 \cup \{ \infty \} $ 
(we denote by
$ \cC_{ E, r } \in \N_0 \cup \{ \infty \} $)
the extended real number given by
\begin{equation}
\begin{split}
\cC_{ ( E, \delta ), r }
=
\inf\biggl(
    \biggl\{
        n \in \N_0 \colon
        \biggl[
        \exists \, A \subseteq E \colon
        \biggl(
        \arraycolsep=0pt \begin{array}{c}
            ( \lvert A \rvert \leq n )
            \land
            ( \forall \, x \in E \colon \\
            \exists \, a \in A \colon \delta( a, x ) \leq r )
        \end{array}
        \biggr)
        \biggr]
    \biggr\}
    \cup \{ \infty \}
\biggr)
.
\end{split}
\end{equation}
\end{definition}

\begin{lemma}
\label{lem:covering_number_cube}
Let
$ d \in \N $,
$ a \in \R $,
$ b \in ( a, \infty ) $,
$ r \in ( 0, \infty ) $,
for every
$ p \in [ 1, \infty ] $
let
$ \delta_p \colon ( [ a, b ]^d ) \times ( [ a, b ]^d ) \to [ 0, \infty ) $
satisfy for all
$ x, y \in [ a, b ]^d $
that
$ \delta_p( x, y ) = \lVert x - y \rVert_p $,
and
let
$ \lceil \cdot \rceil \colon [ 0, \infty ) \to \N_0 $
satisfy for all
$ x \in [ 0, \infty ) $
that
$ \lceil x \rceil = \min( [ x, \infty ) \cap \N_0 ) $
(cf.~\cref{def:p-norm}).
Then
\begin{enumerate}[(i)]
\item
\label{item:lem:covering_number_cube:1}
it holds for all
$ p \in [ 1, \infty ) $
that
\begin{equation}
\cC_{ ( [ a, b ]^d, \delta_p ), r }
\leq
\Bigl(
\Bigl\lceil
    \tfrac{ d^{ \nicefrac{1}{p} } ( b - a ) }{2r}
\Bigr\rceil
\Bigr)^{ \! d }
\leq
\begin{cases}
1
& \colon
r \geq \nicefrac{ d ( b - a ) }{2}
\\
\bigl(
    \tfrac{ d ( b - a ) }{r}
\bigr)^d
& \colon
r < \nicefrac{ d ( b - a ) }{2}
\end{cases}
\end{equation}
and
\item
\label{item:lem:covering_number_cube:2}
it holds that
\begin{equation}
\cC_{ ( [ a, b ]^d, \delta_\infty ), r }
\leq
\bigl(
\bigl\lceil
    \tfrac{ b - a }{2r}
\bigr\rceil
\bigr)^d
\leq
\begin{cases}
1
& \colon
r \geq \nicefrac{ ( b - a ) }{2}
\\
\bigl(
    \tfrac{ b - a }{r}
\bigr)^d
& \colon
r < \nicefrac{ ( b - a ) }{2}
\end{cases}
\end{equation}
\end{enumerate}
(cf.~\cref{def:covering_number}).
\end{lemma}
\begin{cproof2}{lem:covering_number_cube}
Throughout this proof
let
$ ( \fN_p )_{ p \in [ 1, \infty ] } \subseteq \N $
satisfy for all
$ p \in [ 1, \infty ) $
that
\begin{equation}
\label{eq:def_fNp}
\fN_p =
\Bigl\lceil
    \tfrac{ d^{ \nicefrac{1}{p} } ( b - a ) }{2r}
\Bigr\rceil
\qquad
\text{and}
\qquad
\fN_\infty =
\bigl\lceil
    \tfrac{ b - a }{2r}
\bigr\rceil
,
\end{equation}
for every
$ N \in \N $,
$ i \in \{ 1, 2, \ldots, N \} $
let
$ g_{ N, i } \in [ a, b ] $
be given by
$ g_{ N, i } = a + \nicefrac{ ( i - \nicefrac{1}{2} ) ( b - a ) }{N} $,
and
for every
$ p \in [ 1, \infty ] $
let
$ A_p \subseteq [ a, b ]^d $
be given by
$ A_p =
\{ g_{ \fN_p, 1 }, g_{ \fN_p, 2 }, \ldots, g_{ \fN_p, \fN_p } \}^d $.
Observe that it holds for all
$ N \in \N $,
$ i \in \{ 1, 2, \ldots, N \} $,
$ x \in [ a + \nicefrac{ ( i - 1 )( b - a ) }{N}, g_{ N, i } ] $
that
\begin{equation}
\label{eq:grid_estimate1}
\lvert x - g_{ N, i } \rvert
=
a + \tfrac{ ( i - \nicefrac{1}{2} ) ( b - a ) }{N}
-
x
\leq
a + \tfrac{ ( i - \nicefrac{1}{2} ) ( b - a ) }{N}
-
\bigl( a + \tfrac{ ( i - 1 )( b - a ) }{N} \bigr)
=
\tfrac{ b - a }{2N}
.
\end{equation}
In addition,
note that it holds for all
$ N \in \N $,
$ i \in \{ 1, 2, \ldots, N \} $,
$ x \in [ g_{ N, i }, a + \nicefrac{ i ( b - a ) }{N} ] $
that%
\begin{equation}
\label{eq:grid_estimate2}
\lvert x - g_{ N, i } \rvert
=
x
-
\bigl( a + \tfrac{ ( i - \nicefrac{1}{2} ) ( b - a ) }{N} \bigr)
\leq
a + \tfrac{ i ( b - a ) }{N}
-
\bigl( a + \tfrac{ ( i - \nicefrac{1}{2} ) ( b - a ) }{N} \bigr)
=
\tfrac{ b - a }{2N}
.
\end{equation}
Combining~\cref{eq:grid_estimate1} and~\cref{eq:grid_estimate2}
implies for all
$ N \in \N $,
$ i \in \{ 1, 2, \ldots, N \} $,
$ x \in [ a + \nicefrac{ ( i - 1 )( b - a ) }{N}, a + \nicefrac{ i ( b - a ) }{N} ] $
that
$ \lvert x - g_{ N, i } \rvert
\leq
\nicefrac{ ( b - a ) }{(2N)} $.
This proves that
for every
$ N \in \N $,
$ x \in [ a, b ] $
there exists
$ y \in \{ g_{ N, 1 }, g_{ N, 2 }, \ldots, g_{ N, N } \} $
such that
\begin{equation}
\label{eq:grid_estimate_final}
\lvert x - y \rvert
\leq
\tfrac{ b - a }{2N}
.
\end{equation}
This,
in turn,
establishes that
for every
$ p \in [ 1, \infty ) $,
$ x = ( x_1, x_2, \ldots, x_d )
\in [ a, b ]^d $
there exists
$ y = ( y_1, y_2, \ldots, y_d )
\in A_p $
such that
\begin{equation}
\label{eq:covering_Lp}
\delta_p( x, y )
=
\lVert x - y \rVert_p
=
\biggl(
    \smallsum_{i=1}^d
        \lvert x_i - y_i \rvert^p
\biggr)^{ \!\! \nicefrac{1}{p} }
\leq
\biggl(
    \smallsum_{i=1}^d
        \tfrac{ ( b - a )^p }{ ( 2\fN_p )^p }
\biggr)^{ \!\! \nicefrac{1}{p} }
=
\tfrac{ d^{ \nicefrac{1}{p} } ( b - a ) }{2\fN_p}
\leq
\tfrac{ d^{ \nicefrac{1}{p} } ( b - a ) 2r }{ 2 d^{ \nicefrac{1}{p} } ( b - a ) }
=
r
.
\end{equation}
Furthermore,
again
\cref{eq:grid_estimate_final}
shows that
for every
$ x = ( x_1, x_2, \ldots, x_d )
\in [ a, b ]^d $
there exists
$ y = ( y_1, y_2, \ldots, y_d )
\in A_\infty $
such that
\begin{equation}
\label{eq:covering_L_infty}
\delta_\infty( x, y )
=
\lVert x - y \rVert_\infty
=
\max_{ i \in \{ 1, 2, \ldots, d \} }
    \lvert x_i - y_i \rvert
\leq
\tfrac{ b - a }{2\fN_\infty}
\leq
\tfrac{ ( b - a ) 2r }{ 2 ( b - a ) }
=
r
.
\end{equation}
Note that
\cref{eq:covering_Lp},
\cref{eq:def_fNp},
and
the fact that
$ \forall \, x \in [ 0, \infty ) \colon
\lceil x \rceil
\leq
\mathbbm{1}_{ ( 0, 1 ] }( x )
+
2 x \mathbbm{1}_{ ( 1, \infty ) }( x )
=
\mathbbm{1}_{ ( 0, r ] }( r x )
+
2 x \mathbbm{1}_{ ( r, \infty ) }( r x ) $
yield for all
$ p \in [ 1, \infty ) $
that
\begin{equation}
\begin{split}
\cC_{ ( [ a, b ]^d, \delta_p ), r }
& \leq
\lvert A_p \rvert = ( \fN_p )^d
=
\Bigl(
\Bigl\lceil
    \tfrac{ d^{ \nicefrac{1}{p} } ( b - a ) }{2r}
\Bigr\rceil
\Bigr)^{ \! d }
\leq
\bigl(
\bigl\lceil
    \tfrac{ d ( b - a ) }{2r}
\bigr\rceil
\bigr)^d
\\ &
\leq
\bigl(
\mathbbm{1}_{ ( 0, r ] }\bigl( \tfrac{ d ( b - a ) }{2} \bigr)
+
\tfrac{ 2 d ( b - a ) }{2r}
\mathbbm{1}_{ ( r, \infty ) }\bigl( \tfrac{ d ( b - a ) }{2} \bigr)
\bigr)^d
\\ &
=
\mathbbm{1}_{ ( 0, r ] }\bigl( \tfrac{ d ( b - a ) }{2} \bigr)
+
\bigl(
\tfrac{ d ( b - a ) }{r}
\bigr)^d
\mathbbm{1}_{ ( r, \infty ) }\bigl( \tfrac{ d ( b - a ) }{2} \bigr)
.
\end{split}
\end{equation}
This proves~\cref{item:lem:covering_number_cube:1}.
In addition,
\cref{eq:covering_L_infty},
\cref{eq:def_fNp},
and
again
the fact that
$ \forall \, x \in [ 0, \infty ) \colon
\lceil x \rceil
\leq
\mathbbm{1}_{ ( 0, r ] }( r x )
+
2 x \mathbbm{1}_{ ( r, \infty ) }( r x ) $
demonstrate that
\begin{equation}
\cC_{ ( [ a, b ]^d, \delta_\infty ), r }
\leq
\lvert A_\infty \rvert = ( \fN_\infty )^d
=
\bigl(
\bigl\lceil
    \tfrac{ b - a }{2r}
\bigr\rceil
\bigr)^d
\leq
\mathbbm{1}_{ ( 0, r ] }\bigl( \tfrac{ b - a }{2} \bigr)
+
\bigl(
\tfrac{ b - a }{r}
\bigr)^d
\mathbbm{1}_{ ( r, \infty ) }\bigl( \tfrac{ b - a }{2} \bigr)
.
\end{equation}
This implies~\cref{item:lem:covering_number_cube:2}
and thus completes %
\end{cproof2}

\subsection{Convergence rates for the approximation error}

\begin{lemma}
\label{lem:approximation_error}
Let
$ d, \bfd, \bfL \in \N $,
$ L, a \in \R $,
$ b \in ( a, \infty ) $,
$ u \in [ -\infty, \infty ) $,
$ v \in ( u, \infty ] $,
$ \bfl = ( \bfl_0, \bfl_1, \ldots, \bfl_\bfL ) \in \N^{ \bfL + 1 } $,
assume
$ \bfl_0 = d $,
$ \bfl_\bfL = 1 $,
and
$ \bfd \geq \sum_{i=1}^{\bfL} \bfl_i( \bfl_{ i - 1 } + 1 ) $,
and
let
$ f \colon [ a, b ]^d \to ( [ u, v ] \cap \R ) $
satisfy for all
$ x, y \in [ a, b ]^d $
that
$ \lvert f( x ) - f( y ) \rvert \leq L \lVert x - y \rVert_{ 1 } $
(cf.~\cref{def:p-norm}).
Then
there exists $ \vartheta \in \R^\bfd $
such that
$ \lVert \vartheta \rVert_\infty
\leq \sup_{ x \in [ a, b ]^d } \lvert f( x ) \rvert $
and
\begin{equation}
\sup\nolimits_{ x \in [ a, b ]^d }
    \lvert \clippedNN{\vartheta}{\bfl}{u}{v}( x ) - f( x ) \rvert
\leq
\frac{ d L ( b - a ) }{ 2 }
\end{equation}
(cf.~\cref{def:clipped_NN}).
\end{lemma}
\begin{cproof2}{lem:approximation_error}
Throughout this proof
let
$ \fd \in \N $
be given by
$ \fd = \sum_{i=1}^{\bfL} \bfl_i( \bfl_{ i - 1 } + 1 ) $,
let
$ \bfm = ( \bfm_1, \bfm_2, \ldots, \bfm_d ) \in [ a, b ]^d $
satisfy for all
$ i \in \{ 1, 2, \ldots, d \} $
that
$ \bfm_i = \nicefrac{ ( a + b ) }{2} $,
and
let
$ \vartheta = ( \vartheta_1, \vartheta_2, \ldots, \vartheta_\bfd )\in \R^\bfd $
satisfy for all
$ i \in \{ 1, 2, \ldots, \bfd \} \setminus \{ \fd \} $
that
$ \vartheta_i = 0 $
and
$ \vartheta_\fd = f( \bfm ) $.
Observe that
the assumption that
$ \bfl_\bfL = 1 $
and
the fact that
$ \forall \, i \in \{ 1, 2, \ldots, \fd - 1 \} \colon
\vartheta_i = 0 $
show for all
$ x = ( x_1, x_2, \ldots, x_{ \bfl_{ \bfL - 1 } } ) \in \R^{ \bfl_{ \bfL - 1 } } $
that
\begin{equation}
\begin{split}
\cA_{ 1, \bfl_{ \bfL - 1 } }^{ \vartheta, \sum_{i = 1}^{\bfL-1} \bfl_i ( \bfl_{ i - 1 } + 1) }( x )
&
=
\biggl[
    \smallsum_{ i = 1 }^{ \bfl_{ \bfL - 1 } } \vartheta_{ \left[ \sum_{i = 1}^{\bfL-1} \bfl_i ( \bfl_{ i - 1 } + 1) \right] + i } x_i
\biggr]
+ \vartheta_{ \left[ \sum_{i = 1}^{\bfL-1} \bfl_i ( \bfl_{ i - 1 } + 1) \right] + \bfl_{ \bfL - 1 } + 1 }
\\ &
=
\biggl[
    \smallsum_{ i = 1 }^{ \bfl_{ \bfL - 1 } } \vartheta_{ \left[ \sum_{i = 1}^{\bfL} \bfl_i ( \bfl_{ i - 1 } + 1) \right] - ( \bfl_{ \bfL - 1 } - i + 1 ) } x_i
\biggr]
+ \vartheta_{ \sum_{i = 1}^{\bfL} \bfl_i ( \bfl_{ i - 1 } + 1) }
\\ &
=
\biggl[
    \smallsum_{ i = 1 }^{ \bfl_{ \bfL - 1 } } \vartheta_{ \fd - ( \bfl_{ \bfL - 1 } - i + 1 ) } x_i
\biggr]
+ \vartheta_{ \fd }
=
\vartheta_{ \fd }
=
f( \bfm )
\end{split}
\end{equation}
(cf.~\cref{def:affine_function}).
Combining this with
the fact that
$ f( \bfm ) \in [ u, v ] $
ensures for all
$ x \in \R^{ \bfl_{ \bfL - 1 } } $
that
\begin{equation}
\begin{split}
\bigl(
\fC_{ u, v, \bfl_\bfL }
\circ
\cA_{ \bfl_\bfL, \bfl_{ \bfL - 1 } }^{ \vartheta, \sum_{i = 1}^{\bfL-1} \bfl_i ( \bfl_{ i - 1 } + 1) }
\bigr)( x )
& =
\bigl(
\fC_{ u, v, 1 }
\circ
\cA_{ 1, \bfl_{ \bfL - 1 } }^{ \vartheta, \sum_{i = 1}^{\bfL-1} \bfl_i ( \bfl_{ i - 1 } + 1) }
\bigr)( x )
=
\fc_{ u, v }( f( \bfm ) )
\\ &
=
\max\{ u, \min\{ f( \bfm ), v \} \}
=
\max\{ u, f( \bfm ) \}
=
f( \bfm )
\end{split}
\end{equation}
(cf.~\cref{def:clipping_function,def:clipping_multidim}).
This proves for all
$ x \in \R^d $
that
\begin{equation}
\label{eq:NN_constant}
\clippedNN{\vartheta}{\bfl}{u}{v}( x )
=
f( \bfm )
.
\end{equation}
In addition,
note that it holds for all
$ x \in [ a, \bfm_1 ] $,
$ \fx \in [ \bfm_1, b ] $
that
$ \lvert \bfm_1 - x \rvert
= \bfm_1 - x
= \nicefrac{ ( a + b ) }{2} - x
\leq \nicefrac{ ( a + b ) }{2} - a
= \nicefrac{ ( b - a ) }{2} $
and
$ \lvert \bfm_1 - \fx \rvert
= \fx - \bfm_1
= \fx - \nicefrac{ ( a + b ) }{2}
\leq b - \nicefrac{ ( a + b ) }{2}
= \nicefrac{ ( b - a ) }{2} $.
The assumption that
$ \forall \, x, y \in [ a, b ]^d \colon
\lvert f( x ) - f( y ) \rvert \leq L \lVert x - y \rVert_{ 1 } $
and
\cref{eq:NN_constant}
hence
demonstrate for all
$ x = ( x_1, x_2, \ldots, x_d ) \in [ a, b ]^d $
that
\begin{equation}
\begin{split}
\lvert \clippedNN{\vartheta}{\bfl}{u}{v}( x ) - f( x ) \rvert
& =
\lvert f( \bfm ) - f( x ) \rvert
\leq
L \lVert \bfm - x \rVert_{ 1 }
=
L
\smallsum_{i=1}^d
    \lvert \bfm_i - x_i \rvert
\\ &
=
L
\smallsum_{i=1}^d
    \lvert \bfm_1 - x_i \rvert
\leq
\smallsum_{i=1}^d
    \frac{ L ( b - a ) }{2}
=
\frac{ d L ( b - a ) }{2}
.
\end{split}
\end{equation}
This and the fact that
$ \lVert \vartheta \rVert_\infty
=
\max_{ i \in \{ 1, 2, \ldots, \bfd \} }
    \lvert \vartheta_i \rvert
= \lvert f( \bfm ) \rvert
\leq \sup_{ x \in [ a, b ]^d } \lvert f( x ) \rvert $
complete
\end{cproof2}

\begin{proposition}
\label{prop:approximation_error}
Let
$ d, \bfd, \bfL \in \N $,
$ A \in ( 0, \infty ) $,
$ L, a \in \R $,
$ b \in ( a, \infty ) $,
$ u \in [ -\infty, \infty ) $,
$ v \in ( u, \infty ] $,
$ \bfl = ( \bfl_0, \bfl_1, \ldots, \bfl_\bfL ) \in \N^{ \bfL + 1 } $,
assume
$ \bfL \geq \nicefrac{ A \mathbbm{1}_{ \smash{ ( 6^d, \infty ) } }( A ) }{ ( 2d ) } + 1 $,
$ \bfl_0 = d $,
$ \bfl_1 \geq A \mathbbm{1}_{ \smash{ ( 6^d, \infty ) } }( A ) $,
$ \bfl_\bfL = 1 $,
and
$ \bfd \geq \sum_{i=1}^{\bfL} \bfl_i( \bfl_{ i - 1 } + 1 ) $,
assume for all
$ i \in \{ 2, 3, \ldots \} \cap [ 0, \bfL ) $
that
$ \bfl_i \geq \mathbbm{1}_{ \smash{ ( 6^d, \infty ) } }( A ) \max\{ \nicefrac{A}{d} - 2i + 3, 2 \}$,
and
let
$ f \colon [ a, b ]^d \to ( [ u, v ] \cap \R ) $
satisfy for all
$ x, y \in [ a, b ]^d $
that
$ \lvert f( x ) - f( y ) \rvert \leq L \lVert x - y \rVert_{ 1 } $
(cf.~\cref{def:p-norm}).
Then
there exists $ \vartheta \in \R^\bfd $
such that
$ \lVert \vartheta \rVert_\infty
\leq \max\{ 1, L, \lvert a \rvert, \lvert b \rvert, 2[ \sup_{ x \in [ a, b ]^d } \lvert f( x ) \rvert ] \} $
and
\begin{equation}
\sup\nolimits_{ x \in [ a, b ]^d }
    \lvert \clippedNN{\vartheta}{\bfl}{u}{v}( x ) - f( x ) \rvert
\leq
\frac{ 3 d L ( b - a ) }{ A^{ \nicefrac{1}{d} } }
\end{equation}
(cf.~\cref{def:clipped_NN}).
\end{proposition}
\begin{cproof2}{prop:approximation_error}
Throughout this proof
assume w.l.o.g.\ that
$ A > 6^d $
(cf.~Lem\-ma~\ref{lem:approximation_error}),
let $ \fN \in \N $ be given by
\begin{equation}
\label{eq:def_fN}
\fN
=
\max\Bigl\{
    \fn \in \N \colon
    \fn \leq \bigl( \tfrac{ A }{ 2d } \bigr)^{ \nicefrac{1}{d} }
\Bigr\}
,
\end{equation}
let
$ r \in ( 0, \infty ) $
be given by
$ r = \nicefrac{ d ( b - a ) }{ (2\fN) } $,
let
$ \delta \colon ( [ a, b ]^d ) \times ( [ a, b ]^d ) \to [ 0, \infty ) $
satisfy for all
$ x, y \in [ a, b ]^d $
that
$ \delta( x, y ) = \lVert x - y \rVert_1 $,
let
$ \scrD \subseteq [ a, b ]^d $
satisfy
$ \lvert \scrD \rvert
= \max\{ 2, \cC_{ \smash{ ( [ a, b ]^d, \delta ), r } } \} $
and%
\begin{equation}
\label{eq:covering_property}
\sup\nolimits_{ x \in [ a, b ]^d }
\inf\nolimits_{ y \in \scrD }
    \delta( x, y )
\leq
r
\end{equation}
(cf.~\cref{def:covering_number}),
and
let
$ \lceil \cdot \rceil \colon [ 0, \infty ) \to \N_0 $
satisfy for all
$ x \in [ 0, \infty ) $
that
$ \lceil x \rceil = \min( [ x, \infty ) \cap \N_0 ) $.
Note that it holds for all
$ \fd \in \N $
that
\begin{equation}
\label{eq:exponential_estimate}
2 \fd
\leq 2 \cdot 2^{ \fd - 1 }
= 2^\fd
.
\end{equation}
This implies that
$ 3^d
= \nicefrac{6^d}{2^d}
\leq \nicefrac{A}{(2d)} $.
Equation~\cref{eq:def_fN}
hence
demonstrates that
\begin{equation}
\label{eq:fN_lower_bound}
2
\leq
\tfrac{2}{3}
\bigl( \tfrac{ A }{ 2d } \bigr)^{ \nicefrac{1}{d} }
=
\bigl( \tfrac{ A }{ 2d } \bigr)^{ \nicefrac{1}{d} }
-
\tfrac{1}{3}
\bigl( \tfrac{ A }{ 2d } \bigr)^{ \nicefrac{1}{d} }
\leq
\bigl( \tfrac{ A }{ 2d } \bigr)^{ \nicefrac{1}{d} }
-
1
<
\fN
.
\end{equation}
This
and~\cref{item:lem:covering_number_cube:1}
in
\cref{lem:covering_number_cube}
(with
$ \delta_1 \leftarrow \delta $,
$ p \leftarrow 1 $
in the notation of~\cref{item:lem:covering_number_cube:1}
in \cref{lem:covering_number_cube})
establish that
\begin{equation}
\begin{split}
\lvert \scrD \rvert
=
\max\{ 2, \cC_{ ( [ a, b ]^d, \delta ), r } \}
\leq
\max\Bigl\{
    2,
    \Bigl(
    \Bigl\lceil
        \tfrac{ d ( b - a ) }{2r}
    \Bigr\rceil
    \Bigr)^{ \! d }
\Bigr\}
=
\max\{ 2, ( \lceil \fN \rceil )^d \}
=
\fN^d
.
\end{split}
\end{equation}
Combining this with~\cref{eq:def_fN}
proves that
\begin{equation}
\label{eq:fN_upper_bound}
4
\leq
2 d \lvert \scrD \rvert
\leq
2 d \fN^d
\leq
\tfrac{ 2 d A }{2d}
=
A
.
\end{equation}
The fact that
$ \bfL
\geq \nicefrac{ A \mathbbm{1}_{ \smash{ ( 6^d, \infty ) } }( A ) }{ ( 2d ) } + 1
= \nicefrac{A}{(2d)} + 1 $
hence yields that
$ \lvert \scrD \rvert
\leq
\nicefrac{A}{(2d)}
\leq
\bfL - 1 $.
This,
\cref{eq:fN_upper_bound},
and
the facts that
$ \bfl_1 \geq A \mathbbm{1}_{ \smash{ ( 6^d, \infty ) } }( A ) = A $
and
$ \forall \, i \in \{ 2, 3, \ldots \} \cap [ 0, \bfL ) = \{ 2, 3, \ldots, \bfL - 1 \} \colon
\bfl_i
\geq \mathbbm{1}_{ \smash{ ( 6^d, \infty ) } }( A ) \max\{ \nicefrac{A}{d} - 2i + 3, 2 \}
= \max\{ \nicefrac{A}{d} - 2i + 3, 2 \} $
imply
for all
$ i \in \{ 2, 3, \ldots, \lvert \scrD \rvert \} $
that
\begin{equation}
\label{eq:bfl_assumptions1}
\bfL
\geq \lvert \scrD \rvert + 1,
\quad
\bfl_1
\geq A
\geq 2 d \lvert \scrD \rvert,
\quad\text{and}\quad
\bfl_i
\geq \nicefrac{A}{d} - 2i + 3
\geq 2 \lvert \scrD \rvert - 2i + 3
.
\end{equation}
In addition,
the fact that
$ \forall \, i \in \{ 2, 3, \ldots \} \cap [ 0, \bfL ) \colon
\bfl_i \geq \max\{ \nicefrac{A}{d} - 2i + 3, 2 \} $
ensures for all
$ i \in \N \cap ( \lvert \scrD \rvert, \bfL ) $
that
\begin{equation}
\label{eq:bfl_assumptions2}
\bfl_i \geq 2
.
\end{equation}
Furthermore,
observe that it holds
for all
$ x = ( x_1, x_2, \ldots, x_d ), y = ( y_1, y_2, \ldots, y_d ) \in [ a, b ]^d $
that
\begin{equation}
\lvert f( x ) - f( y ) \rvert
\leq
L \lVert x - y \rVert_{ 1 }
=
L
\biggl[
\smallsum_{i=1}^d
    \lvert x_i - y_i \rvert
\biggr]
.
\end{equation}
This,
the assumptions that
$ \bfl_0 = d $,
$ \bfl_\bfL = 1 $,
and
$ \bfd \geq \sum_{i=1}^{\bfL} \bfl_i( \bfl_{ i - 1 } + 1 ) $,
\cref{eq:bfl_assumptions1}--\cref{eq:bfl_assumptions2},
and
Beck, Jentzen, \& Kuckuck~\cite[Corollary~3.8]{BeckJentzenKuckuck2019arXiv}
(with
$ d \leftarrow d $,
$ \fd \leftarrow \bfd $,
$ \fL \leftarrow \bfL $,
$ L \leftarrow L $,
$ u \leftarrow u $,
$ v \leftarrow v $,
$ D \leftarrow [ a, b ]^d $,
$ f \leftarrow f $,
$ \cM \leftarrow \scrD $,
$ l \leftarrow \bfl $
in the notation of \cite[Corollary~3.8]{BeckJentzenKuckuck2019arXiv})
show that
there exists
$ \vartheta \in \R^\bfd $
such that
$ \lVert \vartheta \rVert_{ \infty }
\leq
\max\{ 1, L, \sup_{ x \in \scrD } \lVert x \rVert_{ \infty }, 2 [ \sup_{ x \in \scrD } \lvert f( x ) \rvert ] \} $
and%
\begin{equation}
\label{eq:NN_approximation}
\begin{split}
\sup_{ x \in [ a, b ]^d }
    \, \lvert \clippedNN{\vartheta}{\bfl}{u}{v}( x ) - f( x ) \rvert
& \leq
2 L
\biggl[
    \sup_{ x = ( x_1, x_2, \ldots, x_d ) \in [ a, b ]^d }
    \biggl(
        \inf_{ y = ( y_1, y_2, \ldots, y_d ) \in \scrD }
        \smallsum_{i=1}^d \lvert x_i - y_i \rvert
    \biggr)
\biggr]
\\ &
=
2 L
\biggl[
\sup_{ x \in [ a, b ]^d }
\inf_{ y \in \scrD }
    \lVert x - y \rVert_{ 1 }
\biggr]
=
2 L
\biggl[
\sup_{ x \in [ a, b ]^d }
\inf_{ y \in \scrD }
    \delta( x, y )
\biggr]
.
\end{split}
\end{equation}
Note that
this demonstrates that
\begin{equation}
\label{eq:norm_estimate}
\lVert \vartheta \rVert_{ \infty }
\leq
\max\{ 1, L, \lvert a \rvert, \lvert b \rvert, 2 [ \sup\nolimits_{ x \in [ a, b ]^d } \lvert f( x ) \rvert ] \}
.
\end{equation}
Moreover,
\cref{eq:NN_approximation}
and
\cref{eq:covering_property}--\cref{eq:fN_lower_bound}
prove that
\begin{equation}
\begin{split}
&
\sup\nolimits_{ x \in [ a, b ]^d }
    \lvert \clippedNN{\vartheta}{\bfl}{u}{v}( x ) - f( x ) \rvert
\leq
2 L
\bigl[
\sup\nolimits_{ x \in [ a, b ]^d }
\inf\nolimits_{ y \in \scrD }
    \delta( x, y )
\bigr]
\leq
2 L r
\\ &
=
\frac{ d L ( b - a ) }{ \fN }
\leq
\frac{ d L ( b - a ) }{ \frac{2}{3} \bigl( \frac{ A }{ 2d } \bigr)^{ \nicefrac{1}{d} } }
=
\frac{ ( 2 d )^{ \nicefrac{1}{d} } 3 d L ( b - a ) }{ 2 A^{ \nicefrac{1}{d} } }
\leq
\frac{ 3 d L ( b - a ) }{ A^{ \nicefrac{1}{d} } }
.
\end{split}
\end{equation}
Combining this with
\cref{eq:norm_estimate} completes
\end{cproof2}

\section{Analysis of the generalisation error}
\label{sec:generalisation_error}

In this section
we consider the
\emph{worst-case} generalisation error
arising in deep learning based empirical risk minimisation
with quadratic loss function
for DNNs with a fixed architecture
and weights and biases bounded in size by a fixed constant
(cf.\ \cref{cor:generalisation_error} in \cref{sec:strong_rates_generalisation_error}).
We prove
that this worst-case generalisation error
converges in the probabilistically strong sense
with rate $ \nicefrac{1}{2} $ (up to a logarithmic factor)
with respect to the number of samples used for calculating the empirical risk
and
that the constant in the corresponding upper bound for the worst-case generalisation error
scales favourably (i.e., only very moderately) in terms of
depth and width of the DNNs employed,
cf.~\cref{item:cor:generalisation_error:2} in \cref{cor:generalisation_error}.
\cref{cor:generalisation_error}
is a consequence of
the main result of this section, \cref{prop:generalisation_error} in \cref{sec:strong_rates_generalisation_error},
which provides a similar conclusion in a more general setting.
The proofs of
\cref{prop:generalisation_error} and \cref{cor:generalisation_error}, respectively,
rely on the tools
developed in the two preceding subsections,
\cref{sec:MC_estimates,sec:uniform_strong_error}.

On the one hand,
\cref{sec:MC_estimates}
provides an essentially well-known
estimate for the $ L^p $-error of Monte Carlo-type approximations,
cf.\ \cref{cor:MC_Lp}.
\cref{cor:MC_Lp} is a consequence of
the well-known result stated here as \cref{prop:MC_Lp},
which, in turn, follows directly from, e.g.,
Cox et al.~\cite[Corollary~5.11]{CoxHutzenthalerJentzenVanNeervenWelti2016arXiv}
(with
$ M \leftarrow M $,
$ q \leftarrow 2 $,
$ ( E, \lVert \cdot \rVert_E ) \leftarrow ( \R^d, \lVert \cdot \rVert_2 \vert_{ \R^d } ) $,
$ ( \Omega, \scrF, \P ) \leftarrow ( \Omega, \cF, \P ) $,
$ ( \xi_j )_{ j \in \{ 1, 2, \ldots, M \} } \leftarrow ( X_j )_{ j \in \{ 1, 2, \ldots, M \} } $,
$ p \leftarrow p $
in the notation of~\cite[Corollary~5.11]{CoxHutzenthalerJentzenVanNeervenWelti2016arXiv}
and \cref{prop:MC_Lp}, respectively).
In the proof of \cref{cor:MC_Lp}
we also apply \cref{lem:Kahane_Khintchine},
which is Grohs et al.~\cite[Lemma~2.2]{GrohsHornungJentzenVonWurstemberger2018arXiv}.
In order to make the statements of
\cref{lem:Kahane_Khintchine,prop:MC_Lp}
more accessible for the reader,
we recall in \cref{def:Rademacher_family}
(cf., e.g., \cite[Definition~5.1]{CoxHutzenthalerJentzenVanNeervenWelti2016arXiv})
the notion of a Rademacher family
and in \cref{def:Kahane_Khintchine}
(cf., e.g., \cite[Definition~5.4]{CoxHutzenthalerJentzenVanNeervenWelti2016arXiv}
or Gonon et al.~\cite[Definition~2.1]{GononGrohsJentzenKoflerSiska2019arXiv})
the notion of the $ p $-Kahane--Khintchine constant.

On the other hand,
we derive in \cref{sec:uniform_strong_error}
uniform $ L^p $-estimates for Lipschitz continuous random fields
with a separable metric space as index set
(cf.\ \cref{lem:Lp_sup_centred,lem:Lp_sup_MC,cor:Lp_sup_MC}).
These estimates are uniform in the sense that
the supremum over the index set
is \emph{inside} the expectation belonging to the $ L^p $-norm,
which is necessary since we intend to prove error bounds for the
\emph{worst-case} generalisation error,
as illustrated above.
One of the elementary but crucial arguments in our derivation of such uniform $ L^p $-estimates
is given in \cref{lem:Lp_sup_covering_number}
(cf.\ \cref{lem:Lp_sup_random_field}).
Roughly speaking,
\cref{lem:Lp_sup_covering_number}
illustrates how the
$ L^p $-norm of a supremum
can be bounded from above by
the supremum of certain $ L^p $-norms,
where the $ L^p $-norms are integrating over a general measure space
and where the suprema are taken over a general (bounded) separable metric space.
Furthermore,
the elementary and well-known
\cref{lem:measurability_sup,lem:measurability_sup_centred},
respectively,
follow immediately from
Beck, Jentzen, \& Kuckuck~\cite[(ii) in Lemma~3.13 and (ii) in Lemma~3.14]{BeckJentzenKuckuck2019arXiv}
and ensure that the mathematical statements of
\cref{lem:Lp_sup_random_field,lem:Lp_sup_covering_number,lem:Lp_sup_centred}
do indeed make sense.

The results in \cref{sec:uniform_strong_error,sec:strong_rates_generalisation_error}
are in parts inspired by~\cite[Subsection~3.2]{BeckJentzenKuckuck2019arXiv}
and we refer, e.g., to
\cite{BartlettBousquetMendelson2005,
CuckerSmale2002,
BernerGrohsJentzen2018arXiv,
GyorfiKohlerKrzyzakWalk2002,
Massart2007,
Shalev_ShwartzBen_David2014,
VanDeGeer2000,
EMaWu2019,
EMaWang2019arXiv,
EMaWu2020online}
and the references therein
for further results on the generalisation error.

\subsection{Monte Carlo estimates}
\label{sec:MC_estimates}

\begin{definition}[Rademacher family]
\label{def:Rademacher_family}
Let
$ ( \Omega, \cF, \P ) $
be a probability space
and
let $ J $ be a set.
Then we say that 
$ ( r_j )_{ j \in J } $
is a $ \P $-Rademacher family
if and only if
it holds that
$ r_j \colon \Omega \to \{ -1, 1 \} $,
$ j \in J $,
are independent random variables 
with
$ \forall \, j \in J \colon
\P( r_j = 1 )
= \P( r_j = - 1 ) $.
\end{definition}

\begin{definition}[$ p $-Kahane--Khintchine constant]
\label{def:Kahane_Khintchine}
Let $ p \in ( 0, \infty ) $.
Then
we denote by
$ \fK_p \in ( 0, \infty ] $
the extended real number given by
\begin{equation}
\fK_p =
\sup \mathopen{} \left\{
    c \in [ 0, \infty ) \colon
    \mathopen{} \left[
    \arraycolsep=0pt \begin{array}{c}
        \exists \, \R\text{-Banach space}\ ( E, \lVert \cdot \rVert_E ) \colon \\
        \exists \, \text{probability space}\ ( \Omega, \cF, \P ) \colon \\
        \exists \, \P\text{-Rademacher family}\ ( r_j )_{ j \in \N } \colon \\
        \exists \, k \in \N \colon \exists \, x_1, x_2, \ldots, x_k \in E \setminus \{ 0 \} \colon \\
        \Bigl(
            \E\Bigl[ \bigl\lVert
                \sum_{ j = 1 }^k r_j x_j
            \bigr\rVert_E^p \Bigr]
        \Bigr)^{ \! \nicefrac{1}{p} }
        =
        c \Bigl(
            \E\Bigl[ \bigl\lVert
                \sum_{ j = 1 }^k r_j x_j
            \bigr\rVert_E^2 \Bigr]
        \Bigr)^{ \! \nicefrac{1}{2} }
    \end{array}
    \right]\mathclose{}
\right\}\mathclose{}
\end{equation}
(cf.~\cref{def:Rademacher_family}).
\end{definition}

\begin{lemma}
\label{lem:Kahane_Khintchine}
It holds for all
$ p \in [ 2, \infty ) $
that
$ \fK_p \leq \sqrt{ p - 1 } < \infty $
(cf.~\cref{def:Kahane_Khintchine}).
\end{lemma}

\begin{proposition}
\label{prop:MC_Lp}
Let
$ d, M \in \N $,
$ p \in [ 2, \infty ) $,
let
$ ( \Omega, \cF, \P ) $
be a probability space,
\linebreak
let
$ X_j \colon \Omega \to \R^d $,
$ j \in \{ 1, 2, \ldots, M \} $,
be independent random variables,
and assume
$ \max_{ j \in \{ 1, 2, \ldots, M \} } \E[ \lVert X_j \rVert_2 ] < \infty $
(cf.\ \cref{def:p-norm}).
Then
\begin{equation}
\biggl(
\E\biggl[
    \biggl\lVert
    \biggl[
    \smallsum_{j=1}^M
        X_j
    \biggr]
    -
    \E\biggl[
    \smallsum_{j=1}^M
        X_j
    \biggr]
    \biggr\rVert_2^p
\biggr]
\biggr)^{ \!\! \nicefrac{1}{p} }
\leq
2 \fK_p
\biggl[
\smallsum_{j=1}^M
    \bigl(
    \E\bigl[
        \lVert X_j - \E[ X_j ] \rVert_2^p
    \bigr]
    \bigr)^{ \nicefrac{2}{p} }
\biggr]^{ \nicefrac{1}{2} }
\end{equation}
(cf.~\cref{def:Kahane_Khintchine,lem:Kahane_Khintchine}).
\end{proposition}

\begin{corollary}
\label{cor:MC_Lp}
Let
$ d, M \in \N $,
$ p \in [ 2, \infty ) $,
let
$ ( \Omega, \cF, \P ) $
be a probability space,
\linebreak
let
$ X_j \colon \Omega \to \R^d $,
$ j \in \{ 1, 2, \ldots, M \} $,
be independent random variables,
and assume
$ \max_{ j \in \{ 1, 2, \ldots, M \} } \E[ \lVert X_j \rVert_2 ] < \infty $
(cf.\ \cref{def:p-norm}).
Then
\begin{equation}
\biggl(
\E\biggl[
    \biggl\lVert
    \frac{1}{M}
    \biggl[
    \smallsum_{j=1}^M
        X_j
    \biggr]
    -
    \E\biggl[
    \frac{1}{M}
    \smallsum_{j=1}^M
        X_j
    \biggr]
    \biggr\rVert_2^p
\biggr]
\biggr)^{ \!\! \nicefrac{1}{p} }
\leq
\frac{ 2 \sqrt{ p - 1 } }{ \sqrt{M} }
\biggl[
\max_{ j \in \{ 1, 2, \ldots, M \} }
    \bigl(
    \E\bigl[
        \lVert X_j - \E[ X_j ] \rVert_2^p
    \bigr]
    \bigr)^{ \nicefrac{1}{p} }
\biggr]
.
\end{equation}
\end{corollary}
\begin{cproof}{cor:MC_Lp}
Observe that
\cref{prop:MC_Lp}
and
\cref{lem:Kahane_Khintchine}
imply that
\begin{equation}
\begin{split}
&
\biggl(
\E\biggl[
    \biggl\lVert
    \frac{1}{M}
    \biggl[
    \smallsum_{j=1}^M
        X_j
    \biggr]
    -
    \E\biggl[
    \frac{1}{M}
    \smallsum_{j=1}^M
        X_j
    \biggr]
    \biggr\rVert_2^p
\biggr]
\biggr)^{ \!\! \nicefrac{1}{p} }
\\ &
=
\frac{1}{M}
\biggl(
\E\biggl[
    \biggl\lVert
    \biggl[
    \smallsum_{j=1}^M
        X_j
    \biggr]
    -
    \E\biggl[
    \smallsum_{j=1}^M
        X_j
    \biggr]
    \biggr\rVert_2^p
\biggr]
\biggr)^{ \!\! \nicefrac{1}{p} }
\\ &
\leq
\frac{ 2 \fK_p }{M}
\biggl[
\smallsum_{j=1}^M
    \bigl(
    \E\bigl[
        \lVert X_j - \E[ X_j ] \rVert_2^p
    \bigr]
    \bigr)^{ \nicefrac{2}{p} }
\biggr]^{ \nicefrac{1}{2} }
\\ &
\leq
\frac{ 2 \fK_p }{M}
\biggl[
M
\biggl(
\max_{ j \in \{ 1, 2, \ldots, M \} }
    \bigl(
    \E\bigl[
        \lVert X_j - \E[ X_j ] \rVert_2^p
    \bigr]
    \bigr)^{ \nicefrac{2}{p} }
\biggr)
\biggr]^{ \nicefrac{1}{2} }
\\ &
=
\frac{ 2 \fK_p }{ \sqrt{M} }
\biggl[
\max_{ j \in \{ 1, 2, \ldots, M \} }
    \bigl(
    \E\bigl[
        \lVert X_j - \E[ X_j ] \rVert_2^p
    \bigr]
    \bigr)^{ \nicefrac{1}{p} }
\biggr]
\\ &
\leq
\frac{ 2 \sqrt{ p - 1 } }{ \sqrt{M} }
\biggl[
\max_{ j \in \{ 1, 2, \ldots, M \} }
    \bigl(
    \E\bigl[
        \lVert X_j - \E[ X_j ] \rVert_2^p
    \bigr]
    \bigr)^{ \nicefrac{1}{p} }
\biggr]
\end{split}
\end{equation}
(cf.~\cref{def:Kahane_Khintchine}).
\end{cproof}

\subsection{Uniform strong error estimates for random fields}
\label{sec:uniform_strong_error}

\begin{lemma}
\label{lem:measurability_sup}
Let
$ ( E, \scrE ) $
be a separable topological space,
assume
$ E \neq \emptyset $,
let $ ( \Omega, \cF ) $
be a measurable space,
let
$ f_x \colon \Omega \to \R $, $ x \in E $,
be $ \cF $/$ \cB( \R ) $-measurable functions,
and
assume for all 
$ \omega \in \Omega $ 
that 
$ E \ni x \mapsto f_x( \omega ) \in \R $
is a continuous function.
Then
it holds that the function
\begin{equation}
\Omega \ni \omega \mapsto
\sup\nolimits_{ x \in E } f_x( \omega )
\in \R \cup \{ \infty \}
\end{equation}
is $ \cF $/$ \cB( \R \cup \{ \infty \} ) $-measurable.
\end{lemma}

\begin{lemma}
\label{lem:measurability_sup_centred}
Let
$ ( E, \delta ) $
be a separable metric space,
assume
$ E \neq \emptyset $,
let
$ L \in \R $,
let
$ ( \Omega, \cF, \P ) $
be a probability space,
let
$ Z_x \colon \Omega \to \R $, $ x \in E $,
be random variables,
and assume for all
$ x, y \in E $
that
$ \E[ \lvert Z_x \rvert ] < \infty $
and
$ \lvert Z_x - Z_y \rvert \leq L \delta( x, y ) $.
Then
it holds that the function%
\begin{equation}
\Omega \ni \omega \mapsto
\sup\nolimits_{ x \in E } \lvert Z_x( \omega ) - \E[ Z_x ] \rvert
\in [ 0, \infty ]
\end{equation}
is $ \cF $/$ \cB( [ 0, \infty ] ) $-measurable.
\end{lemma}

\begin{lemma}
\label{lem:Lp_sup_random_field}
Let
$ ( E, \delta ) $
be a separable metric space,
let
$ N \in \N $,
$ p, L, r_1, r_2, \ldots, r_N \in [ 0, \infty ) $,
$ z_1, z_2, \ldots, z_N \in E $
satisfy
$ E \subseteq
\bigcup_{ i = 1 }^N
\{ x \in E \colon \delta( x, z_i ) \leq r_i \} $,
let
$ ( \Omega, \cF, \mu ) $
be a measure space,
let
$ Z_x \colon \Omega \to \R $, $ x \in E $,
be $ \cF $/$ \cB( \R ) $-measurable functions,
and assume for all
$ \omega \in \Omega $,
$ x, y \in E $
that
$ \lvert Z_x( \omega ) - Z_y( \omega ) \rvert \leq L \delta( x, y ) $.
Then
\begin{equation}
\int_\Omega
    \sup_{ x \in E }
    \, \lvert Z_x( \omega ) \rvert^p
\, \mu( \dd \omega )
\leq
\smallsum_{ i = 1 }^N
\int_\Omega
    ( L r_i + \lvert Z_{ z_i }( \omega ) \rvert )^p
\, \mu( \dd \omega )
\end{equation}
(cf.~\cref{lem:measurability_sup}).
\end{lemma}
\begin{cproof}{lem:Lp_sup_random_field}
Throughout this proof
let
$ B_1, B_2, \ldots, B_N \subseteq E $
satisfy for all
$ i \in \{ 1, 2, \ldots, N \} $
that
$ B_i =
\{ x \in E \colon \delta( x, z_i ) \leq r_i \} $.
Note that
the fact that
$ E = \bigcup_{ i = 1 }^N B_i $
shows for all
$ \omega \in \Omega $
that
\begin{equation}
\sup\nolimits_{ x \in E }
\lvert Z_x( \omega ) \rvert
=
\sup\nolimits_{ x \in \left( \bigcup_{ i = 1 }^N B_i \right) }
\lvert Z_x( \omega ) \rvert
=
\max\nolimits_{ i \in \{ 1, 2, \ldots, N \} }
\sup\nolimits_{ x \in B_i }
\lvert Z_x( \omega ) \rvert
.
\end{equation}
This
establishes that
\begin{equation}
\label{eq:sup_power}
\begin{split}
& \int_\Omega
    \sup_{ x \in E }
    \, \lvert Z_x( \omega ) \rvert^p
\, \mu( \dd \omega )
=
\int_\Omega
    \max_{ i \in \{ 1, 2, \ldots, N \} }
    \sup_{ x \in B_i }
    \lvert Z_x( \omega ) \rvert^p
\, \mu( \dd \omega )
\\ &
\leq
\int_\Omega
    \smallsum_{ i = 1 }^N
    \sup_{ x \in B_i }
    \lvert Z_x( \omega ) \rvert^p
\, \mu( \dd \omega )
=
\smallsum_{ i = 1 }^N
\int_\Omega
    \sup_{ x \in B_i }
    \lvert Z_x( \omega ) \rvert^p
\, \mu( \dd \omega )
.
\end{split}
\end{equation}
Furthermore,
the assumption that
$ \forall \, \omega \in \Omega,
\, x, y \in E \colon \lvert Z_x( \omega ) - Z_y( \omega ) \rvert \leq L \delta( x, y ) $
implies for all
$ \omega \in \Omega $,
$ i \in \{ 1, 2, \ldots, N \} $,
$ x \in B_i $
that
\begin{equation}
\begin{split}
\lvert Z_x( \omega ) \rvert
& =
\lvert Z_x( \omega ) - Z_{ z_i }( \omega ) + Z_{ z_i }( \omega ) \rvert
\leq
\lvert Z_x( \omega ) - Z_{ z_i }( \omega ) \rvert + \lvert Z_{ z_i }( \omega ) \rvert
\\ &
\leq
L \delta( x, z_i ) + \lvert Z_{ z_i }( \omega ) \rvert
\leq
L r_i + \lvert Z_{ z_i }( \omega ) \rvert
.
\end{split}
\end{equation}
Combining this with~\cref{eq:sup_power}
proves that
\begin{equation}
\int_\Omega
    \sup_{ x \in E }
    \, \lvert Z_x( \omega ) \rvert^p
\, \mu( \dd \omega )
\leq
\smallsum_{ i = 1 }^N
\int_\Omega
    ( L r_i + \lvert Z_{ z_i }( \omega ) \rvert )^p
\, \mu( \dd \omega )
.
\end{equation}
\end{cproof}

\begin{lemma}
\label{lem:Lp_sup_covering_number}
Let
$ p, L, r \in ( 0, \infty ) $,
let
$ ( E, \delta ) $
be a separable metric space,
let
$ ( \Omega, \cF, \mu ) $
be a measure space,
assume
$ E \neq \emptyset $
and
$ \mu( \Omega ) \neq 0 $,
let
$ Z_x \colon \Omega \to \R $, $ x \in E $,
be $ \cF $/$ \cB( \R ) $-measurable functions,
and assume for all
$ \omega \in \Omega $,
$ x, y \in E $
that
$ \lvert Z_x( \omega ) - Z_y( \omega ) \rvert \leq L \delta( x, y ) $.
Then
\begin{equation}
\int_\Omega
    \sup_{ x \in E }
    \, \lvert Z_x( \omega ) \rvert^p
\, \mu( \dd \omega )
\leq
\cC_{ ( E, \delta ), r }
\biggl[
\sup_{ x \in E }
\int_\Omega
    ( L r + \lvert Z_{ x }( \omega ) \rvert )^p
\, \mu( \dd \omega )
\biggr]
\end{equation}
(cf.~\cref{def:covering_number,lem:measurability_sup}).
\end{lemma}
\begin{cproof}{lem:Lp_sup_covering_number}
Throughout this proof
assume w.l.o.g.\ that
$ \cC_{ ( E, \delta ), r } < \infty $,
let
$ N \in \N $
be given by
$ N = \cC_{ ( E, \delta ), r } $,
and
let
$ z_1, z_2, \ldots, z_N \in E $
satisfy
$ E \subseteq
\bigcup_{ i = 1 }^N
\{ x \in E \colon \delta( x, z_i ) \leq r \} $.
Note that
\cref{lem:Lp_sup_random_field}
(with
$ r_1 \leftarrow r $,
$ r_2 \leftarrow r $,
\ldots,
$ r_N \leftarrow r $
in the notation of \cref{lem:Lp_sup_random_field})
establishes that
\begin{equation}
\begin{split}
& \int_\Omega
    \sup_{ x \in E }
    \, \lvert Z_x( \omega ) \rvert^p
\, \mu( \dd \omega )
\leq
\smallsum_{ i = 1 }^N
\int_\Omega
    ( L r + \lvert Z_{ z_i }( \omega ) \rvert )^p
\, \mu( \dd \omega )
\\ &
\leq
\smallsum_{ i = 1 }^N
\biggl[
\sup_{ x \in E }
\int_\Omega
    ( L r + \lvert Z_{ x }( \omega ) \rvert )^p
\, \mu( \dd \omega )
\biggr]
=
N
\biggl[
\sup_{ x \in E }
\int_\Omega
    ( L r + \lvert Z_{ x }( \omega ) \rvert )^p
\, \mu( \dd \omega )
\biggr]
.
\end{split}
\end{equation}
\end{cproof}

\begin{lemma}
\label{lem:Lp_sup_centred}
Let
$ p \in [ 1, \infty ) $,
$ L, r \in ( 0, \infty ) $,
let
$ ( E, \delta ) $
be a separable metric space,
assume
$ E \neq \emptyset $,
let
$ ( \Omega, \cF, \P ) $
be a probability space,
let
$ Z_x \colon \Omega \to \R $, $ x \in E $,
be random variables,
and assume for all
$ x, y \in E $
that
$ \E[ \lvert Z_x \rvert ] < \infty $
and
$ \lvert Z_x - Z_y \rvert \leq L \delta( x, y ) $.
Then
\begin{equation}
\bigl(
\E\bigl[
    \sup\nolimits_{ x \in E }
    \lvert Z_x - \E[ Z_x ] \rvert^p
\bigr]
\bigr)^{ \nicefrac{1}{p} }
\leq
( \cC_{ ( E, \delta ), r } )^{ \nicefrac{1}{p} }
\Bigl[
2 L r
+
\sup\nolimits_{ x \in E }
\bigl(
\E\bigl[
    \lvert Z_x - \E[ Z_x ] \rvert^p
\bigr]
\bigr)^{ \nicefrac{1}{p} }
\Bigr]
\end{equation}
(cf.~\cref{def:covering_number,lem:measurability_sup_centred}).
\end{lemma}
\begin{cproof}{lem:Lp_sup_centred}
Throughout this proof
let
$ Y_x \colon \Omega \to \R $, $ x \in E $,
satisfy for all
$ x \in E $,
$ \omega \in \Omega $
that
$ Y_x( \omega ) = Z_x( \omega ) - \E[ Z_x ] $.
Note that it holds for all
$ \omega \in \Omega $,
$ x, y \in E $
that
\begin{equation}
\begin{split}
\lvert Y_x( \omega ) - Y_y( \omega ) \rvert
& =
\lvert
    ( Z_x( \omega ) - \E[ Z_x ] )
    -
    ( Z_y( \omega ) - \E[ Z_y ] )
\rvert
\\ &
\leq
\lvert
    Z_x( \omega ) - Z_y( \omega )
\rvert
+
\lvert
    \E[ Z_x ] - \E[ Z_y ]
\rvert
\leq
L \delta( x, y )
+
\E[ \lvert Z_x - Z_y \rvert ]
\\ &
\leq
2 L \delta( x, y )
.
\end{split}
\end{equation}
Combining this
with \cref{lem:Lp_sup_covering_number}
(with
$ L \leftarrow 2 L $,
$ ( \Omega, \cF, \mu ) \leftarrow ( \Omega, \cF, \P ) $,
$ ( Z_x )_{ x \in E }
\leftarrow
( Y_x )_{ x \in E } $
in the notation of \cref{lem:Lp_sup_covering_number})
implies that
\begin{equation}
\begin{split}
& \bigl(
\E\bigl[
    \sup\nolimits_{ x \in E }
    \lvert Z_x - \E[ Z_x ] \rvert^p
\bigr]
\bigr)^{ \nicefrac{1}{p} }
=
\bigl(
\E\bigl[
    \sup\nolimits_{ x \in E }
    \lvert Y_x \rvert^p
\bigr]
\bigr)^{ \nicefrac{1}{p} }
\\ &
\leq
( \cC_{ ( E, \delta ), r } )^{ \nicefrac{1}{p} }
\Bigl[
\sup\nolimits_{ x \in E }
\bigl(
\E\bigl[
    ( 2 L r + \lvert Y_{ x } \rvert )^p
\bigr]
\bigr)^{ \nicefrac{1}{p} }
\Bigr]
\\ &
\leq
( \cC_{ ( E, \delta ), r } )^{ \nicefrac{1}{p} }
\Bigl[
2 L r
+
\sup\nolimits_{ x \in E }
\bigl(
\E\bigl[
    \lvert Y_{ x } \rvert^p
\bigr]
\bigr)^{ \nicefrac{1}{p} }
\Bigr]
\\ &
=
( \cC_{ ( E, \delta ), r } )^{ \nicefrac{1}{p} }
\Bigl[
2 L r
+
\sup\nolimits_{ x \in E }
\bigl(
\E\bigl[
    \lvert Z_x - \E[ Z_x ] \rvert^p
\bigr]
\bigr)^{ \nicefrac{1}{p} }
\Bigr]
.
\end{split}
\end{equation}
\end{cproof}

\begin{lemma}
\label{lem:Lp_sup_MC}
Let
$ M \in \N $,
$ p \in [ 2, \infty ) $,
$ L, r \in ( 0, \infty ) $,
let
$ ( E, \delta ) $
be a separable metric space,
assume
$ E \neq \emptyset $,
let
$ ( \Omega, \cF, \P ) $
be a probability space,
for every
$ x \in E $
let
$ Y_{ x, j } \colon \Omega \to \R $,
$ j \in \{ 1, 2, \ldots, M \} $,
be independent random variables,
assume for all
$ x, y \in E $,
$ j \in \{ 1, 2, \ldots, M \} $
that
$ \E[ \lvert Y_{ x, j } \rvert ] < \infty $
and
$ \lvert Y_{ x, j } - Y_{ y, j } \rvert \leq L \delta( x, y ) $,
and
let
$ Z_x \colon \Omega \to \R $, $ x \in E $,
satisfy for all
$ x \in E $
that
\begin{equation}
Z_x
=
\frac{1}{M}
\biggl[
\smallsum_{j=1}^M
    Y_{ x, j }
\biggr]
.
\end{equation}
Then
\begin{enumerate}[(i)]
\item
\label{item:lem:Lp_sup_MC:1}
it holds for all
$ x \in E $
that
$ \E[ \lvert Z_x \rvert ] < \infty $,
\item
\label{item:lem:Lp_sup_MC:2}
it holds that the function
$ \Omega \ni \omega \mapsto
\sup\nolimits_{ x \in E } \lvert Z_x( \omega ) - \E[ Z_x ] \rvert
\in [ 0, \infty ] $
is $ \cF $/$ \cB( [ 0, \infty ] ) $-measurable,
and
\item
\label{item:lem:Lp_sup_MC:3}
it holds that
\begin{equation}
\begin{split}
&
\bigl(
\E\bigl[
    \sup\nolimits_{ x \in E }
    \lvert Z_x - \E[ Z_x ] \rvert^p
\bigr]
\bigr)^{ \nicefrac{1}{p} }
\\ &
\leq
2
( \cC_{ ( E, \delta ), r } )^{ \nicefrac{1}{p} }
\Bigl[
L r
+
\tfrac{ \sqrt{ p - 1 } }{ \sqrt{M} }
\Bigl(
\sup\nolimits_{ x \in E }
\max\nolimits_{ j \in \{ 1, 2, \ldots, M \} }
    \bigl(
    \E\bigl[
        \lvert Y_{ x, j } - \E[ Y_{ x, j } ] \rvert^p
    \bigr]
    \bigr)^{ \nicefrac{1}{p} }
\Bigr)
\Bigr]
\end{split}
\end{equation}
(cf.~\cref{def:covering_number}).
\end{enumerate}
\end{lemma}
\begin{cproof}{lem:Lp_sup_MC}
Note that
the assumption that
$ \forall \, x \in E, \, j \in \{ 1, 2, \ldots, M \} \colon
$\linebreak$
\E[ \lvert Y_{ x, j } \rvert ] < \infty $
implies for all
$ x \in E $
that
\begin{equation}
\label{eq:max_expectation_finite}
\E[ \lvert Z_x \rvert ]
=
\E\biggl[
\frac{1}{M}
\biggl\lvert
\smallsum_{j=1}^M
    Y_{ x, j }
\biggr\rvert
\biggr]
\leq
\frac{1}{M}
\biggl[
\smallsum_{j=1}^M
    \E[ \lvert Y_{ x, j } \rvert ]
\biggr]
\leq
\max_{ j \in \{ 1, 2, \ldots, M \} }
    \E[ \lvert Y_{ x, j } \rvert ]
< \infty
.
\end{equation}
This proves~\cref{item:lem:Lp_sup_MC:1}.
Next observe that
the assumption that
$ \forall \, x, y \in E, \, j \in \{ 1, 2, \ldots, M \} \colon
\allowbreak
\lvert Y_{ x, j } - Y_{ y, j } \rvert \leq L \delta( x, y ) $
demonstrates for all
$ x, y \in E $
that
\begin{equation}
\label{eq:Lipschitz_Z}
\lvert Z_x - Z_y \rvert
=
\frac{1}{M}
\biggl\lvert
\biggl[
\smallsum_{j=1}^M
    Y_{ x, j }
\biggr]
-
\biggl[
\smallsum_{j=1}^M
    Y_{ y, j }
\biggr]
\biggr\rvert
\leq
\frac{1}{M}
\biggl[
\smallsum_{j=1}^M
    \lvert Y_{ x, j } - Y_{ y, j } \rvert
\biggr]
\leq
L \delta( x, y )
.
\end{equation}
Combining this with~\cref{item:lem:Lp_sup_MC:1}
and
\cref{lem:measurability_sup_centred}
establishes~\cref{item:lem:Lp_sup_MC:2}.
It thus remains to show~\cref{item:lem:Lp_sup_MC:3}.
For this note that
\cref{item:lem:Lp_sup_MC:1},
\cref{eq:Lipschitz_Z},
and
\cref{lem:Lp_sup_centred}
yield that
\begin{equation}
\label{eq:Lp_sup_centred}
\bigl(
\E\bigl[
    \sup\nolimits_{ x \in E }
    \lvert Z_x - \E[ Z_x ] \rvert^p
\bigr]
\bigr)^{ \nicefrac{1}{p} }
\leq
( \cC_{ ( E, \delta ), r } )^{ \nicefrac{1}{p} }
\Bigl[
2 L r
+
\sup\nolimits_{ x \in E }
\bigl(
\E\bigl[
    \lvert Z_x - \E[ Z_x ] \rvert^p
\bigr]
\bigr)^{ \nicefrac{1}{p} }
\Bigr]
.
\end{equation}
Moreover,
\cref{eq:max_expectation_finite}
and
\cref{cor:MC_Lp}
(with
$ d \leftarrow 1 $,
$ ( X_j )_{ j \in \{ 1, 2, \ldots, M \} } \leftarrow ( Y_{ x, j } )_{ j \in \{ 1, 2, \ldots, M \} } $
for $ x \in E $
in the notation of \cref{cor:MC_Lp})
prove for all
$ x \in E $
that
\begin{equation}
\begin{split}
\bigl(
\E\bigl[
    \lvert Z_x - \E[ Z_x ] \rvert^p
\bigr]
\bigr)^{ \nicefrac{1}{p} }
& =
\biggl(
\E\biggl[
    \biggl\lvert
    \frac{1}{M}
    \biggl[
    \smallsum_{j=1}^M
        Y_{ x, j }
    \biggr]
    -
    \E\biggl[
    \frac{1}{M}
    \smallsum_{j=1}^M
        Y_{ x, j }
    \biggr]
    \biggr\rvert^p
\biggr]
\biggr)^{ \!\! \nicefrac{1}{p} }
\\ &
\leq
\frac{ 2 \sqrt{ p - 1 } }{ \sqrt{M} }
\biggl[
\max_{ j \in \{ 1, 2, \ldots, M \} }
    \bigl(
    \E\bigl[
        \lvert Y_{ x, j } - \E[ Y_{ x, j } ] \rvert^p
    \bigr]
    \bigr)^{ \nicefrac{1}{p} }
\biggr]
.
\end{split}
\end{equation}
This
and
\cref{eq:Lp_sup_centred}
imply that
\begin{equation}
\begin{split}
&
\bigl(
\E\bigl[
    \sup\nolimits_{ x \in E }
    \lvert Z_x - \E[ Z_x ] \rvert^p
\bigr]
\bigr)^{ \nicefrac{1}{p} }
\\ &
\leq
( \cC_{ ( E, \delta ), r } )^{ \nicefrac{1}{p} }
\Bigl[
2 L r
+
\tfrac{ 2 \sqrt{ p - 1 } }{ \sqrt{M} }
\Bigl(
\sup\nolimits_{ x \in E }
\max\nolimits_{ j \in \{ 1, 2, \ldots, M \} }
    \bigl(
    \E\bigl[
        \lvert Y_{ x, j } - \E[ Y_{ x, j } ] \rvert^p
    \bigr]
    \bigr)^{ \nicefrac{1}{p} }
\Bigr)
\Bigr]
\\ &
=
2
( \cC_{ ( E, \delta ), r } )^{ \nicefrac{1}{p} }
\Bigl[
L r
+
\tfrac{ \sqrt{ p - 1 } }{ \sqrt{M} }
\Bigl(
\sup\nolimits_{ x \in E }
\max\nolimits_{ j \in \{ 1, 2, \ldots, M \} }
    \bigl(
    \E\bigl[
        \lvert Y_{ x, j } - \E[ Y_{ x, j } ] \rvert^p
    \bigr]
    \bigr)^{ \nicefrac{1}{p} }
\Bigr)
\Bigr]
.
\end{split}
\end{equation}
\end{cproof}

\begin{corollary}
\label{cor:Lp_sup_MC}
Let
$ M \in \N $,
$ p \in [ 2, \infty ) $,
$ L, C \in ( 0, \infty ) $,
let
$ ( E, \delta ) $
be a separable metric space,
assume
$ E \neq \emptyset $,
let
$ ( \Omega, \cF, \P ) $
be a probability space,
for every
$ x \in E $
let
$ Y_{ x, j } \colon \Omega \to \R $,
$ j \in \{ 1, 2, \ldots, M \} $,
be independent random variables,
assume for all
$ x, y \in E $,
$ j \in \{ 1, 2, \ldots, M \} $
that
$ \E[ \lvert Y_{ x, j } \rvert ] < \infty $
and
$ \lvert Y_{ x, j } - Y_{ y, j } \rvert \leq L \delta( x, y ) $,
and
let
$ Z_x \colon \Omega \to \R $, $ x \in E $,
satisfy for all
$ x \in E $
that
\begin{equation}
Z_x
=
\frac{1}{M}
\biggl[
\smallsum_{j=1}^M
    Y_{ x, j }
\biggr]
.
\end{equation}
Then
\begin{enumerate}[(i)]
\item
\label{item:cor:Lp_sup_MC:1}
it holds for all
$ x \in E $
that
$ \E[ \lvert Z_x \rvert ] < \infty $,
\item
\label{item:cor:Lp_sup_MC:2}
it holds that the function
$ \Omega \ni \omega \mapsto
\sup\nolimits_{ x \in E } \lvert Z_x( \omega ) - \E[ Z_x ] \rvert
\in [ 0, \infty ] $
is $ \cF $/$ \cB( [ 0, \infty ] ) $-measurable,
and
\item
\label{item:cor:Lp_sup_MC:3}
it holds that
\begin{equation}
\begin{split}
&
\bigl(
\E\bigl[
    \sup\nolimits_{ x \in E }
    \lvert Z_x - \E[ Z_x ] \rvert^p
\bigr]
\bigr)^{ \nicefrac{1}{p} }
\\ &
\leq
\tfrac{ 2 \sqrt{ p - 1 } }{ \sqrt{M} }
\Bigl( \cC_{ ( E, \delta ), \frac{ C \sqrt{ p - 1 } }{ L \sqrt{M} } } \Bigr)^{ \! \nicefrac{1}{p} }
\Bigl[
C
+
\sup\nolimits_{ x \in E }
\max\nolimits_{ j \in \{ 1, 2, \ldots, M \} }
    \bigl(
    \E\bigl[
        \lvert Y_{ x, j } - \E[ Y_{ x, j } ] \rvert^p
    \bigr]
    \bigr)^{ \nicefrac{1}{p} }
\Bigr]
\end{split}
\end{equation}
(cf.~\cref{def:covering_number}).
\end{enumerate}
\end{corollary}
\begin{cproof2}{cor:Lp_sup_MC}
Note that
\cref{lem:Lp_sup_MC}
shows
\cref{item:cor:Lp_sup_MC:1,item:cor:Lp_sup_MC:2}.
In addition,
Lem\-ma~\ref{lem:Lp_sup_MC}
(with
$ r \leftarrow \nicefrac{ C \sqrt{ p - 1 } }{ ( L \sqrt{M} ) } $
in the notation of
\cref{lem:Lp_sup_MC})
ensures that
\begin{align*}
&
\bigl(
\E\bigl[
    \sup\nolimits_{ x \in E }
    \lvert Z_x - \E[ Z_x ] \rvert^p
\bigr]
\bigr)^{ \nicefrac{1}{p} }
\\ &
\leq
2
\Bigl( \cC_{ ( E, \delta ), \frac{ C \sqrt{ p - 1 } }{ L \sqrt{M} } } \Bigr)^{ \! \nicefrac{1}{p} }
\Bigl[
L \tfrac{ C \sqrt{ p - 1 } }{ L \sqrt{M} }
+
\tfrac{ \sqrt{ p - 1 } }{ \sqrt{M} }
\Bigl(
\sup\nolimits_{ x \in E }
\max\nolimits_{ j \in \{ 1, 2, \ldots, M \} }
    \bigl(
    \E\bigl[
        \lvert Y_{ x, j } - \E[ Y_{ x, j } ] \rvert^p
    \bigr]
    \bigr)^{ \nicefrac{1}{p} }
\Bigr)
\Bigr]
\\ & \yesnumber
=
\tfrac{ 2 \sqrt{ p - 1 } }{ \sqrt{M} }
\Bigl( \cC_{ ( E, \delta ), \frac{ C \sqrt{ p - 1 } }{ L \sqrt{M} } } \Bigr)^{ \! \nicefrac{1}{p} }
\Bigl[
C
+
\sup\nolimits_{ x \in E }
\max\nolimits_{ j \in \{ 1, 2, \ldots, M \} }
    \bigl(
    \E\bigl[
        \lvert Y_{ x, j } - \E[ Y_{ x, j } ] \rvert^p
    \bigr]
    \bigr)^{ \nicefrac{1}{p} }
\Bigr]
.
\end{align*}
This establishes~\cref{item:cor:Lp_sup_MC:3}
and thus completes
\end{cproof2}

\subsection{Strong convergence rates for the generalisation error}
\label{sec:strong_rates_generalisation_error}

\begin{lemma}
\label{lem:abstract_generalisation_error}
Let
$ M \in \N $,
$ p \in [ 2, \infty ) $,
$ L, C, b \in ( 0, \infty ) $,
let
$ ( E, \delta ) $
be a separable metric space,
assume
$ E \neq \emptyset $,
let
$ ( \Omega, \cF, \P ) $
be a probability space,
let
$ X_{ x, j } \colon \Omega \to \R $,
$ j \in \{ 1, 2, \ldots, M \} $,
$ x \in E $,
and
$ Y_j \colon \Omega \to \R $,
$ j \in \{ 1, 2, \ldots, M \} $,
be functions,
assume for every
$ x \in E $
that
$ ( X_{ x, j }, Y_j ) $,
$ j \in \{ 1, 2, \ldots, M \} $,
are i.i.d.\ random variables,
assume for all
$ x, y \in E $,
$ j \in \{ 1, 2, \ldots, M \} $
that
$ \lvert X_{ x, j } - Y_j \rvert \leq b $
and
$ \lvert X_{ x, j } - X_{ y, j } \rvert \leq L \delta( x, y ) $,
let
$ \bfR \colon E \to [ 0, \infty ) $
satisfy for all
$ x \in E $
that
$ \bfR( x )
= \E[ \lvert X_{ x, 1 } - Y_1 \rvert^2 ] $,
and
let
$ \cR \colon E \times \Omega \to [ 0, \infty ) $
satisfy for all
$ x \in E $,
$ \omega \in \Omega $
that
\begin{equation}
\cR( x, \omega )
=
\frac{1}{M}
\biggl[
\smallsum_{j=1}^M
    \lvert X_{ x, j }( \omega ) - Y_j( \omega ) \rvert^2
\biggr]
.
\end{equation}
Then
\begin{enumerate}[(i)]
\item
\label{item:lem:abstract_generalisation_error:1}
it holds that the function
$ \Omega \ni \omega \mapsto
\sup\nolimits_{ x \in E } \lvert \cR( x, \omega ) - \bfR( x ) \rvert
\in [ 0, \infty ] $
is $ \cF $/$ \cB( [ 0, \infty ] ) $-measurable
and
\item
\label{item:lem:abstract_generalisation_error:2}
it holds that
\begin{equation}
\bigl(
\E\bigl[
    \sup\nolimits_{ x \in E }
    \lvert \cR( x ) - \bfR( x ) \rvert^p
\bigr]
\bigr)^{ \nicefrac{1}{p} }
\leq
\Bigl( \cC_{ ( E, \delta ), \frac{ C b \sqrt{ p - 1 } }{ 2 L \sqrt{M} } } \Bigr)^{ \! \nicefrac{1}{p} }
\biggl[ \frac{ 2 ( C + 1 ) b^2 \sqrt{ p - 1 } }{ \sqrt{M} } \biggr]
\end{equation}
(cf.~\cref{def:covering_number}).
\end{enumerate}
\end{lemma}
\begin{cproof2}{lem:abstract_generalisation_error}
Throughout this proof
let
$ \cY_{ x, j } \colon \Omega \to \R $,
$ j \in \{ 1, 2, \ldots, M \} $,
$ x \in E $,
satisfy for all
$ x \in E $,
$ j \in \{ 1, 2, \ldots, M \} $
that
$ \cY_{ x, j } = \lvert X_{ x, j } - Y_j \rvert^2 $.
Note that
the assumption that
for every
$ x \in E $
it holds that
$ ( X_{ x, j }, Y_j ) $,
$ j \in \{ 1, 2, \ldots, M \} $,
are i.i.d.\ random variables
ensures for all
$ x \in E $
that
\begin{equation}
\label{eq:expected_risk}
\E[ \cR( x ) ]
=
\frac{1}{M}
\biggl[
\smallsum_{j=1}^M
    \E\bigl[ \lvert X_{ x, j } - Y_j \rvert^2 \bigr]
\biggr]
=
\frac{ M \, \E\bigl[ \lvert X_{ x, 1 } - Y_1 \rvert^2 \bigr] }{M}
=
\bfR( x )
.
\end{equation}
Furthermore,
the assumption that
$ \forall \, x \in E,
\, j \in \{ 1, 2, \ldots, M \} \colon
\lvert X_{ x, j } - Y_j \rvert \leq b $
shows for all
$ x \in E $,
$ j \in \{ 1, 2, \ldots, M \} $
that
\begin{flalign}
\label{eq:finite_expectation}
&& \E[ \lvert \cY_{ x, j } \rvert ]
& =
\E[ \lvert X_{ x, j } - Y_j \rvert^2 ]
\leq
b^2
< \infty, &&
\\
&& \cY_{ x, j } - \E[ \cY_{ x, j } ]
& =
\lvert X_{ x, j } - Y_j \rvert^2
-
\E\bigl[
    \lvert X_{ x, j } - Y_j \rvert^2
\bigr]
\leq
\lvert X_{ x, j } - Y_j \rvert^2
\leq
b^2,
\\ \nonumber
& \text{and}\hidewidth
\\ \label{eq:bound_centred2}
&& \E[ \cY_{ x, j } ] - \cY_{ x, j }
& =
\E\bigl[
    \lvert X_{ x, j } - Y_j \rvert^2
\bigr]
-
\lvert X_{ x, j } - Y_j \rvert^2
\leq
\E\bigl[
    \lvert X_{ x, j } - Y_j \rvert^2
\bigr]
\leq
b^2
.
\end{flalign}
Combining
\cref{eq:finite_expectation}--\cref{eq:bound_centred2}
implies
for all
$ x \in E $,
$ j \in \{ 1, 2, \ldots, M \} $
that
\begin{equation}
\label{eq:bound_Lp_centred}
\bigl(
\E\bigl[
    \lvert \cY_{ x, j } - \E[ \cY_{ x, j } ] \rvert^p
\bigr]
\bigr)^{ \nicefrac{1}{p} }
\leq
\bigl(
\E\bigl[
    b^{ 2p }
\bigr]
\bigr)^{ \nicefrac{1}{p} }
=
b^2
.
\end{equation}
Moreover,
note that
the assumptions that
$ \forall \, x, y \in E,
\, j \in \{ 1, 2, \ldots, M \} \colon
[ \lvert X_{ x, j } - Y_j \rvert \leq b
\text{ and }
\lvert X_{ x, j }
- X_{ y, j } \rvert \leq L \delta( x, y ) ] $
and
the fact that
$ \forall \, x_1, x_2, y \in \R \colon
( x_1 - y )^2 - ( x_2 - y)^2
= ( x_1 - x_2 )( ( x_1 - y ) + ( x_2 - y ) ) $
establish for all
$ x, y \in E $,
$ j \in \{ 1, 2, \ldots, M \} $
that
\begin{equation}
\begin{split}
\lvert \cY_{ x, j } - \cY_{ y, j } \rvert
& =
\lvert
    ( X_{ x, j } - Y_j )^2
    -
    ( X_{ y, j } - Y_j )^2
\rvert
\\ &
\leq
\lvert
    X_{ x, j } - X_{ y, j }
\rvert
( \lvert X_{ x, j } - Y_j \rvert
    +
\lvert X_{ y, j } - Y_j \rvert )
\\ &
\leq
2 b
\lvert
    X_{ x, j } - X_{ y, j }
\rvert
\leq
2 b
L \delta( x, y )
.
\end{split}
\end{equation}
Combining
this,
\cref{eq:expected_risk},
\cref{eq:finite_expectation},
and
the fact that
for every
$ x \in E $
it holds that
$ \cY_{ x, j } $, 
$ j \in \{ 1, 2, \ldots, M \} $,
are independent random variables
with
\cref{cor:Lp_sup_MC}
(with
$ L \leftarrow 2 b L $,
$ C \leftarrow C b^2 $,
$ ( Y_{ x, j } )_{ x \in E, \, j \in \{ 1, 2, \ldots, M \} }
\leftarrow
( \cY_{ x, j } )_{ x \in E, \, j \in \{ 1, 2, \ldots, M \} } $,
$ ( Z_x )_{ x \in E }
\leftarrow
( \Omega \ni \omega \mapsto \cR( x, \omega ) \in \R )_{ x \in E } $
in the notation of \cref{cor:Lp_sup_MC})
and~\cref{eq:bound_Lp_centred}
proves~\cref{item:lem:abstract_generalisation_error:1}
and
\begin{equation}
\begin{split}
&
\bigl(
\E\bigl[
    \sup\nolimits_{ x \in E }
    \lvert \cR( x ) - \bfR( x ) \rvert^p
\bigr]
\bigr)^{ \nicefrac{1}{p} }
=
\bigl(
\E\bigl[
    \sup\nolimits_{ x \in E }
    \lvert \cR( x ) - \E[ \cR( x ) ] \rvert^p
\bigr]
\bigr)^{ \nicefrac{1}{p} }
\\ &
\leq
\tfrac{ 2 \sqrt{ p - 1 } }{ \sqrt{M} }
\Bigl( \cC_{ ( E, \delta ), \frac{ C b^2 \sqrt{ p - 1 } }{ 2 b L \sqrt{M} } } \Bigr)^{ \! \nicefrac{1}{p} }
\Bigl[
C b^2
+
\sup\nolimits_{ x \in E }
\max\nolimits_{ j \in \{ 1, 2, \ldots, M \} }
    \bigl(
    \E\bigl[
        \lvert \cY_{ x, j } - \E[ \cY_{ x, j } ] \rvert^p
    \bigr]
    \bigr)^{ \nicefrac{1}{p} }
\Bigr]
\\ &
\leq
\tfrac{ 2 \sqrt{ p - 1 } }{ \sqrt{M} }
\Bigl( \cC_{ ( E, \delta ), \frac{ C b \sqrt{ p - 1 } }{ 2 L \sqrt{M} } } \Bigr)^{ \! \nicefrac{1}{p} }
[ C b^2 + b^2 ]
=
\Bigl( \cC_{ ( E, \delta ), \frac{ C b \sqrt{ p - 1 } }{ 2 L \sqrt{M} } } \Bigr)^{ \! \nicefrac{1}{p} }
\biggl[ \frac{ 2 ( C + 1 ) b^2 \sqrt{ p - 1 } }{ \sqrt{M} } \biggr]
.
\end{split}
\end{equation}
This shows~\cref{item:lem:abstract_generalisation_error:2}
and thus completes
\end{cproof2}

\begin{proposition}
\label{prop:generalisation_error}
Let
$ d, \bfd, M \in \N $,
$ L, b \in ( 0, \infty ) $,
$ \alpha \in \R $,
$ \beta \in ( \alpha, \infty ) $,
$ D \subseteq \R^d $,
let
$ ( \Omega, \cF, \P ) $
be a probability space,
let
$ X_j \colon \Omega \to D $,
$ j \in \{ 1, 2, \ldots, M \} $,
and
$ Y_j \colon \Omega \to \R $,
$ j \in \{ 1, 2, \ldots, M \} $,
be functions,
assume that
$ ( X_j, Y_j ) $,
$ j \in \{ 1, 2, \ldots, M \} $,
are i.i.d.\ random variables,
let
$ f = ( f_\theta )_{ \theta \in [ \alpha, \beta ]^\bfd } \colon
[ \alpha, \beta ]^\bfd \to C( D, \R ) $
be a function,
assume for all
$ \theta, \vartheta \in [ \alpha, \beta ]^\bfd $,
$ j \in \{ 1, 2, \ldots, M \} $,
$ x \in D $
that
$ \lvert f_\theta( X_j ) - Y_j \rvert \leq b $
and
$ \lvert f_\theta( x ) - f_\vartheta( x ) \rvert
\leq
L \lVert \theta - \vartheta \rVert_\infty $,
let
$ \bfR \colon [ \alpha, \beta ]^\bfd \to [ 0, \infty ) $
satisfy for all
$ \theta \in [ \alpha, \beta ]^\bfd $
that
$ \bfR( \theta )
= \E[ \lvert f_\theta( X_1 ) - Y_1 \rvert^2 ] $,
and
let
$ \cR \colon [ \alpha, \beta ]^\bfd \times \Omega \to [ 0, \infty ) $
satisfy for all
$ \theta \in [ \alpha, \beta ]^\bfd $,
$ \omega \in \Omega $
that
\begin{equation}
\cR( \theta, \omega )
=
\frac{1}{M}
\biggl[
\smallsum_{j=1}^M
    \lvert f_\theta( X_j( \omega ) ) - Y_j( \omega ) \rvert^2
\biggr]
\end{equation}
(cf.~\cref{def:p-norm}).
Then
\begin{enumerate}[(i)]
\item
\label{item:prop:generalisation_error:1}
it holds that the function
$ \Omega \ni \omega \mapsto
\sup\nolimits_{ \theta \in [ \alpha, \beta ]^\bfd } \lvert \cR( \theta, \omega ) - \bfR( \theta ) \rvert
\in [ 0, \infty ] $
is $ \cF $/\allowbreak$ \cB( [ 0, \infty ] ) $-measurable
and
\item
\label{item:prop:generalisation_error:2}
it holds for all
$ p \in ( 0, \infty ) $
that
\begin{equation}
\begin{split}
&
\bigl(
\E\bigl[
    \sup\nolimits_{ \theta \in [ \alpha, \beta ]^\bfd }
    \lvert \cR( \theta ) - \bfR( \theta ) \rvert^p
\bigr]
\bigr)^{ \nicefrac{1}{p} }
\\ &
\leq
\inf_{ C, \varepsilon \in ( 0, \infty ) }
\Biggl[
\frac{ 2 ( C + 1 ) b^2
    \max\{
    1,
    [ 2 \sqrt{M} L ( \beta - \alpha ) ( C b )^{ -1 } ]^{ \varepsilon }
    \}
\sqrt{ \max\{ 1, p, \nicefrac{\bfd}{ \varepsilon } \} }
}{ \sqrt{M} }
\Biggr]
\\ &
\leq
\inf_{ C \in ( 0, \infty ) }
\Biggl[
\frac{ 2 ( C + 1 ) b^2
\sqrt{ e \max\{ 1, p, \bfd \ln( 4 M L^2 ( \beta - \alpha )^2 ( C b )^{ -2 } ) \} }
}{ \sqrt{M} }
\Biggr]
.
\end{split}
\end{equation}
\end{enumerate}
\end{proposition}
\begin{cproof2}{prop:generalisation_error}
Throughout this proof
let
$ p \in ( 0, \infty ) $,
let
$ ( \kappa_C )_{ C \in ( 0, \infty ) } \subseteq ( 0, \infty ) $
satisfy for all
$ C \in ( 0, \infty ) $
that
$ \nicefrac{ 2 \sqrt{M} L ( \beta - \alpha ) }{ ( C b ) } $,
let
$ \cX_{ \theta, j } \colon \Omega \to \R $,
$ j \in \{ 1, 2, \ldots, M \} $,
$ \theta \in [ \alpha, \beta ]^\bfd $,
satisfy for all
$ \theta \in [ \alpha, \beta ]^\bfd $,
$ j \in \{ 1, 2, \ldots, M \} $
that
$ \cX_{ \theta, j } = f_\theta( X_j ) $,
and
let
$ \delta \colon ( [ \alpha, \beta ]^\bfd ) \times ( [ \alpha, \beta ]^\bfd ) \to [ 0, \infty ) $
satisfy for all
$ \theta, \vartheta \in [ \alpha, \beta ]^\bfd $
that
$ \delta( \theta, \vartheta ) = \lVert \theta - \vartheta \rVert_\infty $.
First of all,
note that
the assumption that
$ \forall \, \theta \in [ \alpha, \beta ]^\bfd,
\, j \in \{ 1, 2, \ldots, M \} \colon
\lvert f_\theta( X_j ) - Y_j \rvert \leq b $
implies
for all
$ \theta \in [ \alpha, \beta ]^\bfd $,
$ j \in \{ 1, 2, \ldots, M \} $
that
\begin{equation}
\label{eq:difference_bounded}
\lvert \cX_{ \theta, j } - Y_j \rvert
=
\lvert f_\theta( X_j ) - Y_j \rvert
\leq b
.
\end{equation}
In addition,
the assumption that
$ \forall \, \theta, \vartheta \in [ \alpha, \beta ]^\bfd,
\, x \in D \colon
\lvert f_\theta( x ) - f_\vartheta( x ) \rvert
\leq
L \lVert \theta - \vartheta \rVert_\infty $
ensures
for all
$ \theta, \vartheta \in [ \alpha, \beta ]^\bfd $,
$ j \in \{ 1, 2, \ldots, M \} $
that
\begin{equation}
\lvert \cX_{ \theta, j } - \cX_{ \vartheta, j } \rvert
=
\lvert f_\theta( X_j ) - f_\vartheta( X_j ) \rvert
\leq
\sup\nolimits_{ x \in D }
    \lvert f_\theta( x ) - f_\vartheta( x ) \rvert
\leq
L \lVert \theta - \vartheta \rVert_\infty
=
L \delta( \theta, \vartheta )
.
\end{equation}
Combining
this,
\cref{eq:difference_bounded},
and
the fact that
for every
$ \theta \in [ \alpha, \beta ]^\bfd $
it holds that
$ ( \cX_{ \theta, j }, Y_j ) $,
$ j \in \{ 1, 2, \ldots, M \} $,
are i.i.d.\ random variables
with
\cref{lem:abstract_generalisation_error}
(with
$ p \leftarrow q $,
$ C \leftarrow C $,
$ ( E, \delta ) \leftarrow ( [ \alpha, \beta ]^\bfd, \delta ) $,
$ ( X_{ x, j } )_{ x \in E, \, j \in \{ 1, 2, \ldots, M \} }
\leftarrow
( \cX_{ \theta, j } )_{ \theta \in [ \alpha, \beta ]^\bfd, \, j \in \{ 1, 2, \ldots, M \} } $
for
$ q \in [ 2, \infty ) $,
$ C \in ( 0, \infty ) $
in the notation of \cref{lem:abstract_generalisation_error})
demonstrates for all
$ C \in ( 0, \infty ) $,
$ q \in [ 2, \infty ) $
that the function
$ \Omega \ni \omega \mapsto
\sup\nolimits_{ \theta \in [ \alpha, \beta ]^\bfd } \lvert \cR( \theta, \omega ) - \bfR( \theta ) \rvert
\in [ 0, \infty ] $
is $ \cF $/$ \cB( [ 0, \infty ] ) $-measurable
and%
\begin{equation}
\label{eq:abstract_generalisation_error}
\bigl(
\E\bigl[
    \sup\nolimits_{ \theta \in [ \alpha, \beta ]^\bfd }
    \lvert \cR( \theta ) - \bfR( \theta ) \rvert^q
\bigr]
\bigr)^{ \nicefrac{1}{q} }
\leq
\Bigl(
    \cC_{ ( [ \alpha, \beta ]^\bfd, \delta ), \frac{ C b \sqrt{ q - 1 } }{ 2 L \sqrt{M} } }
\Bigr)^{ \! \nicefrac{1}{q} }
\biggl[ \frac{ 2 ( C + 1 ) b^2 \sqrt{ q - 1 } }{ \sqrt{M} } \biggr]
\end{equation}%
(cf.~\cref{def:covering_number}).
This finishes the proof of~\cref{item:prop:generalisation_error:1}.
Next observe that
\cref{item:lem:covering_number_cube:2} in \cref{lem:covering_number_cube}
(with
$ d \leftarrow \bfd $,
$ a \leftarrow \alpha $,
$ b \leftarrow \beta $,
$ r \leftarrow r $
for
$ r \in ( 0, \infty ) $
in the notation of \cref{lem:covering_number_cube})
shows for all
$ r \in ( 0, \infty ) $
that
\begin{equation}
\begin{split}
\cC_{ ( [ \alpha, \beta ]^\bfd, \delta ), r }
& \leq
\mathbbm{1}_{ [ 0, r ] }\bigl( \tfrac{ \beta - \alpha }{2} \bigr)
+
\bigl( \tfrac{ \beta - \alpha }{r} \bigr)^\bfd
\mathbbm{1}_{ ( r, \infty ) }\bigl( \tfrac{ \beta - \alpha }{2} \bigr)
\\ &
\leq
\max\Bigl\{
    1,
    \bigl( \tfrac{ \beta - \alpha }{r} \bigr)^\bfd
\Bigr\}
\bigl(
\mathbbm{1}_{ [ 0, r ] }\bigl( \tfrac{ \beta - \alpha }{2} \bigr)
+
\mathbbm{1}_{ ( r, \infty ) }\bigl( \tfrac{ \beta - \alpha }{2} \bigr)
\bigr)
\\ &
=
\max\Bigl\{
    1,
    \bigl( \tfrac{ \beta - \alpha }{r} \bigr)^\bfd
\Bigr\}
.
\end{split}
\end{equation}
This yields for all
$ C \in ( 0, \infty ) $,
$ q \in [ 2, \infty ) $
that
\begin{equation}
\begin{split}
\Bigl(
\cC_{ ( [ \alpha, \beta ]^\bfd, \delta ), \frac{ C b \sqrt{ q - 1 } }{ 2 L \sqrt{M} } }
\Bigr)^{ \! \nicefrac{1}{ q } }
&
\leq
\max\biggl\{
    1,
    \Bigl( \tfrac{ 2 ( \beta - \alpha ) L \sqrt{M} }{ C b \sqrt{ q - 1 } } \Bigr)^{ \! \frac{\bfd}{q} }
\biggr\}
\\ &
\leq
\max\biggl\{
    1,
    \Bigl( \tfrac{ 2 ( \beta - \alpha ) L \sqrt{M} }{ C b } \Bigr)^{ \! \frac{\bfd}{q} }
\biggr\}
=
\max\Bigl\{
    1,
    ( \kappa_C )^{ \frac{\bfd}{q} }
\Bigr\}
.
\end{split}
\end{equation}
Jensen's inequality
and
\cref{eq:abstract_generalisation_error}
hence
prove for all
$ C, \varepsilon \in ( 0, \infty ) $
that
\begin{equation}
\label{eq:generalisation_error_estimate}
\begin{split}
&
\bigl(
\E\bigl[
    \sup\nolimits_{ \theta \in [ \alpha, \beta ]^\bfd }
    \lvert \cR( \theta ) - \bfR( \theta ) \rvert^p
\bigr]
\bigr)^{ \nicefrac{1}{p} }
\\ &
\leq
\bigl(
\E\bigl[
    \sup\nolimits_{ \theta \in [ \alpha, \beta ]^\bfd }
    \lvert \cR( \theta ) - \bfR( \theta ) \rvert^{ \max\{ 2, p, \nicefrac{\bfd}{ \varepsilon } \} }
\bigr]
\bigr)^{ \frac{1}{ \max\{ 2, p, \nicefrac{\bfd}{ \varepsilon } \} } }
\\ &
\leq
\max\Bigl\{
    1,
    ( \kappa_C )^{ \frac{\bfd}{ \max\{ 2, p, \nicefrac{\bfd}{ \varepsilon } \} } }
\Bigr\}
\frac{ 2 ( C + 1 ) b^2 \sqrt{ \max\{ 2, p, \nicefrac{\bfd}{ \varepsilon } \} - 1 } }{ \sqrt{M} }
\\ &
=
\max\bigl\{
    1,
    ( \kappa_C )^{ \min\{ \nicefrac{\bfd}{2}, \nicefrac{\bfd}{p}, \varepsilon \} }
\bigr\}
\frac{ 2 ( C + 1 ) b^2 \sqrt{ \max\{ 1, p - 1, \nicefrac{\bfd}{ \varepsilon } - 1 \} } }{ \sqrt{M} }
\\ &
\leq
\frac{
2 ( C + 1 ) b^2
\max\{ 1, ( \kappa_C )^{ \varepsilon } \}
\sqrt{ \max\{ 1, p, \nicefrac{\bfd}{ \varepsilon } \} }
}{ \sqrt{M} }
.
\end{split}
\end{equation}
Next note that
the fact that
$ \forall \, a \in ( 1, \infty ) \colon
a^{ \nicefrac{1}{ ( 2 \ln( a ) ) }}
= e^{ \nicefrac{ \ln( a ) }{ ( 2 \ln( a ) ) } }
= e^{ \nicefrac{1}{2} }
= \sqrt{e}
\geq 1 $
ensures for all
$ C \in ( 0, \infty ) $
with $ \kappa_C > 1 $
that
\begin{equation}
\label{eq:generalisation_error_epsilon_choice}
\begin{split}
&
\inf_{ \varepsilon \in ( 0, \infty ) }
\Biggl[
\frac{ 2 ( C + 1 ) b^2
    \max\{
    1,
    ( \kappa_C )^{ \varepsilon }
    \}
\sqrt{ \max\{ 1, p, \nicefrac{\bfd}{ \varepsilon } \} }
}{ \sqrt{M} }
\Biggr]
\\ &
\leq
\frac{ 2 ( C + 1 ) b^2
    \max\{
    1,
    ( \kappa_C )^{ \nicefrac{1}{ ( 2 \ln( \kappa_C ) ) } }
    \}
\sqrt{ \max\{ 1, p, 2 \bfd \ln( \kappa_C ) \} }
}{ \sqrt{M} }
\\ &
=
\frac{ 2 ( C + 1 ) b^2
\sqrt{ e \max\{ 1, p, \bfd \ln( [ \kappa_C ]^2 ) \} }
}{ \sqrt{M} }
.
\end{split}
\end{equation}
In addition,
observe that it holds
for all
$ C \in ( 0, \infty ) $
with $ \kappa_C \leq 1 $
that
\begin{equation}
\label{eq:generalisation_error_epsilon_auxiliary}
\begin{split}
&
\inf_{ \varepsilon \in ( 0, \infty ) }
\Biggl[
\frac{ 2 ( C + 1 ) b^2
    \max\{
    1,
    ( \kappa_C )^{ \varepsilon }
    \}
\sqrt{ \max\{ 1, p, \nicefrac{\bfd}{ \varepsilon } \} }
}{ \sqrt{M} }
\Biggr]
\\ &
\leq
\inf_{ \varepsilon \in ( 0, \infty ) }
\Biggl[
\frac{ 2 ( C + 1 ) b^2
\sqrt{ \max\{ 1, p, \nicefrac{\bfd}{ \varepsilon } \} }
}{ \sqrt{M} }
\Biggr]
\leq
\frac{ 2 ( C + 1 ) b^2
\sqrt{ \max\{ 1, p \} }
}{ \sqrt{M} }
\\ &
\leq
\frac{ 2 ( C + 1 ) b^2
\sqrt{ e \max\{ 1, p, \bfd \ln( [ \kappa_C ]^2 ) \} }
}{ \sqrt{M} }
.
\end{split}
\end{equation}
Combining~\cref{eq:generalisation_error_estimate}
with~\cref{eq:generalisation_error_epsilon_choice}
and~\cref{eq:generalisation_error_epsilon_auxiliary}
demonstrates that
\begin{equation}
\begin{split}
&
\bigl(
\E\bigl[
    \sup\nolimits_{ \theta \in [ \alpha, \beta ]^\bfd }
    \lvert \cR( \theta ) - \bfR( \theta ) \rvert^p
\bigr]
\bigr)^{ \nicefrac{1}{p} }
\\ &
\leq
\inf_{ C, \varepsilon \in ( 0, \infty ) }
\Biggl[
\frac{ 2 ( C + 1 ) b^2
    \max\{
    1,
    ( \kappa_C )^{ \varepsilon }
    \}
\sqrt{ \max\{ 1, p, \nicefrac{\bfd}{ \varepsilon } \} }
}{ \sqrt{M} }
\Biggr]
\\ &
=
\inf_{ C, \varepsilon \in ( 0, \infty ) }
\Biggl[
\frac{ 2 ( C + 1 ) b^2
    \max\{
    1,
    [ 2 \sqrt{M} L ( \beta - \alpha ) ( C b )^{ -1 } ]^{ \varepsilon }
    \}
\sqrt{ \max\{ 1, p, \nicefrac{\bfd}{ \varepsilon } \} }
}{ \sqrt{M} }
\Biggr]
\\ &
\leq
\inf_{ C \in ( 0, \infty ) }
\Biggl[
\frac{ 2 ( C + 1 ) b^2
\sqrt{ e \max\{ 1, p, \bfd \ln( [ \kappa_C ]^2 ) \} }
}{ \sqrt{M} }
\Biggr]
\\ &
=
\inf_{ C \in ( 0, \infty ) }
\Biggl[
\frac{ 2 ( C + 1 ) b^2
\sqrt{ e \max\{ 1, p, \bfd \ln( 4 M L^2 ( \beta - \alpha )^2 ( C b )^{ -2 } ) \} }
}{ \sqrt{M} }
\Biggr]
.
\end{split}
\end{equation}
This establishes~\cref{item:prop:generalisation_error:2}
and thus completes
\end{cproof2}

\begin{corollary}
\label{cor:generalisation_error}
Let
$ d, \bfd, \bfL, M \in \N $,
$ B, b \in [ 1, \infty ) $,
$ u \in \R $,
$ v \in [ u + 1, \infty ) $,
$ \bfl = ( \bfl_0, \bfl_1, \ldots, \bfl_\bfL ) \in \N^{ \bfL + 1 } $,
$ D \subseteq [ -b, b ]^d $,
assume
$ \bfl_0 = d $,
$ \bfl_\bfL = 1 $,
and
$ \bfd \geq \sum_{i=1}^{\bfL} \bfl_i( \bfl_{ i - 1 } + 1 ) $,
let
$ ( \Omega, \cF, \P ) $
be a probability space,
let
$ X_j \colon \Omega \to D $,
$ j \in \{ 1, 2, \ldots, M \} $,
and
$ Y_j \colon \Omega \to [ u, v ] $,
$ j \in \{ 1, 2, \ldots, M \} $,
be functions,
assume that
$ ( X_j, Y_j ) $,
$ j \in \{ 1, 2, \ldots, M \} $,
are i.i.d.\ random variables,
let
$ \bfR \colon [ -B, B ]^\bfd \to [ 0, \infty ) $
satisfy for all
$ \theta \in [ -B, B ]^\bfd $
that
$ \bfR( \theta )
= \E[ \lvert \clippedNN{\theta}{\bfl}{u}{v}( X_1 ) - Y_1 \rvert^2 ] $,
and
let
$ \cR \colon [ -B, B ]^\bfd \times \Omega \to [ 0, \infty ) $
satisfy for all
$ \theta \in [ -B, B ]^\bfd $,
$ \omega \in \Omega $
that
\begin{equation}
\cR( \theta, \omega )
=
\frac{1}{M}
\biggl[
\smallsum_{j=1}^M
    \lvert \clippedNN{\theta}{\bfl}{u}{v}( X_j( \omega ) ) - Y_j( \omega ) \rvert^2
\biggr]
\end{equation}
(cf.~\cref{def:clipped_NN}).
Then
\begin{enumerate}[(i)]
\item
\label{item:cor:generalisation_error:1}
it holds that the function
$ \Omega \ni \omega \mapsto
\sup\nolimits_{ \theta \in [ -B, B ]^\bfd } \lvert \cR( \theta, \omega ) - \bfR( \theta ) \rvert
\in [ 0, \infty ] $
is $ \cF $/\allowbreak$ \cB( [ 0, \infty ] ) $-measurable
and
\item
\label{item:cor:generalisation_error:2}
it holds for all
$ p \in ( 0, \infty ) $
that
\begin{equation}
\begin{split}
& \bigl(
\E\bigl[
    \sup\nolimits_{ \theta \in [ -B, B ]^\bfd }
    \lvert \cR( \theta ) - \bfR( \theta ) \rvert^p
\bigr]
\bigr)^{ \nicefrac{1}{p} }
\\ &
\leq
\frac{
    9 ( v - u )^2
    \bfL ( \lVert \bfl \rVert_{ \infty } + 1 )
    \sqrt{
        \max\{
            p,
            \ln( 4 ( M b )^{ \nicefrac{1}{ \bfL } } ( \lVert \bfl \rVert_{ \infty } + 1 ) B )
        \}
    }
}{ \sqrt{M} }
\\ &
\leq
\frac{
    9 ( v - u )^2
    \bfL ( \lVert \bfl \rVert_{ \infty } + 1 )^2
    \max\{
        p,
        \ln( 3 M B b )
    \}
}{ \sqrt{M} }
\end{split}
\end{equation}
(cf.~\cref{def:p-norm}).
\end{enumerate}
\end{corollary}
\begin{cproof}{cor:generalisation_error}
Throughout this proof
let
$ \fd \in \N $
be given by
$ \fd = \sum_{i=1}^{\bfL} \bfl_i( \bfl_{ i - 1 } + 1 ) $,
let
$ L \in ( 0, \infty ) $
be given by
$ L =
b \bfL ( \lVert \bfl \rVert_{ \infty } + 1 )^\bfL B^{ \bfL - 1 } $,
let
$ f = ( f_\theta )_{ \theta \in [ -B, B ]^\fd } \colon
[ -B, B ]^\fd \to C( D, \R ) $
satisfy for all
$ \theta \in [ -B, B ]^\fd $,
$ x \in D $
that
$ f_\theta( x ) = \clippedNN{\theta}{\bfl}{u}{v}( x ) $,
let
$ \scrR \colon [ -B, B ]^\fd \to [ 0, \infty ) $
satisfy for all
$ \theta \in [ -B, B ]^\fd $
that
$ \scrR( \theta )
= \E[ \lvert f_\theta( X_1 ) - Y_1 \rvert^2 ]
= \E[ \lvert \clippedNN{\theta}{\bfl}{u}{v}( X_1 ) - Y_1 \rvert^2 ] $,
and
let
$ R \colon [ -B, B ]^\fd \times \Omega \to [ 0, \infty ) $
satisfy for all
$ \theta \in [ -B, B ]^\fd $,
$ \omega \in \Omega $
that
\begin{equation}
R( \theta, \omega )
=
\frac{1}{M}
\biggl[
\smallsum_{j=1}^M
    \lvert f_\theta( X_j( \omega ) ) - Y_j( \omega ) \rvert^2
\biggr]
=
\frac{1}{M}
\biggl[
\smallsum_{j=1}^M
    \lvert \clippedNN{\theta}{\bfl}{u}{v}( X_j( \omega ) ) - Y_j( \omega ) \rvert^2
\biggr]
.
\end{equation}
Note that
the fact that
$ \forall \, \theta \in \R^\fd,
\, x \in \R^d \colon
\clippedNN{\theta}{\bfl}{u}{v}( x ) \in [ u, v ] $
and the assumption that
$ \forall \, j \in \{ 1, 2, \ldots, M \} \colon
Y_j( \Omega ) \subseteq [ u, v ] $
imply
for all
$ \theta \in [ -B, B ]^\fd $,
$ j \in \{ 1, 2, \ldots, M \} $
that
\begin{equation}
\label{eq:difference_bounded_NN}
\lvert f_\theta( X_j ) - Y_j \rvert
=
\lvert \clippedNN{\theta}{\bfl}{u}{v}( X_j ) - Y_j \rvert
\leq
\sup\nolimits_{ y_1, y_2 \in [ u, v ] }
    \lvert y_1 - y_2 \rvert
=
v - u
.
\end{equation}
Moreover,
the assumptions that
$ D \subseteq [ -b, b ]^d $,
$ \bfl_0 = d $,
and
$ \bfl_\bfL = 1 $,
Beck, Jentzen, \& Kuckuck~\cite[Corollary~2.37]{BeckJentzenKuckuck2019arXiv}
(with
$ a \leftarrow -b $,
$ b \leftarrow b $,
$ u \leftarrow u $,
$ v \leftarrow v $,
$ d \leftarrow \fd $,
$ L \leftarrow \bfL $,
$ l \leftarrow \bfl $
in the notation of~\cite[Corollary~2.37]{BeckJentzenKuckuck2019arXiv}),
and
the assumptions that
$ b \geq 1 $ and $ B \geq 1 $
ensure for all
$ \theta, \vartheta \in [ -B, B ]^\fd $,
$ x \in D $
that
\begin{equation}
\label{eq:Lipschitz_NN}
\begin{split}
\lvert f_\theta( x ) - f_\vartheta( x ) \rvert
& \leq
\sup\nolimits_{ y \in [ -b, b ]^d }
    \lvert \clippedNN{\theta}{\bfl}{u}{v}( y ) - \clippedNN{\vartheta}{\bfl}{u}{v}( y ) \rvert
\\ &
\leq
\bfL
\max\{ 1, b \}
( \lVert \bfl \rVert_{ \infty } + 1 )^\bfL
( \max\{ 1, \lVert \theta \rVert_{ \infty }, \lVert \vartheta \rVert_{ \infty } \} )^{ \bfL - 1 }
\lVert \theta - \vartheta \rVert_{ \infty }
\\ &
\leq
b
\bfL
( \lVert \bfl \rVert_{ \infty } + 1 )^\bfL
B^{ \bfL - 1 }
\lVert \theta - \vartheta \rVert_{ \infty }
=
L \lVert \theta - \vartheta \rVert_\infty
.
\end{split}
\end{equation}
Furthermore,
the facts that
$ \bfd \geq \fd $
and
$ \forall \, \theta = ( \theta_1, \theta_2, \ldots, \theta_\bfd ) \in \R^\bfd \colon
\clippedNN{\theta}{\bfl}{u}{v}
=
\clippedNN{ \smash{ (\theta_1, \theta_2, \ldots, \theta_\fd ) } }{\bfl}{u}{v} $
prove for all
$ \omega \in \Omega $
that
\begin{equation}
\label{eq:independent_parameters}
\sup\nolimits_{ \theta \in [ -B, B ]^\bfd }
    \lvert \cR( \theta, \omega ) - \bfR( \theta ) \rvert
=
\sup\nolimits_{ \theta \in [ -B, B ]^\fd }
    \lvert R( \theta, \omega ) - \scrR( \theta ) \rvert
.
\end{equation}
Next observe that
\cref{eq:difference_bounded_NN},
\cref{eq:Lipschitz_NN},
\cref{prop:generalisation_error}
(with
$ \bfd \leftarrow \fd $,
$ b \leftarrow v - u $,
$ \alpha \leftarrow -B $,
$ \beta \leftarrow B $,
$ \bfR \leftarrow \scrR $,
$ \cR \leftarrow R $
in the notation of \cref{prop:generalisation_error}),
and
the facts that
$ v - u \geq ( u + 1 ) - u = 1 $
and
$ \fd \leq
\bfL \lVert \bfl \rVert_{ \infty }( \lVert \bfl \rVert_{ \infty } + 1 )
\leq
\bfL ( \lVert \bfl \rVert_{ \infty } + 1 )^2 $
demonstrate
for all
$ p \in ( 0, \infty ) $
that the function
$ \Omega \ni \omega \mapsto
\sup\nolimits_{ \theta \in [ -B, B ]^\fd } \lvert R( \theta, \omega ) - \scrR( \theta ) \rvert
\in [ 0, \infty ] $
is $ \cF $/$ \cB( [ 0, \infty ] ) $-measurable
and
\begin{equation}
\label{eq:prop:generalisation_error}
\begin{split}
&
\bigl(
\E\bigl[
    \sup\nolimits_{ \theta \in [ -B, B ]^\fd }
    \lvert R( \theta ) - \scrR( \theta ) \rvert^p
\bigr]
\bigr)^{ \nicefrac{1}{p} }
\\ &
\leq
\inf_{ C \in ( 0, \infty ) }
\Biggl[
\frac{ 2 ( C + 1 ) ( v - u )^2
\sqrt{ e \max\{ 1, p, \fd \ln( 4 M L^2 ( 2 B )^2 ( C [ v - u ] )^{ -2 } ) \} }
}{ \sqrt{M} }
\Biggr]
\\ &
\leq
\inf_{ C \in ( 0, \infty ) }
\Biggl[
\frac{ 2 ( C + 1 ) ( v - u )^2
\sqrt{ e \max\{ 1, p, \bfL ( \lVert \bfl \rVert_{ \infty } + 1 )^2 \ln( 2^4 M L^2 B^2 C^{ -2 } ) \} }
}{ \sqrt{M} }
\Biggr]
.
\end{split}
\end{equation}
This and~\cref{eq:independent_parameters}
establish~\cref{item:cor:generalisation_error:1}.
In addition,
combining
\cref{eq:independent_parameters}--\cref{eq:prop:generalisation_error}
with
the fact that
$ 2^6 \bfL^2
\leq 2^6 \cdot 2^{ 2 ( \bfL - 1 ) }
= 2^{ 4 + 2 \bfL }
\leq 2^{ 4 \bfL + 2 \bfL }
= 2^{ 6 \bfL } $
and
the facts that
$ 3 \geq e $,
$ B \geq 1 $,
$ \bfL \geq 1 $,
$ M \geq 1 $,
and
$ b \geq 1 $
shows for all
$ p \in ( 0, \infty ) $
that
\begin{equation}
\label{eq:generalisation_error_NN1}
\begin{split}
& \bigl(
\E\bigl[
    \sup\nolimits_{ \theta \in [ -B, B ]^\bfd }
    \lvert \cR( \theta ) - \bfR( \theta ) \rvert^p
\bigr]
\bigr)^{ \nicefrac{1}{p} }
=
\bigl(
\E\bigl[
    \sup\nolimits_{ \theta \in [ -B, B ]^\fd }
    \lvert R( \theta ) - \scrR( \theta ) \rvert^p
\bigr]
\bigr)^{ \nicefrac{1}{p} }
\\ &
\leq
\frac{ 2 ( \nicefrac{1}{2} + 1 ) ( v - u )^2
\sqrt{ e \max\{ 1, p, \bfL ( \lVert \bfl \rVert_{ \infty } + 1 )^2 \ln( 2^4 M L^2 B^2 2^{ 2 } ) \} }
}{ \sqrt{M} }
\\ &
=
\frac{ 3 ( v - u )^2
\sqrt{
    e
    \max\{
        p,
        \bfL ( \lVert \bfl \rVert_{ \infty } + 1 )^2
        \ln( 2^6 M b^2 \bfL^2 ( \lVert \bfl \rVert_{ \infty } + 1 )^{ 2 \bfL } B^{ 2 \bfL } )
    \}
    }
}{ \sqrt{M} }
\\ &
\leq
\frac{ 3 ( v - u )^2
\sqrt{
    e
    \max\{
        p,
        3 \bfL^2 ( \lVert \bfl \rVert_{ \infty } + 1 )^2
        \ln( [ 2^{ 6 \bfL } M b^2 ( \lVert \bfl \rVert_{ \infty } + 1 )^{ 2 \bfL } B^{ 2 \bfL } ]^{ \nicefrac{1}{ ( 3 \bfL ) } } )
    \}
    }
}{ \sqrt{M} }
\\ &
\leq
\frac{ 3 ( v - u )^2
\sqrt{
    3
    \max\{
        p,
        3 \bfL^2 ( \lVert \bfl \rVert_{ \infty } + 1 )^2
        \ln( 2^2 ( M b^2 )^{ \nicefrac{1}{ ( 3 \bfL ) } } ( \lVert \bfl \rVert_{ \infty } + 1 ) B )
    \}
    }
}{ \sqrt{M} }
\\ &
\leq
\frac{ 9 ( v - u )^2
\bfL ( \lVert \bfl \rVert_{ \infty } + 1 )
\sqrt{
    \max\{
        p,
        \ln( 4 ( M b )^{ \nicefrac{1}{ \bfL } } ( \lVert \bfl \rVert_{ \infty } + 1 ) B )
    \}
    }
}{ \sqrt{M} }
.
\end{split}
\end{equation}
Furthermore,
note that
the fact that
$ \forall \, n \in \N \colon
n \leq 2^{ n - 1 } $
and
the fact that
$ \lVert \bfl \rVert_{ \infty } \geq 1 $
imply that
\begin{equation}
4 ( \lVert \bfl \rVert_{ \infty } + 1 )
\leq 2^2 \cdot 2^{ ( \lVert \bfl \rVert_{ \infty } + 1 ) - 1 }
= 2^3 \cdot 2^{ ( \lVert \bfl \rVert_{ \infty } + 1 ) - 2 }
\leq 3^2 \cdot 3^{ ( \lVert \bfl \rVert_{ \infty } + 1 ) - 2 }
= 3^{ ( \lVert \bfl \rVert_{ \infty } + 1 ) }
.
\end{equation}
This
demonstrates for all
$ p \in ( 0, \infty ) $
that
\begin{equation}
\begin{split}
&
\frac{ 9 ( v - u )^2
\bfL ( \lVert \bfl \rVert_{ \infty } + 1 )
\sqrt{
    \max\{
        p,
        \ln( 4 ( M b )^{ \nicefrac{1}{ \bfL } } ( \lVert \bfl \rVert_{ \infty } + 1 ) B )
    \}
    }
}{ \sqrt{M} }
\\ &
\leq
\frac{ 9 ( v - u )^2
\bfL ( \lVert \bfl \rVert_{ \infty } + 1 )
\sqrt{
    \max\{
        p,
        ( \lVert \bfl \rVert_{ \infty } + 1 )
        \ln( [ 3^{ ( \lVert \bfl \rVert_{ \infty } + 1 ) } ( M b )^{ \nicefrac{1}{ \bfL } } B ]^{ \nicefrac{1}{ ( \lVert \bfl \rVert_{ \infty } + 1 ) } } )
    \}
    }
}{ \sqrt{M} }
\\ &
\leq
\frac{
    9 ( v - u )^2
    \bfL ( \lVert \bfl \rVert_{ \infty } + 1 )^2
    \max\{
        p,
        \ln( 3 M B b )
    \}
}{ \sqrt{M} }
.
\end{split}
\end{equation}
Combining this with~\cref{eq:generalisation_error_NN1}
shows~\cref{item:cor:generalisation_error:2}.
\end{cproof}

\section{Analysis of the optimisation error}
\label{sec:optimisation_error}

The main result of this section, \cref{prop:minimum_MC_rate},
establishes that
the optimisation error of the Minimum Monte Carlo method
applied to a Lipschitz continuous random field with
a $ \bfd $-dimensional hypercube as index set,
where $ \bfd \in \N $,
converges in the probabilistically strong sense
with rate $ \nicefrac{1}{\bfd} $
with respect to the number of samples used,
provided that the sample indices are continuous uniformly drawn
from the index hypercube
(cf.~\cref{item:prop:minimum_MC_rate:2} in \cref{prop:minimum_MC_rate}).
We refer to
Beck, Jentzen, \& Kuckuck~\cite[Lemmas~3.22--3.23]{BeckJentzenKuckuck2019arXiv}
for analogous results for convergence in probability instead of strong convergence
and
to Beck et al.~\cite[Lemma~3.5]{BeckBeckerGrohsJaafariJentzen2018arXiv}
for a related result.
\cref{cor:minimum_MC_rate} below
specialises \cref{prop:minimum_MC_rate} to the case
where the empirical risk from deep learning based empirical risk minimisation
with quadratic loss function
indexed by a hypercube of DNN parameter vectors
plays the role of the random field under consideration.
In the proof of \cref{cor:minimum_MC_rate} we make use of the elementary and well-known fact
that this choice for the random field is indeed Lipschitz continuous,
which is the assertion of \cref{lem:Lipschitz_risk}.
Further results on the optimisation error in the context of stochastic approximation
can be found, e.g., in
\cite{AroraDuHuLiWang2019,
BachMoulines2013,
BercuFort2013,
ChauMoulinesRasonyiSabanisZhang2019arXiv,
DereichKassing2019arXiv,
DereichMuller_Gronbach2019,
DuLeeLiWangZhai2019,
DuZhaiPoczosSingh2018arXiv,
FehrmanGessJentzen2019arXiv,
JentzenKuckuckNeufeldVonWurstemberger2018arXiv,
JentzenvonWurstemberger2020,
KarimiMiasojedowMoulinesWai2019arXiv,
LeiHuLiTang2019arXiv,
Shamir2019,
ZhangMartensGrosse2019,
ZouCaoZhouGu2019}
and the references therein.

The proof of the main result of this section, \cref{prop:minimum_MC_rate},
crucially relies
(cf.\ \cref{lem:minimum_MC_Lp})
on the complementary distribution function formula
(cf., e.g., Elbr\"achter et al.~\cite[Lemma~2.2]{ElbraechterGrohsJentzenSchwab2018arXiv})
and
the elementary estimate for the beta function
given in \cref{cor:Beta_function}.
In order to prove \cref{cor:Beta_function},
we first collect 
a few basic facts about the gamma and the beta function
in the elementary and well-known \cref{lem:Gamma_basic}
and
derive from these in \cref{prop:Gamma_function}
further elementary and essentially well-known properties of the gamma function.
In particular,
the inequalities in~\cref{eq:Gamma_ratio} in \cref{prop:Gamma_function} below
are slightly reformulated versions
of the well-known inequalities
called
\emph{Wendel's double inequality}
(cf.\ Wendel~\cite{Wendel1948})
or
\emph{Gautschi's double inequality}
(cf.\ Gautschi~\cite{Gautschi1959});
cf., e.g., Qi~\cite[Subsection~2.1 and Subsection~2.4]{Qi2010}.

\subsection{Properties of the gamma and the beta function}

\begin{lemma}
\label{lem:Gamma_basic}
Let
$ \Gamma \colon ( 0, \infty ) \to ( 0, \infty ) $
satisfy for all
$ x \in ( 0, \infty ) $
that
$ \Gamma( x ) = \int_0^{ \infty } t^{ x - 1 } e^{ - t } \ud t $
and
let
$ \bbB \colon ( 0, \infty )^2 \to ( 0, \infty ) $
satisfy for all
$ x, y \in ( 0, \infty ) $
that
$ \bbB( x, y )
=
\int_{0}^{1}
    t^{ x - 1 } ( 1 - t )^{ y - 1 }
\ud t $.
Then
\begin{enumerate}[(i)]
\item
\label{item:lem:Gamma_basic:1}
it holds for all
$ x \in ( 0, \infty ) $
that
$ \Gamma( x + 1 ) = x \, \Gamma( x ) $,
\item
\label{item:lem:Gamma_basic:2}
it holds that
$ \Gamma(1) = \Gamma(2) = 1 $,
and
\item
\label{item:lem:Gamma_basic:3}
it holds for all
$ x, y \in ( 0, \infty ) $
that
$ \bbB( x, y )
=
\frac{ \Gamma( x ) \Gamma( y ) }{ \Gamma( x + y ) } $.
\end{enumerate}
\end{lemma}

\begin{lemma}
\label{lem:unit_interval_basic}
It holds for all
$ \alpha, x \in [ 0, 1 ] $
that
$ ( 1 - x )^\alpha \leq 1 - \alpha x $.
\end{lemma}
\begin{cproof}{lem:unit_interval_basic}
Note that
the fact that
for every
$ y \in [ 0, \infty ) $
it holds that
the function
$ [ 0, \infty ) \ni z \mapsto y^z \in [ 0, \infty ) $
is a convex function
implies for all
$ \alpha, x \in [ 0, 1 ] $
that
\begin{equation}
\begin{split}
( 1 - x )^\alpha
& =
( 1 - x )^{ \alpha \cdot 1 + ( 1 - \alpha ) \cdot 0 }
\\ &
\leq
\alpha ( 1 - x )^1
+
( 1 - \alpha ) ( 1 - x )^0
\\ &
=
\alpha - \alpha x
+
1 - \alpha
=
1 - \alpha x
.
\end{split}
\end{equation}
\end{cproof}

\begin{proposition}
\label{prop:Gamma_function}
Let
$ \Gamma \colon ( 0, \infty ) \to ( 0, \infty ) $
satisfy for all
$ x \in ( 0, \infty ) $
that
$ \Gamma( x ) =
$\linebreak$
\int_0^{ \infty } t^{ x - 1 } e^{ - t } \ud t $
and
let
$ \llfloor \cdot \rrfloor \colon ( 0, \infty ) \to \N_0 $
satisfy for all
$ x \in ( 0, \infty ) $
that
$ \llfloor x \rrfloor = \max( [ 0, x ) \cap \N_0 ) $.
Then
\begin{enumerate}[(i)]
\item
\label{item:prop:Gamma_function:1}
it holds that
$ \Gamma \colon ( 0, \infty ) \to ( 0, \infty ) $
is a convex function,
\item
\label{item:prop:Gamma_function:2}
it holds for all
$ x \in ( 0, \infty ) $
that
$
\Gamma( x + 1 )
=
x \, \Gamma( x )
\leq
x^{ \llfloor x \rrfloor }
\leq
\max\{ 1, x^x \}
$,
\item
\label{item:prop:Gamma_function:3}
it holds for all
$ x \in ( 0, \infty ) $,
$ \alpha \in [ 0, 1 ] $
that
\begin{equation}
\label{eq:Gamma_ratio}
( \max\{ x + \alpha - 1, 0 \} )^{ \alpha }
\leq
\frac{x}{( x + \alpha )^{ 1 - \alpha }}
\leq
\frac{ \Gamma( x + \alpha ) }{ \Gamma( x ) }
\leq
x^\alpha
,
\end{equation}
and
\item
\label{item:prop:Gamma_function:4}
it holds for all
$ x \in ( 0, \infty ) $,
$ \alpha \in [ 0, \infty ) $
that
\begin{equation}
( \max\{ x + \min\{ \alpha - 1, 0 \}, 0 \} )^{ \alpha }
\leq
\frac{ \Gamma( x + \alpha ) }{ \Gamma( x ) }
\leq
( x + \max\{ \alpha - 1, 0 \} )^{ \alpha }
.
\end{equation}
\end{enumerate}
\end{proposition}
\begin{cproof}{prop:Gamma_function}
First,
observe that
the fact that
for every
$ t \in ( 0, \infty ) $
it holds that
the function
$ \R \ni x \mapsto t^x \in ( 0, \infty ) $
is a convex function
implies for all
$ x, y \in ( 0, \infty ) $,
$ \alpha \in [ 0, 1 ] $
that
\begin{equation}
\begin{split}
\Gamma( \alpha x + ( 1 - \alpha ) y )
& =
\int_0^{ \infty } t^{ \alpha x + ( 1 - \alpha ) y - 1 } e^{ - t } \ud t
=
\int_0^{ \infty } t^{ \alpha x + ( 1 - \alpha ) y } t^{ -1 } e^{ - t } \ud t
\\ &
\leq
\int_0^{ \infty } ( \alpha t^x + ( 1 - \alpha ) t^y ) t^{ -1 } e^{ - t } \ud t
\\ &
=
\alpha
\int_0^{ \infty } t^{ x - 1 } e^{ - t } \ud t
+
( 1 - \alpha )
\int_0^{ \infty } t^{ y - 1 } e^{ - t } \ud t
\\ &
=
\alpha
\, \Gamma( x )
+
( 1 - \alpha )
\Gamma( y )
.
\end{split}
\end{equation}
This shows~\cref{item:prop:Gamma_function:1}.

Second,
note that
\cref{item:lem:Gamma_basic:2} in \cref{lem:Gamma_basic}
and~\cref{item:prop:Gamma_function:1} establish for all
$ \alpha \in [ 0, 1 ] $
that
\begin{equation}
\Gamma( \alpha + 1 )
=
\Gamma( \alpha \cdot 2 + ( 1 - \alpha ) \cdot 1 )
\leq
\alpha
\, \Gamma( 2 )
+
( 1 - \alpha )
\Gamma( 1 )
=
\alpha + ( 1 - \alpha )
= 1
.
\end{equation}
This yields for all
$ x \in ( 0, 1 ] $
that
\begin{equation}
\label{eq:Gamma_bound_small}
\Gamma( x + 1 )
\leq
1
=
x^{ \llfloor x \rrfloor }
=
\max\{ 1, x^x \}
.
\end{equation}
Induction,
\cref{item:lem:Gamma_basic:1} in \cref{lem:Gamma_basic},
and
the fact that
$ \forall \, x \in ( 0, \infty ) \colon
x - \llfloor x \rrfloor \in ( 0, 1 ] $
hence
ensure for all
$ x \in [ 1, \infty ) $
that
\begin{equation}
\Gamma( x + 1 )
=
\biggl[
    \smallprod_{i=1}^{ \llfloor x \rrfloor }
    ( x - i + 1 )
\biggr]
\Gamma( x - \llfloor x \rrfloor + 1 )
\leq
x^{ \llfloor x \rrfloor }
\Gamma( x - \llfloor x \rrfloor + 1 )
\leq
x^{ \llfloor x \rrfloor }
\leq
x^x
=
\max\{ 1, x^x \}
.
\end{equation}
Combining this with
again~\cref{item:lem:Gamma_basic:1} in \cref{lem:Gamma_basic}
and
\cref{eq:Gamma_bound_small}
establishes~\cref{item:prop:Gamma_function:2}.

Third,
note that
H\"older's inequality
and
\cref{item:lem:Gamma_basic:1} in \cref{lem:Gamma_basic}
prove for all
$ x \in ( 0, \infty ) $,
$ \alpha \in [ 0, 1 ] $
that
\begin{equation}
\label{eq:Gamma_ratio_UB}
\begin{split}
\Gamma( x + \alpha )
& =
\int_0^{ \infty } t^{ x + \alpha - 1 } e^{ - t } \ud t
=
\int_0^{ \infty }
    t^{ \alpha x } e^{ - \alpha t }
    t^{ ( 1 - \alpha ) x - ( 1 - \alpha ) } e^{ - ( 1 - \alpha ) t }
\ud t
\\ &
=
\int_0^{ \infty }
    [ t^{ x } e^{ - t } ]^\alpha
    [ t^{ x - 1 } e^{ - t } ]^{ 1 - \alpha }
\ud t
\\ &
\leq
\biggl(
    \int_0^{ \infty } t^{ x } e^{ - t } \ud t
\biggr)^{ \!\! \alpha }
\biggl(
    \int_0^{ \infty } t^{ x - 1 } e^{ - t } \ud t
\biggr)^{ \!\! 1 - \alpha }
\\ &
=
[ \Gamma( x + 1 ) ]^{ \alpha }
[ \Gamma( x ) ]^{ 1 - \alpha }
=
x^\alpha [ \Gamma( x ) ]^{ \alpha }
[ \Gamma( x ) ]^{ 1 - \alpha }
\\ &
=
x^\alpha
\Gamma( x )
.
\end{split}
\end{equation}
This
and
again \cref{item:lem:Gamma_basic:1} in \cref{lem:Gamma_basic}
demonstrate for all
$ x \in ( 0, \infty ) $,
$ \alpha \in [ 0, 1 ] $
that
\begin{equation}
\label{eq:Gamma_ratio_LB}
x \, \Gamma( x )
=
\Gamma( x + 1 )
=
\Gamma( x + \alpha + ( 1 - \alpha ) )
\leq
( x + \alpha )^{ 1 - \alpha }
\Gamma( x + \alpha )
.
\end{equation}
Combining~\cref{eq:Gamma_ratio_UB} and~\cref{eq:Gamma_ratio_LB}
yields for all
$ x \in ( 0, \infty ) $,
$ \alpha \in [ 0, 1 ] $
that
\begin{equation}
\label{eq:sharper_bound}
\frac{x}{( x + \alpha )^{ 1 - \alpha }}
\leq
\frac{ \Gamma( x + \alpha ) }{ \Gamma( x ) }
\leq
x^\alpha
.
\end{equation}
Furthermore,
observe that~\cref{item:lem:Gamma_basic:1} in \cref{lem:Gamma_basic}
and~\cref{eq:sharper_bound}
imply for all
$ x \in ( 0, \infty ) $,
$ \alpha \in [ 0, 1 ] $
that
\begin{equation}
\frac{ \Gamma( x + \alpha ) }{ \Gamma( x + 1 ) }
=
\frac{ \Gamma( x + \alpha ) }{ x \, \Gamma( x ) }
\leq
x^{ \alpha - 1 }
.
\end{equation}
This shows for all
$ \alpha \in [ 0, 1 ] $,
$ x \in ( \alpha, \infty ) $
that
\begin{equation}
\frac{ \Gamma( x ) }{ \Gamma( x + ( 1 - \alpha ) ) }
=
\frac{ \Gamma( ( x - \alpha ) + \alpha ) }{ \Gamma( ( x - \alpha ) + 1 ) }
\leq
( x - \alpha )^{ \alpha - 1 }
=
\frac{1}{( x - \alpha )^{ 1 - \alpha }}
.
\end{equation}
This, in turn,
ensures for all
$ \alpha \in [ 0, 1 ] $,
$ x \in ( 1 - \alpha, \infty ) $
that
\begin{equation}
( x + \alpha - 1 )^{ \alpha }
=
( x - ( 1 - \alpha ) )^{ \alpha }
\leq
\frac{ \Gamma( x + \alpha ) }{ \Gamma( x ) }
.
\end{equation}
Next
note that
\cref{lem:unit_interval_basic}
proves for all
$ x \in ( 0, \infty ) $,
$ \alpha \in [ 0, 1 ] $
that
\begin{equation}
\begin{split}
( \max\{ x + \alpha - 1, 0 \} )^{ \alpha }
& =
( x + \alpha )^\alpha
\biggl(
    \frac{ \max\{ x + \alpha - 1, 0 \} }{ x + \alpha }
\biggr)^{ \!\! \alpha }
\\ &
=
( x + \alpha )^\alpha
\biggl(
    \max\biggl\{ 1 - \frac{1}{ x + \alpha }, 0 \biggr\}
\biggr)^{ \!\! \alpha }
\\ &
\leq
( x + \alpha )^\alpha
\biggl(
    1 - \frac{\alpha}{ x + \alpha }
\biggr)
=
( x + \alpha )^\alpha
\biggl(
    \frac{x}{ x + \alpha }
\biggr)
\\ &
=
\frac{x}{( x + \alpha )^{ 1 - \alpha }}
.
\end{split}
\end{equation}
This and~\cref{eq:sharper_bound}
establish~\cref{item:prop:Gamma_function:3}.

Fourth,
we show~\cref{item:prop:Gamma_function:4}.
For this
let
$ \lfloor \cdot \rfloor \colon [ 0, \infty ) \to \N_0 $
satisfy for all
$ x \in [ 0, \infty ) $
that
$ \lfloor x \rfloor = \max( [ 0, x ] \cap \N_0 ) $.
Observe that
induction,
\cref{item:lem:Gamma_basic:1} in \cref{lem:Gamma_basic},
the fact that
$ \forall \, \alpha \in [ 0, \infty ) \colon
\alpha - \lfloor \alpha \rfloor \in [ 0, 1 ) $,
and
\cref{item:prop:Gamma_function:3}
demonstrate for all
$ x \in ( 0, \infty ) $,
$ \alpha \in [ 0, \infty ) $
that
\begin{equation}
\label{eq:Gamma_ratio_UB_general}
\begin{split}
\frac{ \Gamma( x + \alpha ) }{ \Gamma( x ) }
& =
\biggl[
    \smallprod_{i=1}^{ \lfloor \alpha \rfloor }
    ( x + \alpha - i )
\biggr]
\frac{ \Gamma( x + \alpha - \lfloor \alpha \rfloor ) }{ \Gamma( x ) }
\leq
\biggl[
    \smallprod_{i=1}^{ \lfloor \alpha \rfloor }
    ( x + \alpha - i )
\biggr]
x^{ \alpha - \lfloor \alpha \rfloor }
\\ &
\leq
( x + \alpha - 1 )^{ \lfloor \alpha \rfloor }
x^{ \alpha - \lfloor \alpha \rfloor }
\\ &
\leq
( x + \max\{ \alpha - 1, 0 \} )^{ \lfloor \alpha \rfloor }
( x + \max\{ \alpha - 1, 0 \} )^{ \alpha - \lfloor \alpha \rfloor }
\\ &
=
( x + \max\{ \alpha - 1, 0 \} )^{ \alpha }
.
\end{split}
\end{equation}
Furthermore,
again
the fact that
$ \forall \, \alpha \in [ 0, \infty ) \colon
\alpha - \lfloor \alpha \rfloor \in [ 0, 1 ) $,
\cref{item:prop:Gamma_function:3},
induction,
and
\cref{item:lem:Gamma_basic:1} in \cref{lem:Gamma_basic}
imply for all
$ x \in ( 0, \infty ) $,
$ \alpha \in [ 0, \infty ) $
that
\begin{equation}
\begin{split}
\frac{ \Gamma( x + \alpha ) }{ \Gamma( x ) }
& =
\frac{ \Gamma( x + \lfloor \alpha \rfloor + \alpha - \lfloor \alpha \rfloor ) }{ \Gamma( x ) }
\\ &
\geq
( \max\{ x + \lfloor \alpha \rfloor + \alpha - \lfloor \alpha \rfloor - 1, 0 \} )^{ \alpha - \lfloor \alpha \rfloor }
\biggl[
    \frac{ \Gamma( x + \lfloor \alpha \rfloor ) }{ \Gamma( x ) }
\biggr]
\\ &
=
( \max\{ x + \alpha - 1, 0 \} )^{ \alpha - \lfloor \alpha \rfloor }
\biggl[
    \smallprod_{i=1}^{ \lfloor \alpha \rfloor }
    ( x + \lfloor \alpha \rfloor - i )
\biggr]
\frac{ \Gamma( x ) }{ \Gamma( x ) }
\\ &
\geq
( \max\{ x + \alpha - 1, 0 \} )^{ \alpha - \lfloor \alpha \rfloor }
x^{ \lfloor \alpha \rfloor }
\\ &
=
( \max\{ x + \alpha - 1, 0 \} )^{ \alpha - \lfloor \alpha \rfloor }
( \max\{ x, 0 \} )^{ \lfloor \alpha \rfloor }
\\ &
\geq
( \max\{ x + \min\{ \alpha - 1, 0 \}, 0 \} )^{ \alpha - \lfloor \alpha \rfloor }
( \max\{ x + \min\{ \alpha - 1, 0 \}, 0 \} )^{ \lfloor \alpha \rfloor }
\\ &
=
( \max\{ x + \min\{ \alpha - 1, 0 \}, 0 \} )^{ \alpha }
.
\end{split}
\end{equation}
Combining this
with~\cref{eq:Gamma_ratio_UB_general}
shows~\cref{item:prop:Gamma_function:4}.
\end{cproof}

\begin{corollary}
\label{cor:Beta_function}
Let
$ \bbB \colon ( 0, \infty )^2 \to ( 0, \infty ) $
satisfy for all
$ x, y \in ( 0, \infty ) $
that
$ \bbB( x, y )
=
\int_{0}^{1}
    t^{ x - 1 } ( 1 - t )^{ y - 1 }
\ud t $
and
let
$ \Gamma \colon ( 0, \infty ) \to ( 0, \infty ) $
satisfy for all
$ x \in ( 0, \infty ) $
that
$ \Gamma( x ) = \int_0^{ \infty } t^{ x - 1 } e^{ - t } \ud t $.
Then it holds for all
$ x, y \in ( 0, \infty ) $
with
$ x + y > 1 $
that
\begin{equation}
\label{eq:cor:Beta_function}
\frac{\Gamma( x )}{ ( y + \max\{ x - 1, 0 \} )^{ x } }
\leq
\bbB( x, y )
\leq
\frac{\Gamma( x )}{ ( y + \min\{ x - 1, 0 \} )^{ x } }
\leq
\frac{ \max\{ 1, x^{ x } \} }{ x ( y + \min\{ x - 1, 0 \} )^{ x } }
.
\end{equation}
\end{corollary}
\begin{cproof}{cor:Beta_function}
Note that
\cref{item:lem:Gamma_basic:3} in \cref{lem:Gamma_basic}
ensures for all
$ x, y \in ( 0, \infty ) $
that
\begin{equation}
\label{eq:Beta_Gamma_relation}
\bbB( x, y )
=
\frac{ \Gamma( x ) \Gamma( y ) }{ \Gamma( y + x ) }
.
\end{equation}
In addition,
observe that it holds for all
$ x, y \in ( 0, \infty ) $
with
$ x + y > 1 $
that
$ y + \min\{ x - 1, 0 \} > 0 $.
This and
\cref{item:prop:Gamma_function:4} in \cref{prop:Gamma_function}
demonstrate for all
$ x, y \in ( 0, \infty ) $
with
$ x + y > 1 $
that%
\begin{equation}
0 <
( y + \min\{ x - 1, 0 \} )^{ x }
\leq
\frac{ \Gamma( y + x ) }{ \Gamma( y ) }
\leq
( y + \max\{ x - 1, 0 \} )^{ x }
.
\end{equation}
Combining this
with~\cref{eq:Beta_Gamma_relation}
and
\cref{item:prop:Gamma_function:2} in \cref{prop:Gamma_function}
shows for all
$ x, y \in ( 0, \infty ) $
with
$ x + y > 1 $
that
\begin{equation}
\frac{\Gamma( x )}{ ( y + \max\{ x - 1, 0 \} )^{ x } }
\leq
\bbB( x, y )
\leq
\frac{\Gamma( x )}{ ( y + \min\{ x - 1, 0 \} )^{ x } }
\leq
\frac{ \max\{ 1, x^{ x } \} }{ x ( y + \min\{ x - 1, 0 \} )^{ x } }
.
\end{equation}
\end{cproof}

\subsection{Strong convergence rates for the optimisation error}

\begin{lemma}
\label{lem:minimum_MC_Lp}
Let
$ K \in \N $,
$ p, L \in ( 0, \infty ) $,
let
$ ( E, \delta ) $
be a metric space,
let
$ ( \Omega, \cF, \P ) $
be a probability space,
let
$ \cR \colon E \times \Omega \to \R $
be a
$ ( \cB( E ) \otimes \cF ) $/$ \cB( \R ) $-measurable function,
assume for all
$ x, y \in E $,
$ \omega \in \Omega $
that
$ \lvert \cR( x, \omega ) - \cR( y, \omega ) \rvert
\leq L \delta( x, y ) $,
and let
$ X_k \colon \Omega \to E $, $ k \in \{ 1, 2, \ldots, K \} $,
be i.i.d.\ random variables.
Then it holds for all
$ x \in E $
that
\begin{equation}
\label{eq:minimum_MC_Lp}
\E\bigl[
    \min\nolimits_{ k \in \{ 1, 2, \ldots, K \} } \lvert \cR( X_k ) - \cR( x ) \rvert^p
\bigr]
\leq
L^p
\int_{0}^{\infty}
    [ \P( \delta( X_1, x ) > \varepsilon^{ \nicefrac{1}{p} } ) ]^K
\ud \varepsilon
.
\end{equation}
\end{lemma}
\begin{cproof}{lem:minimum_MC_Lp}
Throughout this proof let
$ x \in E $
and
let
$ Y \colon \Omega \to [ 0, \infty ) $
be the function which satisfies for all
$ \omega \in \Omega $
that
$ Y( \omega ) = \min\nolimits_{ k \in \{ 1, 2, \ldots, K \} } [ \delta( X_k( \omega ), x ) ]^p $.
Observe that
the fact that
$ Y $ is a random variable,
the assumption that
$ \forall \, x, y \in E,
\, \omega \in \Omega \colon
\lvert \cR( x, \omega ) - \cR( y, \omega ) \rvert
\leq L \delta( x, y ) $,
and
the
complementary distribution function formula
(see, e.g., Elbr\"achter et al.~\cite[Lemma~2.2]{ElbraechterGrohsJentzenSchwab2018arXiv})
demonstrate that
\begin{equation}
\label{eq:expectation_min}
\begin{split}
&
\E\bigl[
    \min\nolimits_{ k \in \{ 1, 2, \ldots, K \} } \lvert \cR( X_k ) - \cR( x ) \rvert^p
\bigr]
\leq
L^p
\, \E\bigl[
    \min\nolimits_{ k \in \{ 1, 2, \ldots, K \} } [ \delta( X_k, x ) ]^p
\bigr]
\\ &
=
L^p
\, \E[ Y ]
=
L^p
\int_{0}^{\infty}
y
\, \P_Y( \dd y )
=
L^p
\int_{0}^{\infty}
\P_Y( ( \varepsilon, \infty ) )
\ud \varepsilon
\\ &
=
L^p
\int_{0}^{\infty}
\P( Y > \varepsilon )
\ud \varepsilon
=
L^p
\int_{0}^{\infty}
    \P\bigl( \min\nolimits_{ k \in \{ 1, 2, \ldots, K \} } [ \delta( X_k, x ) ]^p > \varepsilon \bigr)
\ud \varepsilon
.
\end{split}
\end{equation}
Moreover,
the assumption that
$ \Theta_k $, $ k \in \{ 1, 2, \ldots, K \} $,
are i.i.d.\ random variables
shows for all
$ \varepsilon \in ( 0, \infty ) $
that
\begin{equation}
\label{eq:probability_min}
\begin{split}
&
\P\bigl( \min\nolimits_{ k \in \{ 1, 2, \ldots, K \} } [ \delta( X_k, x ) ]^p > \varepsilon \bigr)
=
\P\bigl(
    \forall \, k \in \{ 1, 2, \ldots, K \} \colon
    [ \delta( X_k, x ) ]^p > \varepsilon
\bigr)
\\ &
=
\smallprod_{ k = 1 }^K
    \P( [ \delta( X_k, x ) ]^p > \varepsilon )
=
[ \P( [ \delta( X_1, x ) ]^p > \varepsilon ) ]^K
=
[ \P( \delta( X_1, x ) > \varepsilon^{ \nicefrac{1}{p} } ) ]^K
.
\end{split}
\end{equation}
Combining~\cref{eq:expectation_min} with~\cref{eq:probability_min}
proves~\cref{eq:minimum_MC_Lp}.
\end{cproof}

\begin{proposition}
\label{prop:minimum_MC_rate}
Let
$ \bfd, K \in \N $,
$ L, \alpha \in \R $,
$ \beta \in ( \alpha, \infty ) $,
let
$ ( \Omega, \cF, \P ) $
be a probability space,
let
$ \cR \colon [ \alpha, \beta ]^\bfd \times \Omega \to \R $
be a random field,
assume for all
$ \theta, \vartheta \in [ \alpha, \beta ]^\bfd $,
$ \omega \in \Omega $
that
$ \lvert \cR( \theta, \omega ) - \cR( \vartheta, \omega ) \rvert
\leq L \lVert \theta - \vartheta \rVert_{ \infty } $,
let
$ \Theta_k \colon \Omega \to [ \alpha, \beta ]^\bfd $, $ k \in \{ 1, 2, \ldots, K \} $,
be i.i.d.\ random variables,
and assume that
$ \Theta_1 $ is continuous uniformly distributed on $ [ \alpha, \beta ]^\bfd $
(cf.~\cref{def:p-norm}).
Then
\begin{enumerate}[(i)]
\item
\label{item:prop:minimum_MC_rate:1}
it holds that
$ \cR $
is a
$ ( \cB( [ \alpha, \beta ]^\bfd ) \otimes \cF ) $/$ \cB( \R ) $-measurable function
and
\item
\label{item:prop:minimum_MC_rate:2}
it holds for all
$ \theta \in [ \alpha, \beta ]^\bfd $,
$ p \in ( 0, \infty ) $
that
\begin{equation}
\begin{split}
&
\bigl(
\E\bigl[
    \min\nolimits_{ k \in \{ 1, 2, \ldots, K \} } \lvert \cR( \Theta_k ) - \cR( \theta ) \rvert^p
\bigr]
\bigr)^{ \nicefrac{1}{p} }
\\ &
\leq
\frac{ L ( \beta - \alpha ) \max\{ 1, ( \nicefrac{p}{\bfd} )^{ \nicefrac{1}{\bfd} } \} }{ K^{ \nicefrac{1}{\bfd} } }
\leq
\frac{ L ( \beta - \alpha ) \max\{ 1, p \} }{ K^{ \nicefrac{1}{\bfd} } }
.
\end{split}
\end{equation}
\end{enumerate}
\end{proposition}
\begin{cproof2}{prop:minimum_MC_rate}
Throughout this proof
assume w.l.o.g.\ that
$ L > 0 $,
let
$ \delta \colon
\allowbreak
( [ \alpha, \beta ]^\bfd ) \times ( [ \alpha, \beta ]^\bfd ) \to [ 0, \infty ) $
satisfy for all
$ \theta, \vartheta \in [ \alpha, \beta ]^\bfd $
that
$ \delta( \theta, \vartheta ) = \lVert \theta - \vartheta \rVert_\infty $,
let
$ \bbB \colon ( 0, \infty )^2 \to ( 0, \infty ) $
satisfy for all
$ x, y \in ( 0, \infty ) $
that
$ \bbB( x, y )
=
\int_{0}^{1}
    t^{ x - 1 } ( 1 - t )^{ y - 1 }
\ud t $,
and
let
$ \Theta_{ 1, 1 }, \Theta_{ 1, 2 }, \ldots, \Theta_{ 1, \bfd } \colon \Omega \to [ \alpha, \beta ] $
satisfy
$ \Theta_1 = ( \Theta_{ 1, 1 }, \Theta_{ 1, 2 }, \ldots, \Theta_{ 1, \bfd } ) $.
First of all,
note that
the assumption that
$ \forall \, \theta, \vartheta \in [ \alpha, \beta ]^\bfd,
\, \omega \in \Omega \colon
\lvert \cR( \theta, \omega ) - \cR( \vartheta, \omega ) \rvert
\leq L \lVert \theta - \vartheta \rVert_{ \infty } $
ensures for all
$ \omega \in \Omega $
that the function
$ [ \alpha, \beta ]^\bfd \ni \theta
\mapsto \cR( \theta, \omega ) \in \R $
is continuous.
Combining this with
the fact that
$ ( [ \alpha, \beta ]^\bfd, \delta ) $
is a separable metric space,
the fact that
for every
$ \theta \in [ \alpha, \beta ]^\bfd $
it holds that the function
$ \Omega \ni \omega
\mapsto \cR( \theta, \omega ) \in \R $
is $ \cF $/$ \cB( \R ) $-measurable,
and, e.g.,
Aliprantis \& Border~\cite[Lemma~4.51]{AliprantisBorder2006}
(see also, e.g.,
Beck et al.~\cite[Lemma~2.4]{BeckBeckerGrohsJaafariJentzen2018arXiv})
proves~\cref{item:prop:minimum_MC_rate:1}.
Next observe that it holds for all
$ \theta \in [ \alpha, \beta ] $,
$ \varepsilon \in [ 0, \infty ) $
that
\begin{equation}
\begin{split}
& \min\{ \theta + \varepsilon, \beta \} - \max\{ \theta - \varepsilon, \alpha \}
=
\min\{ \theta + \varepsilon, \beta \} + \min\{ \varepsilon - \theta, -\alpha \}
\\ &
=
\min\bigl\{
    \theta + \varepsilon + \min\{ \varepsilon - \theta, -\alpha \},
    \beta + \min\{ \varepsilon - \theta, -\alpha \}
\bigr\}
\\ &
=
\min\bigl\{
    \min\{ 2 \varepsilon, \theta - \alpha + \varepsilon \},
    \min\{ \beta - \theta + \varepsilon, \beta - \alpha \}
\bigr\}
\\ &
\geq
\min\bigl\{
    \min\{ 2 \varepsilon, \alpha - \alpha + \varepsilon \},
    \min\{ \beta - \beta + \varepsilon, \beta - \alpha \}
\bigr\}
\\ &
=
\min\{
    2 \varepsilon, \varepsilon,
    \varepsilon, \beta - \alpha
\}
=
\min\{ \varepsilon, \beta - \alpha \}
.
\end{split}
\end{equation}
The assumption that
$ \Theta_1 $ is continuous uniformly distributed on $ [ \alpha, \beta ]^\bfd $
hence
shows for all
$ \theta = ( \theta_1, \theta_2, \ldots, \theta_\bfd ) \in [ \alpha, \beta ]^\bfd $,
$ \varepsilon \in [ 0, \infty ) $
that
\begin{equation}
\begin{split}
& \P( \lVert \Theta_1 - \theta \rVert_{ \infty } \leq \varepsilon )
=
\P\bigl(
    \max\nolimits_{ i \in \{ 1, 2, \ldots, \bfd \} } \lvert \Theta_{ 1, i } - \theta_i \rvert \leq \varepsilon
\bigr)
\\ &
=
\P\bigl(
    \forall \, i \in \{ 1, 2, \ldots, \bfd \} \colon
    -\varepsilon \leq \Theta_{ 1, i } - \theta_i \leq \varepsilon
\bigr)
\\ &
=
\P\bigl(
    \forall \, i \in \{ 1, 2, \ldots, \bfd \} \colon
    \theta_i - \varepsilon \leq \Theta_{ 1, i } \leq \theta_i + \varepsilon
\bigr)
\\ &
=
\P\bigl(
    \forall \, i \in \{ 1, 2, \ldots, \bfd \} \colon
    \max\{ \theta_i - \varepsilon, \alpha \} \leq \Theta_{ 1, i } \leq \min\{ \theta_i + \varepsilon, \beta \}
\bigr)
\\ &
=
\P\bigl(
    \Theta_1 \in \bigl[ \times_{ i = 1 }^\bfd [ \max\{ \theta_i - \varepsilon, \alpha \}, \min\{ \theta_i + \varepsilon, \beta \} ] \bigr]
\bigr)
\\ &
=
\tfrac{1}{ ( \beta - \alpha )^\bfd }
\smallprod_{ i = 1 }^\bfd
    ( \min\{ \theta_i + \varepsilon, \beta \} - \max\{ \theta_i - \varepsilon, \alpha \} )
\\ &
\geq
\tfrac{1}{ ( \beta - \alpha )^\bfd }
[ \min\{ \varepsilon, \beta - \alpha \} ]^\bfd
=
\min\Bigl\{
    1, \tfrac{ \varepsilon^\bfd }{ ( \beta - \alpha )^\bfd }
\Bigr\}
.
\end{split}
\end{equation}
Therefore, we obtain
for all
$ \theta \in [ \alpha, \beta ]^\bfd $,
$ p \in ( 0, \infty ) $,
$ \varepsilon \in [ 0, \infty ) $
that
\begin{equation}
\begin{split}
&
\P( \lVert \Theta_1 - \theta \rVert_{ \infty } > \varepsilon^{ \nicefrac{1}{p} } )
=
1 - \P( \lVert \Theta_1 - \theta \rVert_{ \infty } \leq \varepsilon^{ \nicefrac{1}{p} } )
\\ &
\leq
1
-
\min\Bigl\{
    1, \tfrac{ \varepsilon^{ \nicefrac{\bfd}{p} } }{ ( \beta - \alpha )^\bfd }
\Bigr\}
=
\max\Bigl\{
    0, 1 - \tfrac{ \varepsilon^{ \nicefrac{\bfd}{p} } }{ ( \beta - \alpha )^\bfd }
\Bigr\}
.
\end{split}
\end{equation}
This,
\cref{item:prop:minimum_MC_rate:1},
the assumption that
$ \forall \, \theta, \vartheta \in [ \alpha, \beta ]^\bfd,
\, \omega \in \Omega \colon
\lvert \cR( \theta, \omega ) - \cR( \vartheta, \omega ) \rvert
\leq L \lVert \theta - \vartheta \rVert_{ \infty } $,
the assumption that
$ \Theta_k $, $ k \in \{ 1, 2, \ldots, K \} $,
are i.i.d.\ random variables,
and
\cref{lem:minimum_MC_Lp}
(with
$ ( E, \delta ) \leftarrow ( [ \alpha, \beta ]^\bfd, \delta ) $,
$ ( X_k )_{ k \in \{ 1, 2, \ldots, K \} }
\leftarrow
( \Theta_k )_{ k \in \{ 1, 2, \ldots, K \} } $
in the notation of \cref{lem:minimum_MC_Lp})
establish for all
$ \theta \in [ \alpha, \beta ]^\bfd $,
$ p \in ( 0, \infty ) $
that
\begin{equation}
\begin{split}
&
\E\bigl[
    \min\nolimits_{ k \in \{ 1, 2, \ldots, K \} } \lvert \cR( \Theta_k ) - \cR( \theta ) \rvert^p
\bigr]
\leq
L^p
\int_{0}^{\infty}
    [ \P( \lVert \Theta_1 - \theta \rVert_{ \infty } > \varepsilon^{ \nicefrac{1}{p} } ) ]^K
\ud \varepsilon
\\ &
\leq
L^p
\int_{0}^{\infty}
    \Bigl[ \max\Bigl\{
        0, 1 - \tfrac{ \varepsilon^{ \nicefrac{\bfd}{p} } }{ ( \beta - \alpha )^\bfd }
    \Bigr\} \Bigr]^K
\ud \varepsilon
=
L^p
\int_{0}^{ ( \beta - \alpha )^p }
    \Bigl(
        1 - \tfrac{ \varepsilon^{ \nicefrac{\bfd}{p} } }{ ( \beta - \alpha )^\bfd }
    \Bigr)^{ \! K }
\ud \varepsilon
\\ &
=
\tfrac{p}{\bfd}
L^p
( \beta - \alpha )^p
\int_{0}^{ 1 }
    t^{ \nicefrac{p}{\bfd} - 1 }
    ( 1 - t )^K
\ud t
=
\tfrac{p}{\bfd}
L^p
( \beta - \alpha )^p
\int_{0}^{ 1 }
    t^{ \nicefrac{p}{\bfd} - 1 }
    ( 1 - t )^{ K + 1 - 1 }
\ud t
\\ &
=
\tfrac{p}{\bfd}
L^p
( \beta - \alpha )^p
\, \bbB( \nicefrac{p}{\bfd}, K + 1 )
.
\end{split}
\end{equation}
\cref{cor:Beta_function}
(with
$ x \leftarrow \nicefrac{p}{\bfd} $,
$ y \leftarrow K + 1 $
for
$ p \in ( 0, \infty ) $
in the notation of \cref{eq:cor:Beta_function} in \cref{cor:Beta_function})
hence demonstrates for all
$ \theta \in [ \alpha, \beta ]^\bfd $,
$ p \in ( 0, \infty ) $
that
\begin{equation}
\begin{split}
&
\E\bigl[
    \min\nolimits_{ k \in \{ 1, 2, \ldots, K \} } \lvert \cR( \Theta_k ) - \cR( \theta ) \rvert^p
\bigr]
\\ &
\leq
\frac{ \tfrac{p}{\bfd} L^p ( \beta - \alpha )^p \max\{ 1, ( \nicefrac{p}{\bfd} )^{ \nicefrac{p}{\bfd} } \}
}{
\tfrac{p}{\bfd} ( K + 1 + \min\{ \nicefrac{p}{\bfd} - 1, 0 \} )^{ \nicefrac{p}{\bfd} } }
\leq
\frac{ L^p ( \beta - \alpha )^p \max\{ 1, ( \nicefrac{p}{\bfd} )^{ \nicefrac{p}{\bfd} } \}
}{
K^{ \nicefrac{p}{\bfd} } }
.
\end{split}
\end{equation}
This implies for all
$ \theta \in [ \alpha, \beta ]^\bfd $,
$ p \in ( 0, \infty ) $
that
\begin{equation}
\begin{split}
&
\bigl(
\E\bigl[
    \min\nolimits_{ k \in \{ 1, 2, \ldots, K \} } \lvert \cR( \Theta_k ) - \cR( \theta ) \rvert^p
\bigr]
\bigr)^{ \nicefrac{1}{p} }
\\ &
\leq
\frac{ L ( \beta - \alpha ) \max\{ 1, ( \nicefrac{p}{\bfd} )^{ \nicefrac{1}{\bfd} } \} }{ K^{ \nicefrac{1}{\bfd} } }
\leq
\frac{ L ( \beta - \alpha ) \max\{ 1, p \} }{ K^{ \nicefrac{1}{\bfd} } }
.
\end{split}
\end{equation}
This shows~\cref{item:prop:minimum_MC_rate:2}
and thus completes
\end{cproof2}

\begin{lemma}
\label{lem:Lipschitz_risk}
Let
$ d, \bfd, \bfL, M \in \N $,
$ B, b \in [ 1, \infty ) $,
$ u \in \R $,
$ v \in ( u, \infty ) $,
$ \bfl = ( \bfl_0, \bfl_1, \ldots, \bfl_\bfL ) \in \N^{ \bfL + 1 } $,
$ D \subseteq [ -b, b ]^d $,
assume
$ \bfl_0 = d $,
$ \bfl_\bfL = 1 $,
and
$ \bfd \geq \sum_{i=1}^{\bfL} \bfl_i( \bfl_{ i - 1 } + 1 ) $,
let
$ \Omega $
be a set,
let
$ X_j \colon \Omega \to D $,
$ j \in \{ 1, 2, \ldots, M \} $,
and
$ Y_j \colon \Omega \to [ u, v ] $,
$ j \in \{ 1, 2, \ldots, M \} $,
be functions,
and
let
$ \cR \colon [ -B, B ]^\bfd \times \Omega \to [ 0, \infty ) $
satisfy for all
$ \theta \in [ -B, B ]^\bfd $,
$ \omega \in \Omega $
that
\begin{equation}
\cR( \theta, \omega )
=
\frac{1}{M}
\biggl[
\smallsum_{j=1}^M
    \lvert \clippedNN{\theta}{\bfl}{u}{v}( X_j( \omega ) ) - Y_j( \omega ) \rvert^2
\biggr]
\end{equation}
(cf.~\cref{def:clipped_NN}).
Then it holds for all
$ \theta, \vartheta \in [ -B, B ]^\bfd $,
$ \omega \in \Omega $
that
\begin{equation}
\lvert \cR( \theta, \omega ) - \cR( \vartheta, \omega ) \rvert
\leq
2 ( v - u ) b
\bfL
( \lVert \bfl \rVert_{ \infty } + 1 )^\bfL
B^{ \bfL - 1 }
\lVert \theta - \vartheta \rVert_{ \infty }
\end{equation}
(cf.~\cref{def:p-norm}).
\end{lemma}
\begin{cproof}{lem:Lipschitz_risk}
Observe that
the fact that
$ \forall \, x_1, x_2, y \in \R \colon
( x_1 - y )^2 - ( x_2 - y)^2
= ( x_1 - x_2 )( ( x_1 - y ) + ( x_2 - y ) ) $,
the fact that
$ \forall \, \theta \in \R^\bfd,
\, x \in \R^d \colon
\clippedNN{\theta}{\bfl}{u}{v}( x ) \in [ u, v ] $,
and the assumption that
$ \forall \, j \in \{ 1, 2, \ldots, M \},
\, \omega \in \Omega \colon
Y_j( \omega ) \in [ u, v ] $
prove for all
$ \theta, \vartheta \in [ -B, B ]^\bfd $,
$ \omega \in \Omega $
that
\begin{equation}
\label{eq:difference_cR}
\begin{split}
& \lvert \cR( \theta, \omega ) - \cR( \vartheta, \omega ) \rvert
\\ &
=
\frac{1}{M}
\biggl\lvert
\biggl[
\smallsum_{j=1}^M
    \lvert \clippedNN{\theta}{\bfl}{u}{v}( X_j( \omega ) ) - Y_j( \omega ) \rvert^2
\biggr]
-
\biggl[
\smallsum_{j=1}^M
    \lvert \clippedNN{\vartheta}{\bfl}{u}{v}( X_j( \omega ) ) - Y_j( \omega ) \rvert^2
\biggr]
\biggr\rvert
\\ &
\leq
\frac{1}{M}
\biggl[
\smallsum_{j=1}^M
    \bigl\lvert
    [ \clippedNN{\theta}{\bfl}{u}{v}( X_j( \omega ) ) - Y_j( \omega ) ]^2
    -
    [ \clippedNN{\vartheta}{\bfl}{u}{v}( X_j( \omega ) ) - Y_j( \omega ) ]^2
    \bigr\rvert
\biggr]
\\ &
=
\frac{1}{M}
\biggl[
\smallsum_{j=1}^M
    \bigl(
    \bigl\lvert
        \clippedNN{\theta}{\bfl}{u}{v}( X_j( \omega ) ) - \clippedNN{\vartheta}{\bfl}{u}{v}( X_j( \omega ) )
    \bigr\rvert
\\ &
\hphantom{
\; =
\frac{1}{M}
\biggl[
\smallsum_{j=1}^M
    \bigl(
}
    \cdot
    \bigl\lvert
    [ \clippedNN{\theta}{\bfl}{u}{v}( X_j( \omega ) ) - Y_j( \omega ) ]
    +
    [ \clippedNN{\vartheta}{\bfl}{u}{v}( X_j( \omega ) ) - Y_j( \omega ) ]
    \bigr\rvert
    \bigr)
\biggr]
\\ &
\leq
\frac{2}{M}
\biggl[
\smallsum_{j=1}^M
    \bigl(
    \bigl[
    \sup\nolimits_{ x \in D }
    \lvert
        \clippedNN{\theta}{\bfl}{u}{v}( x ) - \clippedNN{\vartheta}{\bfl}{u}{v}( x )
    \rvert
    \bigr]
    \bigl[
    \sup\nolimits_{ y_1, y_2 \in [ u, v ] }
        \lvert y_1 - y_2 \rvert
    \bigr]
    \bigr)
\biggr]
\\ &
=
2 ( v - u )
\bigl[
\sup\nolimits_{ x \in D }
    \lvert
        \clippedNN{\theta}{\bfl}{u}{v}( x ) - \clippedNN{\vartheta}{\bfl}{u}{v}( x )
    \rvert
\bigr]
.
\end{split}
\end{equation}
In addition,
combining
the assumptions that
$ D \subseteq [ -b, b ]^d $,
$ \bfd \geq \sum_{i=1}^{\bfL} \bfl_i( \bfl_{ i - 1 } + 1 ) $,
$ \bfl_0 = d $,
$ \bfl_\bfL = 1 $,
$ b \geq 1 $,
and
$ B \geq 1 $
with
Beck, Jentzen, \& Kuckuck~\cite[Corollary~2.37]{BeckJentzenKuckuck2019arXiv}
(with
$ a \leftarrow -b $,
$ b \leftarrow b $,
$ u \leftarrow u $,
$ v \leftarrow v $,
$ d \leftarrow \bfd $,
$ L \leftarrow \bfL $,
$ l \leftarrow \bfl $
in the notation of~\cite[Corollary~2.37]{BeckJentzenKuckuck2019arXiv})
shows for all
$ \theta, \vartheta \in [ -B, B ]^\bfd $
that
\begin{equation}
\begin{split}
&
\sup\nolimits_{ x \in D }
    \lvert
        \clippedNN{\theta}{\bfl}{u}{v}( x ) - \clippedNN{\vartheta}{\bfl}{u}{v}( x )
    \rvert
\leq
\sup\nolimits_{ x \in [ -b, b ]^d }
    \lvert
        \clippedNN{\theta}{\bfl}{u}{v}( x ) - \clippedNN{\vartheta}{\bfl}{u}{v}( x )
    \rvert
\\ &
\leq
\bfL
\max\{ 1, b \}
( \lVert \bfl \rVert_{ \infty } + 1 )^\bfL
( \max\{ 1, \lVert \theta \rVert_{ \infty }, \lVert \vartheta \rVert_{ \infty } \} )^{ \bfL - 1 }
\lVert \theta - \vartheta \rVert_{ \infty }
\\ &
\leq
b
\bfL
( \lVert \bfl \rVert_{ \infty } + 1 )^\bfL
B^{ \bfL - 1 }
\lVert \theta - \vartheta \rVert_{ \infty }
.
\end{split}
\end{equation}
This and~\cref{eq:difference_cR}
imply for all
$ \theta, \vartheta \in [ -B, B ]^\bfd $,
$ \omega \in \Omega $
that
\begin{equation}
\lvert \cR( \theta, \omega ) - \cR( \vartheta, \omega ) \rvert
\leq
2 ( v - u ) b
\bfL
( \lVert \bfl \rVert_{ \infty } + 1 )^\bfL
B^{ \bfL - 1 }
\lVert \theta - \vartheta \rVert_{ \infty }
.
\end{equation}
\end{cproof}

\begin{corollary}
\label{cor:minimum_MC_rate}
Let
$ d, \bfd, \fd, \bfL, M, K \in \N $,
$ B, b \in [ 1, \infty ) $,
$ u \in \R $,
$ v \in ( u, \infty ) $,
$ \bfl = ( \bfl_0, \bfl_1, \ldots, \bfl_\bfL ) \in \N^{ \bfL + 1 } $,
$ D \subseteq [ -b, b ]^d $,
assume
$ \bfl_0 = d $,
$ \bfl_\bfL = 1 $,
and
$ \bfd \geq \fd = \sum_{i=1}^{\bfL} \bfl_i( \bfl_{ i - 1 } + 1 ) $,
let
$ ( \Omega, \cF, \P ) $
be a probability space,
let
$ \Theta_k \colon \Omega \to [ -B, B ]^\bfd $, $ k \in \{ 1, 2, \ldots, K \} $,
be i.i.d.\ random variables,
assume that
$ \Theta_1 $ is continuous uniformly distributed on $ [ -B, B ]^\bfd $,
let
$ X_j \colon \Omega \to D $,
$ j \in \{ 1, 2, \ldots, M \} $,
and
$ Y_j \colon \Omega \to [ u, v ] $,
$ j \in \{ 1, 2, \ldots, M \} $,
be random variables,
and
let
$ \cR \colon [ -B, B ]^\bfd \times \Omega \to [ 0, \infty ) $
satisfy for all
$ \theta \in [ -B, B ]^\bfd $,
$ \omega \in \Omega $
that
\begin{equation}
\cR( \theta, \omega )
=
\frac{1}{M}
\biggl[
\smallsum_{j=1}^M
    \lvert \clippedNN{\theta}{\bfl}{u}{v}( X_j( \omega ) ) - Y_j( \omega ) \rvert^2
\biggr]
\end{equation}
(cf.~\cref{def:clipped_NN}).
Then
\begin{enumerate}[(i)]
\item
\label{item:cor:minimum_MC_rate:1}
it holds that
$ \cR $
is a
$ ( \cB( [ -B, B ]^\bfd ) \otimes \cF ) $/$ \cB( [ 0, \infty ) ) $-measurable function
and
\item
\label{item:cor:minimum_MC_rate:2}
it holds
for all
$ \theta \in [ -B, B ]^\bfd $,
$ p \in ( 0, \infty ) $
that
\begin{align*}
& \yesnumber
\bigl(
\E\bigl[
    \min\nolimits_{ k \in \{ 1, 2, \ldots, K \} } \lvert \cR( \Theta_k ) - \cR( \theta ) \rvert^p
\bigr]
\bigr)^{ \nicefrac{1}{p} }
\\ &
\leq
\frac{
4 ( v - u ) b
\bfL
( \lVert \bfl \rVert_{ \infty } + 1 )^\bfL
B^{ \bfL }
\sqrt{ \max\{ 1, \nicefrac{p}{\fd} \} }
}{
K^{ \nicefrac{1}{\fd} }
}
\leq
\frac{
4 ( v - u ) b
\bfL
( \lVert \bfl \rVert_{ \infty } + 1 )^\bfL
B^{ \bfL }
\max\{ 1, p \}
}{
K^{ [ \bfL^{-1} ( \lVert \bfl \rVert_{ \infty } + 1 )^{-2} ] }
}
\end{align*}
(cf.~\cref{def:p-norm}).
\end{enumerate}
\end{corollary}
\begin{cproof}{cor:minimum_MC_rate}
Throughout this proof
let
$ L \in \R $
be given by
$ L =
2 ( v - u ) b
\bfL
( \lVert \bfl \rVert_{ \infty } \allowbreak + 1 )^\bfL
B^{ \bfL - 1 } $,
let
$ P \colon [ -B, B ]^\bfd \to [ -B, B ]^\fd $
satisfy for all
$ \theta = ( \theta_1, \theta_2, \ldots, \theta_\bfd ) \in [ -B, B ]^\bfd $
that
$ P( \theta ) = ( \theta_1, \theta_2, \ldots, \theta_\fd ) $,
and
let
$ R \colon [ -B, B ]^\fd \times \Omega \to \R $
satisfy for all
$ \theta \in [ -B, B ]^\fd $,
$ \omega \in \Omega $
that
\begin{equation}
R( \theta, \omega )
=
\frac{1}{M}
\biggl[
\smallsum_{j=1}^M
    \lvert \clippedNN{\theta}{\bfl}{u}{v}( X_j( \omega ) ) - Y_j( \omega ) \rvert^2
\biggr]
.
\end{equation}
Note that
the fact that
$ \forall \, \theta \in [ -B, B ]^\bfd \colon
\clippedNN{\theta}{\bfl}{u}{v}
=
\clippedNN{ \smash{ P( \theta ) } }{\bfl}{u}{v} $
implies for all
$ \theta \in [ -B, B ]^\bfd $,
$ \omega \in \Omega $
that
\begin{equation}
\label{eq:clipped_parameters}
\begin{split}
\cR( \theta, \omega )
& =
\frac{1}{M}
\biggl[
\smallsum_{j=1}^M
    \lvert \clippedNN{\theta}{\bfl}{u}{v}( X_j( \omega ) ) - Y_j( \omega ) \rvert^2
\biggr]
\\ &
=
\frac{1}{M}
\biggl[
\smallsum_{j=1}^M
    \lvert \clippedNN{ P( \theta ) }{\bfl}{u}{v}( X_j( \omega ) ) - Y_j( \omega ) \rvert^2
\biggr]
=
R( P( \theta ), \omega )
.
\end{split}
\end{equation}
Furthermore,
\cref{lem:Lipschitz_risk}
(with
$ \bfd \leftarrow \fd $,
$ \cR \leftarrow
( [ -B, B ]^\fd \times \Omega \ni ( \theta, \omega )
\mapsto R( \theta, \omega ) \in [ 0, \infty ) ) $
in the notation of \cref{lem:Lipschitz_risk})
demonstrates for all
$ \theta, \vartheta \in [ -B, B ]^\fd $,
$ \omega \in \Omega $
that
\begin{equation}
\label{eq:Lipschitz_R}
\lvert R( \theta, \omega ) - R( \vartheta, \omega ) \rvert
\leq
2 ( v - u ) b
\bfL
( \lVert \bfl \rVert_{ \infty } + 1 )^\bfL
B^{ \bfL - 1 }
\lVert \theta - \vartheta \rVert_{ \infty }
=
L
\lVert \theta - \vartheta \rVert_{ \infty }
.
\end{equation}
Moreover,
observe that
the assumption that
$ X_j $,
$ j \in \{ 1, 2, \ldots, M \} $,
and
$ Y_j $,
$ j \in \{ 1, 2, \ldots, \allowbreak M \} $,
are random variables
ensures that
$ R \colon [ -B, B ]^\fd \times \Omega \to \R $
is a random field.
This,
\cref{eq:Lipschitz_R},
the fact that
$ P \circ \Theta_k \colon \Omega \to [ -B, B ]^\fd $, $ k \in \{ 1, 2, \ldots, K \} $,
are i.i.d.\ random variables,
the fact that
$ P \circ \Theta_1 $ is continuous uniformly distributed on $ [ -B, B ]^\fd $,
and
\cref{prop:minimum_MC_rate}
(with
$ \bfd \leftarrow \fd $,
$ \alpha \leftarrow -B $,
$ \beta \leftarrow B $,
$ \cR \leftarrow R $,
$ ( \Theta_k )_{ k \in \{ 1, 2, \ldots, K \} }
\leftarrow
( P \circ \Theta_k )_{ k \in \{ 1, 2, \ldots, K \} } $
in the notation of \cref{prop:minimum_MC_rate})
prove for all
$ \theta \in [ -B, B ]^\bfd $,
$ p \in ( 0, \infty ) $
that
$ R $
is a
$ ( \cB( [ -B, B ]^\fd ) \otimes \cF ) $/$ \cB( \R ) $-measurable function
and
\begin{equation}
\label{eq:prop:minimum_MC_rate}
\begin{split}
&
\bigl(
\E\bigl[
    \min\nolimits_{ k \in \{ 1, 2, \ldots, K \} } \lvert R( P( \Theta_k ) ) - R( P( \theta ) ) \rvert^p
\bigr]
\bigr)^{ \nicefrac{1}{p} }
\\ &
\leq
\frac{ L ( 2 B ) \max\{ 1, ( \nicefrac{p}{\fd} )^{ \nicefrac{1}{\fd} } \} }{ K^{ \nicefrac{1}{\fd} } }
=
\frac{
4 ( v - u ) b
\bfL
( \lVert \bfl \rVert_{ \infty } + 1 )^\bfL
B^{ \bfL } 
\max\{ 1, ( \nicefrac{p}{\fd} )^{ \nicefrac{1}{\fd} } \}
}{
K^{ \nicefrac{1}{\fd} }
}
.
\end{split}
\end{equation}
The fact that
$ P $ is a
$ \cB( [ -B, B ]^\bfd ) $/$ \cB( [ -B, B ]^\fd ) $-measurable function
and
\cref{eq:clipped_parameters}
hence\linebreak
show~\cref{item:cor:minimum_MC_rate:1}.
In addition,
\cref{eq:clipped_parameters},
\cref{eq:prop:minimum_MC_rate},
and
the fact that
$ 2
\leq \fd
=
\sum_{i=1}^{\bfL} \bfl_i( \bfl_{ i - 1 } + 1 )
\leq
\bfL
( \lVert \bfl \rVert_{ \infty } + 1 )^2 $
yield for all
$ \theta \in [ -B, B ]^\bfd $,
$ p \in ( 0, \infty ) $
that
\begin{align*}
&
\bigl(
\E\bigl[
    \min\nolimits_{ k \in \{ 1, 2, \ldots, K \} } \lvert \cR( \Theta_k ) - \cR( \theta ) \rvert^p
\bigr]
\bigr)^{ \nicefrac{1}{p} }
\\ & \yesnumber
=
\bigl(
\E\bigl[
    \min\nolimits_{ k \in \{ 1, 2, \ldots, K \} } \lvert R( P( \Theta_k ) ) - R( P( \theta ) ) \rvert^p
\bigr]
\bigr)^{ \nicefrac{1}{p} }
\\ &
\leq
\frac{
4 ( v - u ) b
\bfL
( \lVert \bfl \rVert_{ \infty } + 1 )^\bfL
B^{ \bfL }
\sqrt{ \max\{ 1, \nicefrac{p}{\fd} \} }
}{
K^{ \nicefrac{1}{\fd} }
}
\leq
\frac{
4 ( v - u ) b
\bfL
( \lVert \bfl \rVert_{ \infty } + 1 )^\bfL
B^{ \bfL }
\max\{ 1, p \}
}{
K^{ [ \bfL^{-1} ( \lVert \bfl \rVert_{ \infty } + 1 )^{-2} ] }
}
.
\end{align*}
This establishes~\cref{item:cor:minimum_MC_rate:2}.
\end{cproof}

\section{Analysis of the overall error}
\label{sec:DL_risk_minimisation}

In \cref{sec:overall_error_strong_rate} below
we present the main result of this article, \cref{thm:main},
that provides an estimate for the overall $ L^2 $-error
arising in deep learning based empirical risk minimisation
with quadratic loss function
in the probabilistically strong sense
and
that covers the case
where the underlying DNNs are trained using
a general stochastic optimisation algorithm
with random initialisation.

In order to prove \cref{thm:main},
we require a link to combine the results from \cref{sec:approximation_error,sec:generalisation_error,sec:optimisation_error},
which is given in \cref{sec:overall_error} below.
More specifically,
\cref{prop:error_decomposition} in \cref{sec:overall_error}
shows that the overall error
can be decomposed
into three different error sources:
the \emph{approximation error} (cf.\ \cref{sec:approximation_error}),
the \emph{worst-case generalisation error} (cf.\ \cref{sec:generalisation_error}),
and the \emph{optimisation error} (cf.\ \cref{sec:optimisation_error}).
\cref{prop:error_decomposition}
is a consequence of the well-known
bias--variance decomposition
(cf., e.g.,
Beck, Jentzen, \& Kuckuck~\cite[Lemma~4.1]{BeckJentzenKuckuck2019arXiv}
or
Berner, Grohs, \& Jentzen~\cite[Lemma~2.2]{BernerGrohsJentzen2018arXiv})
and is very similar to \cite[Lemma~4.3]{BeckJentzenKuckuck2019arXiv}.

Thereafter,
\cref{sec:overall_error_strong_rate}
is devoted to strong convergence results
for deep learning based empirical risk minimisation
with quadratic loss function
where a general stochastic approximation algorithm
with random initialisation
is allowed to be the employed optimisation method.
Apart from the main result (cf.\ \cref{thm:main}),
\cref{sec:overall_error_strong_rate}
also
includes on the one hand
\cref{prop:main},
which combines
the overall error decomposition (cf.\ \cref{prop:error_decomposition})
with our convergence result for the generalisation error (cf.\ \cref{cor:generalisation_error} in \cref{sec:generalisation_error})
and our convergence result for the optimisation error (cf.\ \cref{cor:minimum_MC_rate} in \cref{sec:optimisation_error}),
and on the other hand
\cref{cor:min_NN_architecture},
which replaces the architecture parameter $ A \in ( 0, \infty ) $
in \cref{thm:main} (cf.\ \cref{prop:approximation_error})
by the minimum of the depth parameter $ \bfL \in \N $
and the hidden layer sizes $ \bfl_1, \bfl_2, \ldots, \bfl_{ \bfL - 1 } \in \N $
of the trained DNN
(cf.\ \cref{eq:choice_A} below).

Finally,
in \cref{sec:SGD}
we present three more strong convergence results for the special case
where SGD
with random initialisation
is the employed optimisation method.
In particular,
\cref{cor:SGD_full_error}
specifies \cref{cor:min_NN_architecture}
to this special case,
\cref{cor:SGD_L1}
provides a convergence estimate
for the expectation of the $ L^1 $-distance
between the trained DNN and the target function,
and
\cref{cor:SGD_simplfied} reaches
an analogous conclusion in a simplified setting.

\subsection{Overall error decomposition}
\label{sec:overall_error}

\begin{proposition}
\label{prop:error_decomposition}
Let
$ d, \bfd, \bfL, M, K, N \in \N $,
$ B \in [ 0, \infty ) $,
$ u \in \R $,
$ v \in ( u, \infty ) $,
$ \bfl = ( \bfl_0, \bfl_1, \ldots, \bfl_\bfL ) \in \N^{ \bfL + 1 } $,
$ \bfN \subseteq \{ 0, 1, \ldots, N \} $,
$ D \subseteq \R^d $,
assume
$ 0 \in \bfN $,
$ \bfl_0 = d $,
$ \bfl_\bfL = 1 $,
and
$ \bfd \geq \sum_{i=1}^{\bfL} \bfl_i( \bfl_{ i - 1 } + 1 ) $,
let
$ ( \Omega, \cF, \P ) $
be a probability space,
let
$ X_j \colon \Omega \to D $,
$ j \in \{ 1, 2, \ldots, M \} $,
and
$ Y_j \colon \Omega \to [ u, v ] $,
$ j \in \{ 1, 2, \ldots, M \} $,
be random variables,
let
$ \cE \colon D \to [ u, v ] $
be a
$ \cB( D ) $/$ \cB( [ u, v ] ) $-measurable function,
assume that it holds
$ \P $-a.s.\
that
$ \cE( X_1 )
= \E[ Y_1 \vert X_1 ] $,
let
$ \Theta_{ k, n } \colon \Omega \to \R^{\bfd} $,
$ k, n \in \N_0 $,
satisfy
$ \bigl( \bigcup_{ k = 1 }^{ \infty }
\Theta_{ k, 0 }( \Omega ) \bigr)
\subseteq [ -B, B ]^\bfd $,
let
$ \bfR \colon \R^\bfd \to [ 0, \infty ) $
satisfy for all
$ \theta \in \R^\bfd $
that
$ \bfR( \theta )
= \E[ \lvert \clippedNN{\theta}{\bfl}{u}{v}( X_1 ) - Y_1 \rvert^2 ] $,
and
let
$ \cR \colon \R^{\bfd} \times \Omega \to [ 0, \infty ) $
and
$ \bfk \colon \Omega \to ( \N_0 )^2 $
satisfy for all
$ \theta \in \R^{\bfd} $,
$ \omega \in \Omega $
that
\begin{align}
\cR( \theta, \omega )
=
\frac{1}{M}
\biggl[
\smallsum_{j=1}^M
    \lvert \clippedNN{\theta}{\bfl}{u}{v}( X_j( \omega ) ) - Y_j( \omega ) \rvert^2
\biggr]
\qquad
\text{and}
\\
\label{eq:choice_NN_parameters}
\bfk( \omega ) \in
\argmin\nolimits_{ ( k, n ) \in \{ 1, 2, \ldots, K \} \times \bfN, \, \lVert \Theta_{ k, n }( \omega ) \rVert_{ \infty } \leq B }
\cR( \Theta_{ k, n }( \omega ), \omega )
\end{align}
(cf.~\cref{def:clipped_NN,def:p-norm}).
Then
it holds for all
$ \vartheta \in [ -B, B ]^\bfd $
that
\begin{equation}
\begin{split}
&
\int_D
    \lvert \clippedNN{\Theta_\bfk}{\bfl}{u}{v}( x ) - \cE( x ) \rvert^2
\, \P_{ X_1 }( \dd x )
\\ &
\leq
\bigl[
\sup\nolimits_{ x \in D }
    \lvert \clippedNN{\vartheta}{\bfl}{u}{v}( x ) - \cE( x ) \rvert^2
\bigr]
+
2\bigl[
\sup\nolimits_{ \theta \in [ -B, B ]^\bfd }
\lvert
    \cR( \theta )
    -
    \bfR( \theta )
\rvert
\bigr]
\\ & \quad
+
\min\nolimits_{ ( k, n ) \in \{ 1, 2, \ldots, K \} \times \bfN, \, \lVert \Theta_{ k, n } \rVert_{ \infty } \leq B }
\lvert \cR( \Theta_{ k, n } )
-
\cR( \vartheta ) \rvert
.
\end{split}
\end{equation}
\end{proposition}
\begin{cproof}{prop:error_decomposition}
Throughout this proof
let
$ \scrR \colon \cL^2( \P_{ \smash{ X_1 } }; \R ) \to [ 0, \infty ) $
satisfy for all
$ f \in \cL^2( \P_{ \smash{ X_1 } }; \R ) $
that
$ \scrR( f ) = \E[ \lvert f( X_1 ) - Y_1 \rvert^2 ] $.
Observe that
the assumption that
$ \forall \, \omega \in \Omega \colon
Y_1( \omega ) \in [ u, v ] $
and
the fact that
$ \forall \, \theta \in \R^\bfd,
\, x \in \R^d \colon
\clippedNN{\theta}{\bfl}{u}{v}( x ) \in [ u, v ] $
ensure for all
$ \theta \in \R^\bfd $
that
$ \E[ \lvert Y_1 \rvert^2 ]
\leq v^2 < \infty $
and
\begin{equation}
\int_D
        \lvert \clippedNN{\theta}{\bfl}{u}{v}( x ) \rvert^2
\, \P_{ X_1 }( \dd x )
=
\E\bigl[
    \lvert \clippedNN{\theta}{\bfl}{u}{v}( X_1 ) \rvert^2
\bigr]
\leq
v^2 < \infty
.
\end{equation}
The bias--variance decomposition
(cf., e.g.,
Beck, Jentzen, \& Kuckuck~\cite[(iii) in Lem\-ma~4.1]{BeckJentzenKuckuck2019arXiv}
with
$ ( \Omega, \cF, \P ) \leftarrow ( \Omega, \cF, \P ) $,
$ ( S, \cS ) \leftarrow ( D, \cB( D ) ) $,
$ X \leftarrow X_1 $,
$ Y \leftarrow
( \Omega \ni \omega \mapsto Y_1( \omega ) \in \R ) $,
$ \cE \leftarrow \scrR $,
$ f \leftarrow \clippedNN{\theta}{\bfl}{u}{v} \vert_D $,
$ g \leftarrow \clippedNN{\vartheta}{\bfl}{u}{v} \vert_D $
for
$ \theta, \vartheta \in \R^\bfd $
in the notation of~\cite[(iii) in Lemma~4.1]{BeckJentzenKuckuck2019arXiv})
hence
proves for all
$ \theta, \vartheta \in \R^\bfd $
that
\begin{align*}
&
\int_D
    \lvert \clippedNN{\theta}{\bfl}{u}{v}( x ) - \cE( x ) \rvert^2
\, \P_{ X_1 }( \dd x )
\\ &
=
\E\bigl[
    \lvert \clippedNN{\theta}{\bfl}{u}{v}( X_1 )
    -
    \cE( X_1 ) \rvert^2
\bigr]
=
\E\bigl[
    \lvert \clippedNN{\theta}{\bfl}{u}{v}( X_1 )
    -
    \E[ Y_1 \vert X_1 ] \rvert^2
\bigr]
\\ & \yesnumber
=
\E\bigl[
    \lvert \clippedNN{\vartheta}{\bfl}{u}{v}( X_1 )
    -
    \E[ Y_1 \vert X_1 ] \rvert^2
\bigr]
+
\scrR( \clippedNN{\theta}{\bfl}{u}{v} \vert_D )
-
\scrR( \clippedNN{\vartheta}{\bfl}{u}{v} \vert_D )
\\ &
=
\E\bigl[
    \lvert \clippedNN{\vartheta}{\bfl}{u}{v}( X_1 )
    -
    \cE( X_1 ) \rvert^2
\bigr]
+
\E\bigl[
    \lvert \clippedNN{\theta}{\bfl}{u}{v}( X_1 )
    -
    Y_1 \rvert^2
\bigr]
-
\E\bigl[
    \lvert \clippedNN{\vartheta}{\bfl}{u}{v}( X_1 )
    -
    Y_1 \rvert^2
\bigr]
\\ &
=
\int_D
    \lvert \clippedNN{\vartheta}{\bfl}{u}{v}( x ) - \cE( x ) \rvert^2
\, \P_{ X_1 }( \dd x )
+
\bfR( \theta )
-
\bfR( \vartheta )
.
\end{align*}
This
implies for all
$ \theta, \vartheta \in \R^\bfd $
that
\begin{align*}
\label{eq:first_error_estimate}
&
\int_D
    \lvert \clippedNN{\theta}{\bfl}{u}{v}( x ) - \cE( x ) \rvert^2
\, \P_{ X_1 }( \dd x )
\\ &
=
\int_D
    \lvert \clippedNN{\vartheta}{\bfl}{u}{v}( x ) - \cE( x ) \rvert^2
\, \P_{ X_1 }( \dd x )
-
[ \cR( \theta )
- \bfR( \theta ) ]
+
\cR( \vartheta ) - \bfR( \vartheta )
+
\cR( \theta )
-
\cR( \vartheta )
\\ &
\leq
\int_D
    \lvert \clippedNN{\vartheta}{\bfl}{u}{v}( x ) - \cE( x ) \rvert^2
\, \P_{ X_1 }( \dd x )
+
\lvert
    \cR( \theta )
    -
    \bfR( \theta )
\rvert
+
\lvert
    \cR( \vartheta )
    -
    \bfR( \vartheta )
\rvert
+
\cR( \theta )
-
\cR( \vartheta )
\\ & \yesnumber
\leq
\int_D
    \lvert \clippedNN{\vartheta}{\bfl}{u}{v}( x ) - \cE( x ) \rvert^2
\, \P_{ X_1 }( \dd x )
+
2\bigl[
\max\nolimits_{ \eta \in \{ \theta, \vartheta \} }
\lvert
    \cR( \eta )
    -
    \bfR( \eta )
\rvert
\bigr]
+
\cR( \theta )
-
\cR( \vartheta )
.
\end{align*}
Next note that
the fact that
$ \forall \, \omega \in \Omega \colon
\lVert \Theta_{ \bfk( \omega ) }( \omega ) \rVert_{ \infty } \leq B $
ensures for all
$ \omega \in \Omega $
that
$ \Theta_{ \bfk( \omega ) }( \omega ) \in [ -B, B ]^\bfd $.
Combining
\cref{eq:first_error_estimate}
with
\cref{eq:choice_NN_parameters}
hence
establishes for all
$ \vartheta \in [ -B, B ]^\bfd $
that
\begin{equation}
\begin{split}
&
\int_D
    \lvert \clippedNN{\Theta_\bfk}{\bfl}{u}{v}( x ) - \cE( x ) \rvert^2
\, \P_{ X_1 }( \dd x )
\\ &
\leq
\int_D
    \lvert \clippedNN{\vartheta}{\bfl}{u}{v}( x ) - \cE( x ) \rvert^2
\, \P_{ X_1 }( \dd x )
+
2\bigl[
\sup\nolimits_{ \theta \in [ -B, B ]^\bfd }
\lvert
    \cR( \theta )
    -
    \bfR( \theta )
\rvert
\bigr]
+
\cR( \Theta_\bfk )
-
\cR( \vartheta )
\\ &
=
\int_D
    \lvert \clippedNN{\vartheta}{\bfl}{u}{v}( x ) - \cE( x ) \rvert^2
\, \P_{ X_1 }( \dd x )
+
2\bigl[
\sup\nolimits_{ \theta \in [ -B, B ]^\bfd }
\lvert
    \cR( \theta )
    -
    \bfR( \theta )
\rvert
\bigr]
\\ & \quad
+
\min\nolimits_{ ( k, n ) \in \{ 1, 2, \ldots, K \} \times \bfN, \, \lVert \Theta_{ k, n } \rVert_{ \infty } \leq B }
[ \cR( \Theta_{ k, n } )
-
\cR( \vartheta ) ]
\\ &
\leq
\bigl[
\sup\nolimits_{ x \in D }
    \lvert \clippedNN{\vartheta}{\bfl}{u}{v}( x ) - \cE( x ) \rvert^2
\bigr]
+
2\bigl[
\sup\nolimits_{ \theta \in [ -B, B ]^\bfd }
\lvert
    \cR( \theta )
    -
    \bfR( \theta )
\rvert
\bigr]
\\ & \quad
+
\min\nolimits_{ ( k, n ) \in \{ 1, 2, \ldots, K \} \times \bfN, \, \lVert \Theta_{ k, n } \rVert_{ \infty } \leq B }
\lvert \cR( \Theta_{ k, n } )
-
\cR( \vartheta ) \rvert
.
\end{split}
\end{equation}
\end{cproof}

\subsection{Full strong error analysis for the training of DNNs}
\label{sec:overall_error_strong_rate}

\begin{lemma}
\label{lem:measurability}
Let
$ d, \bfd, \bfL \in \N $,
$ p \in [ 0, \infty ) $,
$ u \in [ -\infty, \infty ) $,
$ v \in ( u, \infty ] $,
$ \bfl = ( \bfl_0, \bfl_1, \ldots, \bfl_\bfL ) \in \N^{ \bfL + 1 } $,
$ D \subseteq \R^d $,
assume
$ \bfl_0 = d $,
$ \bfl_\bfL = 1 $,
and
$ \bfd \geq \sum_{i=1}^{\bfL} \bfl_i( \bfl_{ i - 1 } + 1 ) $,
let
$ \cE \colon D \to \R $
be a
$ \cB( D ) $/$ \cB( \R ) $-measurable function,
let
$ ( \Omega, \cF, \P ) $
be a probability space,
and
let
$ X \colon \Omega \to D $,
$ \bfk \colon \Omega \to ( \N_0 )^2 $,
and
$ \Theta_{ k, n } \colon \Omega \to \R^{\bfd} $,
$ k, n \in \N_0 $,
be random variables.
Then
\begin{enumerate}[(i)]
\item
\label{item:lem:measurability:1}
it holds that
the function
$ \R^\bfd \times \R^d \ni ( \theta, x )
\mapsto
\clippedNN{\theta}{\bfl}{u}{v}( x )
\in \R $
is $ ( \cB( \R^\bfd ) \otimes \cB( \R^d ) ) $/$ \cB( \R ) $-measurable,
\item
\label{item:lem:measurability:2}
it holds that
the function
$ \Omega \ni \omega
\mapsto
\Theta_{ \bfk( \omega ) }( \omega )
\in \R^\bfd $
is
$ \cF $/$ \cB( \R^\bfd ) $-measurable,
and
\item
\label{item:lem:measurability:3}
it holds that
the function
\begin{equation}
\Omega \ni \omega
\mapsto
\int_D
    \lvert \clippedNN{ \smash{ \Theta_{ \smash{ \bfk( \omega ) } }( \omega ) } }{ \bfl }{u}{v}( x ) - \cE( x ) \rvert^p
\, \P_{ X }( \dd x )
\in [ 0, \infty ]
\end{equation}
is
$ \cF $/$ \cB( [ 0, \infty ] ) $-measurable
\end{enumerate}
(cf.~\cref{def:clipped_NN}).
\end{lemma}
\begin{cproof}{lem:measurability}
First,
observe that
Beck, Jentzen, \& Kuckuck~\cite[Corollary~2.37]{BeckJentzenKuckuck2019arXiv}
(with
$ a \leftarrow -\lVert x \rVert_\infty $,
$ b \leftarrow \lVert x \rVert_\infty $,
$ u \leftarrow u $,
$ v \leftarrow v $,
$ d \leftarrow \bfd $,
$ L \leftarrow \bfL $,
$ l \leftarrow \bfl $
for
$ x \in \R^d $
in the notation of~\cite[Corollary~2.37]{BeckJentzenKuckuck2019arXiv})
demonstrates for all
$ x \in \R^d $,
$ \theta, \vartheta \in \R^\bfd $
that
\begin{equation}
\begin{split}
&
\lvert
    \clippedNN{\theta}{\bfl}{u}{v}( x ) - \clippedNN{\vartheta}{\bfl}{u}{v}( x )
\rvert
\leq
\sup\nolimits_{ y \in [ -\lVert x \rVert_\infty, \lVert x \rVert_\infty ]^{ \bfl_0 } }
    \lvert
        \clippedNN{\theta}{\bfl}{u}{v}( y ) - \clippedNN{\vartheta}{\bfl}{u}{v}( y )
    \rvert
\\ &
\leq
\bfL
\max\{ 1, \lVert x \rVert_\infty \}
( \lVert \bfl \rVert_{ \infty } + 1 )^\bfL
( \max\{ 1, \lVert \theta \rVert_{ \infty }, \lVert \vartheta \rVert_{ \infty } \} )^{ \bfL - 1 }
\lVert \theta - \vartheta \rVert_{ \infty }
\end{split}
\end{equation}
(cf.~\cref{def:p-norm}).
This
implies for all
$ x \in \R^d $
that the function
\begin{equation}
\label{eq:NN_continuous}
\R^\bfd \ni \theta
\mapsto \clippedNN{\theta}{\bfl}{u}{v}( x ) \in \R
\end{equation}
is continuous.
In addition,
the fact that
$ \forall \, \theta \in \R^\bfd \colon
\clippedNN{\theta}{\bfl}{u}{v} \in C( \R^d, \R ) $
ensures for all
$ \theta \in \R^\bfd $
that
the function
$ \R^d \ni x
\mapsto \clippedNN{\theta}{\bfl}{u}{v}( x ) \in \R $
is $ \cB( \R^d ) $/$ \cB( \R ) $-measurable.
This,
\cref{eq:NN_continuous},
the fact that
$ ( \R^\bfd, \lVert \cdot \rVert_\infty \vert_{ \R^\bfd } ) $
is a separable normed $ \R $-vector space,
and, e.g.,
Aliprantis \& Border~\cite[Lemma~4.51]{AliprantisBorder2006}
(see also, e.g.,
Beck et al.~\cite[Lemma~2.4]{BeckBeckerGrohsJaafariJentzen2018arXiv})
show~\cref{item:lem:measurability:1}.

Second,
we prove~\cref{item:lem:measurability:2}.
For this
let
$ \Xi \colon \Omega \to \R^\bfd $
satisfy for all
$ \omega \in \Omega $
that
$ \Xi( \omega ) = \Theta_{ \bfk( \omega ) }( \omega ) $.
Observe that
the assumption that
$ \Theta_{ k, n } \colon \Omega \to \R^{\bfd} $,
$ k, n \in \N_0 $,
and
$ \bfk \colon \Omega \to ( \N_0 )^2 $
are random variables
establishes for all
$ U \in \cB( \R^\bfd ) $
that
\begin{equation}
\begin{split}
\Xi^{ -1 }( U )
& =
\{ \omega \in \Omega \colon
\Xi( \omega ) \in U \}
=
\{ \omega \in \Omega \colon
\Theta_{ \bfk( \omega ) }( \omega ) \in U \}
\\ &
=
\bigl\{ \omega \in \Omega \colon
\bigl[ \exists \, k, n \in \N_0 \colon
( [ \Theta_{ k, n }( \omega ) \in U ]
\land
[ \bfk( \omega ) = ( k, n ) ]
) \bigr] \bigr\}
\\ &
=
\smallbigcup_{ k = 0 }^\infty
\smallbigcup_{ n = 0 }^\infty
    \bigl(
        \{ \omega \in \Omega \colon \Theta_{ k, n }( \omega ) \in U \}
        \cap
        \{ \omega \in \Omega \colon \bfk( \omega ) = ( k, n ) \}    
    \bigr)
\\ &
=
\smallbigcup_{ k = 0 }^\infty
\smallbigcup_{ n = 0 }^\infty
    \bigl(
         [ ( \Theta_{ k, n } )^{ -1 }( U ) ]
        \cap
        [ \bfk^{ -1 }( \{ ( k, n ) \} ) ]
    \bigr)
\in
\cF
.
\end{split}
\end{equation}
This implies~\cref{item:lem:measurability:2}.

Third,
note that
\cref{item:lem:measurability:1}--\cref{item:lem:measurability:2}
yield that
the function
$ \Omega \times \R^d \ni ( \omega, x )
\mapsto
\clippedNN{ \smash{ \Theta_{ \smash{ \bfk( \omega ) } }( \omega ) } }{ \smash{ \bfl } }{u}{v}( x )
\in \R $
is $ ( \cF \otimes \cB( \R^d ) ) $/$ \cB( \R ) $-measurable.
This
and
the assumption that
$ \cE \colon D \to \R $
is
$ \cB( D ) $/$ \cB( \R ) $-measurable
demonstrate that
the function
$ \Omega \times D \ni ( \omega, x )
\mapsto
\lvert \clippedNN{ \smash{ \Theta_{ \smash{ \bfk( \omega ) } }( \omega ) } }{ \smash{ \bfl } }{u}{v}( x ) - \cE( x ) \rvert^p
\in [ 0, \infty ) $
is $ ( \cF \otimes \cB( D ) ) $/$ \cB( [ 0, \infty ) ) $-measurable.
Tonelli's theorem
hence establishes~\cref{item:lem:measurability:3}.
\end{cproof}

\begin{proposition}
\label{prop:main}
Let
$ d, \bfd, \bfL, M, K, N \in \N $,
$ b, c \in [ 1, \infty ) $,
$ B \in [ c, \infty ) $,
$ u \in \R $,
$ v \in ( u, \infty ) $,
$ \bfl = ( \bfl_0, \bfl_1, \ldots, \bfl_\bfL ) \in \N^{ \bfL + 1 } $,
$ \bfN \subseteq \{ 0, 1, \ldots, N \} $,
$ D \subseteq [ -b, b ]^d $,
assume
$ 0 \in \bfN $,
$ \bfl_0 = d $,
$ \bfl_\bfL = 1 $,
and
$ \bfd \geq \sum_{i=1}^{\bfL} \bfl_i( \bfl_{ i - 1 } + 1 ) $,
let
$ ( \Omega, \cF, \P ) $
be a probability space,
let
$ X_j \colon \Omega \to D $,
$ j \in \N $,
and
$ Y_j \colon \Omega \to [ u, v ] $,
$ j \in \N $,
be functions,
assume that
$ ( X_j, Y_j ) $,
$ j \in \{ 1, 2, \ldots, M \} $,
are i.i.d.\ random variables,
let
$ \cE \colon D \to [ u, v ] $
be a
$ \cB( D ) $/$ \cB( [ u, v ] ) $-measurable function,
assume that it holds
$ \P $-a.s.\
that
$ \cE( X_1 )
= \E[ Y_1 \vert X_1 ] $,
let
$ \Theta_{ k, n } \colon \Omega \to \R^{\bfd} $,
$ k, n \in \N_0 $,
and
$ \bfk \colon \Omega \to ( \N_0 )^2 $
be random variables,
assume
$ \bigl( \bigcup_{ k = 1 }^{ \infty }
\Theta_{ k, 0 }( \Omega ) \bigr)
\subseteq [ -B, B ]^\bfd $,
assume that
$ \Theta_{ k, 0 } $,
$ k \in \{ 1, 2, \ldots, K \} $,
are i.i.d.,
assume that
$ \Theta_{ 1, 0 } $ is continuous uniformly distributed on $ [ -c, c ]^\bfd $,
and let
$ \cR \colon \R^{\bfd} \times \Omega \to [ 0, \infty ) $
satisfy for all
$ \theta \in \R^{\bfd} $,
$ \omega \in \Omega $
that
\begin{align}
\cR( \theta, \omega )
=
\frac{1}{M}
\biggl[
\smallsum_{j=1}^M
    \lvert \clippedNN{\theta}{\bfl}{u}{v}( X_j( \omega ) ) - Y_j( \omega ) \rvert^2
\biggr]
\qquad
\text{and}
\\
\bfk( \omega ) \in
\argmin\nolimits_{ ( k, n ) \in \{ 1, 2, \ldots, K \} \times \bfN, \, \lVert \Theta_{ k, n }( \omega ) \rVert_{ \infty } \leq B }
\cR( \Theta_{ k, n }( \omega ), \omega )
\end{align}
(cf.~\cref{def:clipped_NN,def:p-norm}).
Then
it holds for all
$ p \in ( 0, \infty ) $
that
\begin{equation}
\begin{split}
&
\Bigl(
\E\Bigl[
\Bigl( \medint{D}
        \lvert \clippedNN{\Theta_{\bfk}}{\bfl}{u}{v}( x ) - \cE( x ) \rvert^2
    \, \P_{ X_1 }( \dd x )
\Bigr)^{ \! p \, }
\Bigr]
\Bigr)^{ \! \nicefrac{1}{p} }
\\ &
\leq
\bigl[
\inf\nolimits_{ \theta \in [ -c, c ]^\bfd }
\sup\nolimits_{ x \in D }
    \lvert \clippedNN{\theta}{\bfl}{u}{v}( x ) - \cE( x ) \rvert^{2}
\bigr]
+
\frac{
4 ( v - u ) b
\bfL
( \lVert \bfl \rVert_{ \infty } + 1 )^\bfL
c^{ \bfL }
\max\{ 1, p \}
}{
K^{ [ \bfL^{-1} ( \lVert \bfl \rVert_{ \infty } + 1 )^{-2} ] }
}
\\ & \quad
+
\frac{
    18 \max\{ 1, ( v - u )^2 \}
    \bfL ( \lVert \bfl \rVert_{ \infty } + 1 )^2
    \max\{
        p,
        \ln( 3 M B b )
    \}
}{ \sqrt{M} }
\\ &
\leq
\bigl[
\inf\nolimits_{ \theta \in [ -c, c ]^\bfd }
\sup\nolimits_{ x \in D }
    \lvert \clippedNN{\theta}{\bfl}{u}{v}( x ) - \cE( x ) \rvert^{2}
\bigr]
\\ & \quad
+
\frac{
    20 \max\{ 1, ( v - u )^2 \}
    b \bfL ( \lVert \bfl \rVert_{ \infty } + 1 )^{ \bfL + 1 } B^\bfL
    \max\{ p, \ln( 3 M ) \}
}{
\min\{
    \sqrt{M},
    K^{ [ \bfL^{-1} ( \lVert \bfl \rVert_{ \infty } + 1 )^{-2} ] }
\} }
\end{split}
\end{equation}
(cf.~\cref{item:lem:measurability:3} in \cref{lem:measurability}).
\end{proposition}
\begin{cproof2}{prop:main}
Throughout this proof
let
$ \bfR \colon \R^\bfd \to [ 0, \infty ) $
satisfy for all
$ \theta \in \R^\bfd $
that
$ \bfR( \theta )
= \E[ \lvert \clippedNN{\theta}{\bfl}{u}{v}( X_1 ) - Y_1 \rvert^2 ] $.
First of all,
observe that
the assumption that
$ \bigl( \bigcup_{ k = 1 }^{ \infty }
\Theta_{ k, 0 }( \Omega ) \bigr)
\subseteq [ -B, B ]^\bfd $,
the assumption that
$ 0 \in \bfN $,
and
\cref{prop:error_decomposition}
show for all
$ \vartheta \in [ -B, B ]^\bfd $
that%
\begin{equation}
\label{eq:reduce_to_minimum_MC}
\begin{split}
&
\int_D
    \lvert \clippedNN{\Theta_\bfk}{\bfl}{u}{v}( x ) - \cE( x ) \rvert^2
\, \P_{ X_1 }( \dd x )
\\ &
\leq
\bigl[
\sup\nolimits_{ x \in D }
    \lvert \clippedNN{\vartheta}{\bfl}{u}{v}( x ) - \cE( x ) \rvert^2
\bigr]
+
2\bigl[
\sup\nolimits_{ \theta \in [ -B, B ]^\bfd }
\lvert
    \cR( \theta )
    -
    \bfR( \theta )
\rvert
\bigr]
\\ & \quad
+
\min\nolimits_{ ( k, n ) \in \{ 1, 2, \ldots, K \} \times \bfN, \, \lVert \Theta_{ k, n } \rVert_{ \infty } \leq B }
\lvert \cR( \Theta_{ k, n } )
-
\cR( \vartheta ) \rvert
\\ &
\leq
\bigl[
\sup\nolimits_{ x \in D }
    \lvert \clippedNN{\vartheta}{\bfl}{u}{v}( x ) - \cE( x ) \rvert^2
\bigr]
+
2\bigl[
\sup\nolimits_{ \theta \in [ -B, B ]^\bfd }
\lvert
    \cR( \theta )
    -
    \bfR( \theta )
\rvert
\bigr]
\\ & \quad
+
\min\nolimits_{ k \in \{ 1, 2, \ldots, K \}, \, \lVert \Theta_{ k, 0 } \rVert_{ \infty } \leq B }
\lvert \cR( \Theta_{ k, 0 } )
-
\cR( \vartheta ) \rvert
\\ &
=
\bigl[
\sup\nolimits_{ x \in D }
    \lvert \clippedNN{\vartheta}{\bfl}{u}{v}( x ) - \cE( x ) \rvert^2
\bigr]
+
2\bigl[
\sup\nolimits_{ \theta \in [ -B, B ]^\bfd }
\lvert
    \cR( \theta )
    -
    \bfR( \theta )
\rvert
\bigr]
\\ & \quad
+
\min\nolimits_{ k \in \{ 1, 2, \ldots, K \} }
\lvert \cR( \Theta_{ k, 0 } )
-
\cR( \vartheta ) \rvert
.
\end{split}
\end{equation}
Minkowski's inequality
hence establishes for all
$ p \in [ 1, \infty ) $,
$ \vartheta \in [ -c, c ]^\bfd
\subseteq [ -B, B ]^\bfd $
that%
\begin{equation}
\label{eq:strong_error_decomposition}
\begin{split}
&
\Bigl(
\E\Bigl[
\Bigl( \medint{D}
        \lvert \clippedNN{\Theta_{\bfk}}{\bfl}{u}{v}( x ) - \cE( x ) \rvert^2
    \, \P_{ X_1 }( \dd x )
\Bigr)^{ \! p \, }
\Bigr]
\Bigr)^{ \! \nicefrac{1}{p} }
\\ &
\leq
\bigl(
\E\bigl[
\sup\nolimits_{ x \in D }
    \lvert \clippedNN{\vartheta}{\bfl}{u}{v}( x ) - \cE( x ) \rvert^{2p}
\bigr]
\bigr)^{ \nicefrac{1}{p} }
+
2 \bigl(
\E\bigl[
\sup\nolimits_{ \theta \in [ -B, B ]^\bfd }
\lvert
    \cR( \theta )
    -
    \bfR( \theta )
\rvert^p
\bigr]
\bigr)^{ \nicefrac{1}{p} }
\\ & \quad
+
\bigl(
\E\bigl[
\min\nolimits_{ k \in \{ 1, 2, \ldots, K \} }
\lvert \cR( \Theta_{ k, 0 } )
-
\cR( \vartheta ) \rvert^p
\bigr]
\bigr)^{ \nicefrac{1}{p} }
\\ &
\leq
\bigl[
\sup\nolimits_{ x \in D }
    \lvert \clippedNN{\vartheta}{\bfl}{u}{v}( x ) - \cE( x ) \rvert^{2}
\bigr]
+
2 \bigl(
\E\bigl[
\sup\nolimits_{ \theta \in [ -B, B ]^\bfd }
\lvert
    \cR( \theta )
    -
    \bfR( \theta )
\rvert^p
\bigr]
\bigr)^{ \nicefrac{1}{p} }
\\ & \quad
+
\sup\nolimits_{ \theta \in [ -c, c ]^\bfd }
\bigl(
\E\bigl[
\min\nolimits_{ k \in \{ 1, 2, \ldots, K \} }
\lvert \cR( \Theta_{ k, 0 } )
-
\cR( \theta ) \rvert^p
\bigr]
\bigr)^{ \nicefrac{1}{p} }
\end{split}
\end{equation}
(cf.~\cref{item:cor:generalisation_error:1} in \cref{cor:generalisation_error}
and
\cref{item:cor:minimum_MC_rate:1} in \cref{cor:minimum_MC_rate}).
Next note that
\cref{cor:generalisation_error}
(with
$ v \leftarrow \max\{ u + 1, v \} $,
$ \bfR \leftarrow \bfR \vert_{ [ -B, B ]^\bfd } $,
$ \cR \leftarrow \cR \vert_{ [ -B, B ]^\bfd \times \Omega } $
in the notation of \cref{cor:generalisation_error})
proves for all
$ p \in ( 0, \infty ) $
that
\begin{equation}
\label{eq:strong_generalisation_error}
\begin{split}
& \bigl(
\E\bigl[
    \sup\nolimits_{ \theta \in [ -B, B ]^\bfd }
    \lvert \cR( \theta ) - \bfR( \theta ) \rvert^p
\bigr]
\bigr)^{ \nicefrac{1}{p} }
\\ &
\leq
\frac{
    9 ( \max\{ u + 1, v \} - u )^2
    \bfL ( \lVert \bfl \rVert_{ \infty } + 1 )^2
    \max\{
        p,
        \ln( 3 M B b )
    \}
}{ \sqrt{M} }
\\ &
=
\frac{
    9 \max\{ 1, ( v - u )^2 \}
    \bfL ( \lVert \bfl \rVert_{ \infty } + 1 )^2
    \max\{
        p,
        \ln( 3 M B b )
    \}
}{ \sqrt{M} }
.
\end{split}
\end{equation}
In addition,
observe that
\cref{cor:minimum_MC_rate}
(with
$ \fd \leftarrow \sum_{i=1}^{\bfL} \bfl_i( \bfl_{ i - 1 } + 1 ) $,
$ B \leftarrow c $,
$ ( \Theta_k )_{ k \in \{ 1, 2, \ldots, K \} }
\allowbreak
\leftarrow
( \Omega \ni \omega \mapsto
\mathbbm{1}_{ \{ \Theta_{ k, 0 } \in [ -c, c ]^\bfd \} }( \omega ) \Theta_{ k, 0 }( \omega ) \in [ -c, c ]^\bfd )_{ k \in \{ 1, 2, \ldots, K \} } $,
$ \cR \leftarrow \cR \vert_{ [ -c, c ]^\bfd \times \Omega } $
in the notation of \cref{cor:minimum_MC_rate})
implies for all
$ p \in ( 0, \infty ) $
that
\begin{equation}
\begin{split}
&
\sup\nolimits_{ \theta \in [ -c, c ]^\bfd }
\bigl(
\E\bigl[
\min\nolimits_{ k \in \{ 1, 2, \ldots, K \} }
    \lvert \cR( \Theta_{ k, 0 } )
    -
    \cR( \theta ) \rvert^p
\bigr]
\bigr)^{ \nicefrac{1}{p} }
\\ &
=
\sup\nolimits_{ \theta \in [ -c, c ]^\bfd }
\bigl(
\E\bigl[
\min\nolimits_{ k \in \{ 1, 2, \ldots, K \} }
    \lvert \cR( \mathbbm{1}_{ \{ \Theta_{ k, 0 } \in [ -c, c ]^\bfd \} } \Theta_{ k, 0 } )
    -
    \cR( \theta ) \rvert^p
\bigr]
\bigr)^{ \nicefrac{1}{p} }
\\ &
\leq
\frac{
4 ( v - u ) b
\bfL
( \lVert \bfl \rVert_{ \infty } + 1 )^\bfL
c^{ \bfL }
\max\{ 1, p \}
}{
K^{ [ \bfL^{-1} ( \lVert \bfl \rVert_{ \infty } + 1 )^{-2} ] }
}
.
\end{split}
\end{equation}
Combining
this,
\cref{eq:strong_error_decomposition},
\cref{eq:strong_generalisation_error},
and the fact that
$ \ln( 3 M B b ) \geq 1 $
with
Jensen's inequality
demonstrates for all
$ p \in ( 0, \infty ) $
that
\begin{equation}
\label{eq:prop:main_derivation}
\begin{split}
&
\Bigl(
\E\Bigl[
\Bigl( \medint{D}
        \lvert \clippedNN{\Theta_{\bfk}}{\bfl}{u}{v}( x ) - \cE( x ) \rvert^2
    \, \P_{ X_1 }( \dd x )
\Bigr)^{ \! p \, }
\Bigr]
\Bigr)^{ \! \nicefrac{1}{p} }
\\ &
\leq
\biggl(
\E\biggl[
\Bigl( \medint{D}
        \lvert \clippedNN{\Theta_{\bfk}}{\bfl}{u}{v}( x ) - \cE( x ) \rvert^2
    \, \P_{ X_1 }( \dd x )
\Bigr)^{ \! \max\{ 1, p \} }
\biggr]
\biggr)^{ \!\! \frac{1}{\max\{ 1, p \}} }
\\ &
\leq
\bigl[
\inf\nolimits_{ \theta \in [ -c, c ]^\bfd }
\sup\nolimits_{ x \in D }
    \lvert \clippedNN{\theta}{\bfl}{u}{v}( x ) - \cE( x ) \rvert^{2}
\bigr]
\\ & \quad
+
\sup\nolimits_{ \theta \in [ -c, c ]^\bfd }
\bigl(
\E\bigl[
\min\nolimits_{ k \in \{ 1, 2, \ldots, K \} }
\lvert \cR( \Theta_{ k, 0 } )
-
\cR( \theta ) \rvert^{ \max\{ 1, p \} }
\bigr]
\bigr)^{ \frac{1}{\max\{ 1, p \}} }
\\ & \quad
+
2 \bigl(
\E\bigl[
\sup\nolimits_{ \theta \in [ -B, B ]^\bfd }
\lvert
    \cR( \theta )
    -
    \bfR( \theta )
\rvert^{ \max\{ 1, p \} }
\bigr]
\bigr)^{ \frac{1}{\max\{ 1, p \}} }
\\ &
\leq
\bigl[
\inf\nolimits_{ \theta \in [ -c, c ]^\bfd }
\sup\nolimits_{ x \in D }
    \lvert \clippedNN{\theta}{\bfl}{u}{v}( x ) - \cE( x ) \rvert^{2}
\bigr]
+
\frac{
4 ( v - u ) b
\bfL
( \lVert \bfl \rVert_{ \infty } + 1 )^\bfL
c^{ \bfL }
\max\{ 1, p \}
}{
K^{ [ \bfL^{-1} ( \lVert \bfl \rVert_{ \infty } + 1 )^{-2} ] }
}
\\ & \quad
+
\frac{
    18 \max\{ 1, ( v - u )^2 \}
    \bfL ( \lVert \bfl \rVert_{ \infty } + 1 )^2
    \max\{
        p,
        \ln( 3 M B b )
    \}
}{ \sqrt{M} }
.
\end{split}
\end{equation}
Moreover,
note that the fact that
$ \forall \, x \in [ 0, \infty ) \colon
x + 1 \leq e^{ x } \leq 3^{ x } $
and the facts that
$ B b \geq 1 $
and
$ M \geq 1 $
ensure that
\begin{equation}
\ln( 3 M B b )
\leq
\ln( 3 M 3^{ B b - 1 } )
=
\ln( 3^{ B b } M )
=
B b \ln( [ 3^{ B b } M ]^{ \nicefrac{1}{ ( B b ) }} )
\leq
B b \ln( 3 M )
.
\end{equation}
The facts that
$ \lVert \bfl \rVert_{ \infty } + 1 \geq 2 $,
$ B \geq c \geq 1 $,
$ \ln( 3 M ) \geq 1 $,
$ b \geq 1 $,
and
$ \bfL \geq 1 $
hence
show for all
$ p \in ( 0, \infty ) $
that
\begin{equation}
\begin{split}
&
\frac{
4 ( v - u ) b
\bfL
( \lVert \bfl \rVert_{ \infty } + 1 )^\bfL
c^{ \bfL }
\max\{ 1, p \}
}{
K^{ [ \bfL^{-1} ( \lVert \bfl \rVert_{ \infty } + 1 )^{-2} ] }
}
\\ &
+
\frac{
    18 \max\{ 1, ( v - u )^2 \}
    \bfL ( \lVert \bfl \rVert_{ \infty } + 1 )^2
    \max\{
        p,
        \ln( 3 M B b )
    \}
}{ \sqrt{M} }
\\ &
\leq
\frac{ 2 ( \lVert \bfl \rVert_{ \infty } + 1 ) \max\{ 1, ( v - u )^2 \}
b \bfL ( \lVert \bfl \rVert_{ \infty } + 1 )^\bfL B^{ \bfL }
\max\{ p, \ln( 3 M ) \}
}{
K^{ [ \bfL^{-1} ( \lVert \bfl \rVert_{ \infty } + 1 )^{-2} ] }
}
\\ & \quad
+
\frac{
    18 \max\{ 1, ( v - u )^2 \}
    b \bfL ( \lVert \bfl \rVert_{ \infty } + 1 )^2 B
    \max\{ p, \ln( 3 M ) \}
}{ \sqrt{M} }
\\ &
\leq
\frac{
    20 \max\{ 1, ( v - u )^2 \}
    b \bfL ( \lVert \bfl \rVert_{ \infty } + 1 )^{ \bfL + 1 } B^\bfL
    \max\{ p, \ln( 3 M ) \}
}{
\min\{
    \sqrt{M},
    K^{ [ \bfL^{-1} ( \lVert \bfl \rVert_{ \infty } + 1 )^{-2} ] }
\} }
.
\end{split}
\end{equation}
This and~\cref{eq:prop:main_derivation}
complete
\end{cproof2}

\begin{lemma}
\label{lem:ln_estimate}
Let
$ a, x, p \in ( 0, \infty ) $,
$ M, c \in [ 1, \infty ) $,
$ B \in [ c, \infty ) $.
Then
\begin{enumerate}[(i)]
\item
\label{item:lem:ln_estimate:1}
it holds that
$ a x^p
\leq
\exp\bigl( a^{ \nicefrac{1}{p} } \tfrac{ p x }{ e } \bigr) $
and
\item
\label{item:lem:ln_estimate:2}
it holds that
$ \ln( 3 M B c )
\leq
\tfrac{23 B}{18}
\ln( e M ) $.
\end{enumerate}
\end{lemma}
\begin{cproof}{lem:ln_estimate}
First,
note that
the fact that
$ \forall \, y \in \R \colon
y + 1 \leq e^{ y } $
demonstrates that%
\begin{equation}
a x^p
=
( a^{ \nicefrac{1}{p} } x )^p
=
\bigl[ e \bigl( a^{ \nicefrac{1}{p} } \tfrac{ x }{ e } - 1 + 1 \bigr) \bigr]^p
\leq
\bigl[ e \exp\bigl( a^{ \nicefrac{1}{p} } \tfrac{ x }{ e }  - 1 \bigr) \bigr]^p
=
\exp\bigl( a^{ \nicefrac{1}{p} } \tfrac{ p x }{ e } \bigr)
.
\end{equation}
This proves~\cref{item:lem:ln_estimate:1}.

Second,
observe that
\cref{item:lem:ln_estimate:1}
and
the fact that
$ \nicefrac{ 2 \sqrt{3} }{e}
\leq 
\nicefrac{23}{18} $
ensure that
\begin{equation}
3 B^2
\leq
\exp\bigl( \sqrt{3} \tfrac{ 2 B }{e} \bigr)
=
\exp\bigl( \tfrac{ 2 \sqrt{3} B }{e} \bigr)
\leq
\exp\bigl( \tfrac{ 23 B }{18} \bigr)
.
\end{equation}
The facts that
$ B \geq c \geq 1 $
and
$ M \geq 1 $
hence
imply that
\begin{equation}
\ln( 3 M B c )
\leq
\ln( 3 B^2 M )
\leq
\ln( [ e M ]^{ \nicefrac{23 B}{18} } )
=
\tfrac{23 B}{18}
\ln( e M )
.
\end{equation}
This establishes~\cref{item:lem:ln_estimate:2}.
\end{cproof}

\begin{theorem}
\label{thm:main}
Let
$ d, \bfd, \bfL, M, K, N \in \N $,
$ A \in ( 0, \infty ) $,
$ L, a, u \in \R $,
$ b \in ( a, \infty ) $,
$ v \in ( u, \infty ) $,
$ c \in
[ \max\{ 1, L, \lvert a \rvert, \lvert b \rvert, 2 \lvert u \rvert, 2 \lvert v \rvert \}, \infty ) $,
$ B \in [ c, \infty ) $,
$ \bfl = ( \bfl_0, \bfl_1, \ldots, \bfl_\bfL ) \in \N^{ \bfL + 1 } $,
$ \bfN \subseteq \{ 0, 1, \ldots, N \} $,
assume
$ 0 \in \bfN $,
$ \bfL \geq \nicefrac{ A \mathbbm{1}_{ \smash{ ( 6^d, \infty ) } }( A ) }{ ( 2d ) } + 1 $,
$ \bfl_0 = d $,
$ \bfl_1 \geq A \mathbbm{1}_{ \smash{ ( 6^d, \infty ) } }( A ) $,
$ \bfl_\bfL = 1 $,
and
$ \bfd \geq \sum_{i=1}^{\bfL} \bfl_i( \bfl_{ i - 1 } + 1 ) $,
assume for all
$ i \in \{ 2, 3, \ldots \} \cap [ 0, \bfL ) $
that
$ \bfl_i \geq \mathbbm{1}_{ \smash{ ( 6^d, \infty ) } }( A ) \max\{ \nicefrac{A}{d} - 2i + 3, 2 \} $,
let
$ ( \Omega, \cF, \P ) $
be a probability space,
let
$ X_j \colon \Omega \to [ a, b ]^d $,
$ j \in \N $,
and
$ Y_j \colon \Omega \to [ u, v ] $,
$ j \in \N $,
be functions,
assume that
$ ( X_j, Y_j ) $,
$ j \in \{ 1, 2, \ldots, M \} $,
are i.i.d.\ random variables,
let
$ \cE \colon [ a, b ]^d \to [ u, v ] $
satisfy
$ \P $-a.s.\
that
$ \cE( X_1 )
= \E[ Y_1 \vert X_1 ] $,
assume for all
$ x, y \in [ a, b ]^d $
that
$ \lvert \cE( x ) - \cE( y ) \rvert \leq L \lVert x - y \rVert_{ 1 } $,
let
$ \Theta_{ k, n } \colon \Omega \to \R^{\bfd} $,
$ k, n \in \N_0 $,
and
$ \bfk \colon \Omega \to ( \N_0 )^2 $
be random variables,
assume
$ \bigl( \bigcup_{ k = 1 }^{ \infty }
\Theta_{ k, 0 }( \Omega ) \bigr)
\subseteq [ -B, B ]^\bfd $,
assume that
$ \Theta_{ k, 0 } $,
$ k \in \{ 1, 2, \ldots, K \} $,
are i.i.d.,
assume that
$ \Theta_{ 1, 0 } $ is continuous uniformly distributed on $ [ -c, c ]^\bfd $,
and let
$ \cR \colon \R^{\bfd} \times \Omega \to [ 0, \infty ) $
satisfy for all
$ \theta \in \R^{\bfd} $,
$ \omega \in \Omega $
that
\begin{align}
\cR( \theta, \omega )
=
\frac{1}{M}
\biggl[
\smallsum_{j=1}^M
    \lvert \clippedNN{\theta}{\bfl}{u}{v}( X_j( \omega ) ) - Y_j( \omega ) \rvert^2
\biggr]
\qquad
\text{and}
\\
\bfk( \omega ) \in
\argmin\nolimits_{ ( k, n ) \in \{ 1, 2, \ldots, K \} \times \bfN, \, \lVert \Theta_{ k, n }( \omega ) \rVert_{ \infty } \leq B }
\cR( \Theta_{ k, n }( \omega ), \omega )
\end{align}
(cf.~\cref{def:clipped_NN,def:p-norm}).
Then
it holds for all
$ p \in ( 0, \infty ) $
that
\begin{equation}
\label{eq:thm:main}
\begin{split}
&
\Bigl(
\E\Bigl[
\Bigl( \medint{[ a, b ]^d}
        \lvert \clippedNN{\Theta_{\bfk}}{\bfl}{u}{v}( x ) - \cE( x ) \rvert^2
    \, \P_{ X_1 }( \dd x )
\Bigr)^{ \! p \, }
\Bigr]
\Bigr)^{ \! \nicefrac{1}{p} }
\\ &
\leq
\frac{ 9 d^2 L^2 ( b - a )^2 }{ A^{ \nicefrac{2}{d} } }
+
\frac{
4 ( v - u )
\bfL
( \lVert \bfl \rVert_{ \infty } + 1 )^\bfL
c^{ \bfL + 1 }
\max\{ 1, p \}
}{
K^{ [ \bfL^{-1} ( \lVert \bfl \rVert_{ \infty } + 1 )^{-2} ] }
}
\\ & \quad
+
\frac{
    18 \max\{ 1, ( v - u )^2 \}
    \bfL ( \lVert \bfl \rVert_{ \infty } + 1 )^2
    \max\{
        p,
        \ln( 3 M B c )
    \}
}{ \sqrt{M} }
\\ &
\leq
\frac{ 36 d^2 c^4 }{ A^{ \nicefrac{2}{d} } }
+
\frac{
4
\bfL
( \lVert \bfl \rVert_{ \infty } + 1 )^\bfL
c^{ \bfL + 2 }
\max\{ 1, p \}
}{
K^{ [ \bfL^{-1} ( \lVert \bfl \rVert_{ \infty } + 1 )^{-2} ] }
}
+
\frac{
    23 B^3
    \bfL ( \lVert \bfl \rVert_{ \infty } + 1 )^2
    \max\{ p, \ln( e M ) \}
}{ \sqrt{M} }
\end{split}
\end{equation}
(cf.~\cref{item:lem:measurability:3} in \cref{lem:measurability}).
\end{theorem}
\begin{cproof}{thm:main}
First of all,
note that
the assumption that
$ \forall \, x, y \in [ a, b ]^d \colon
\lvert \cE( x ) - \cE( y ) \rvert \leq L \lVert x - y \rVert_{ 1 } $
ensures that
$ \cE \colon [ a, b ]^d \to [ u, v ] $
is a
$ \cB( [ a, b ]^d ) $/$ \cB( [ u, v ] ) $-measurable function.
The fact that
$ \max\{ 1, \lvert a \rvert, \lvert b \rvert \}
\leq c $
and
\cref{prop:main}
(with
$ b \leftarrow \max\{ 1, \lvert a \rvert, \lvert b \rvert \} $,
$ D \leftarrow [ a, b ]^d $
in the notation of \cref{prop:main})
hence
show for all
$ p \in ( 0, \infty ) $
that
\begin{align*}
\label{eq:prop:main_estimate}
&
\Bigl(
\E\Bigl[
\Bigl( \medint{[ a, b ]^d}
        \lvert \clippedNN{\Theta_{\bfk}}{\bfl}{u}{v}( x ) - \cE( x ) \rvert^2
    \, \P_{ X_1 }( \dd x )
\Bigr)^{ \! p \, }
\Bigr]
\Bigr)^{ \! \nicefrac{1}{p} }
\\ &
\leq
\bigl[
\inf\nolimits_{ \theta \in [ -c, c ]^\bfd }
\sup\nolimits_{ x \in [ a, b ]^d }
    \lvert \clippedNN{\theta}{\bfl}{u}{v}( x ) - \cE( x ) \rvert^{2}
\bigr]
\\ & \quad
+
\frac{
4 ( v - u ) \max\{ 1, \lvert a \rvert, \lvert b \rvert \}
\bfL
( \lVert \bfl \rVert_{ \infty } + 1 )^\bfL
c^{ \bfL }
\max\{ 1, p \}
}{
K^{ [ \bfL^{-1} ( \lVert \bfl \rVert_{ \infty } + 1 )^{-2} ] }
}
\\ & \quad \yesnumber
+
\frac{
    18 \max\{ 1, ( v - u )^2 \}
    \bfL ( \lVert \bfl \rVert_{ \infty } + 1 )^2
    \max\{
        p,
        \ln( 3 M B \max\{ 1, \lvert a \rvert, \lvert b \rvert \} )
    \}
}{ \sqrt{M} }
\\ &
\leq
\bigl[
\inf\nolimits_{ \theta \in [ -c, c ]^\bfd }
\sup\nolimits_{ x \in [ a, b ]^d }
    \lvert \clippedNN{\theta}{\bfl}{u}{v}( x ) - \cE( x ) \rvert^{2}
\bigr]
+
\frac{
4 ( v - u )
\bfL
( \lVert \bfl \rVert_{ \infty } + 1 )^\bfL
c^{ \bfL + 1 }
\max\{ 1, p \}
}{
K^{ [ \bfL^{-1} ( \lVert \bfl \rVert_{ \infty } + 1 )^{-2} ] }
}
\\ & \quad
+
\frac{
    18 \max\{ 1, ( v - u )^2 \}
    \bfL ( \lVert \bfl \rVert_{ \infty } + 1 )^2
    \max\{
        p,
        \ln( 3 M B c )
    \}
}{ \sqrt{M} }
.
\end{align*}
Furthermore,
observe that
\cref{prop:approximation_error}
(with
$ f \leftarrow \cE $
in the notation of \cref{prop:approximation_error})
proves that
there exists $ \vartheta \in \R^\bfd $
such that
$ \lVert \vartheta \rVert_\infty
\leq \max\{ 1, L, \lvert a \rvert, \lvert b \rvert, 2[ \sup_{ x \in [ a, b ]^d } \lvert \cE( x ) \rvert ] \} $
and
\begin{equation}
\label{eq:prop:approximation_error}
\sup\nolimits_{ x \in [ a, b ]^d }
    \lvert \clippedNN{\vartheta}{\bfl}{u}{v}( x ) - \cE( x ) \rvert
\leq
\frac{ 3 d L ( b - a ) }{ A^{ \nicefrac{1}{d} } }
.
\end{equation}
The fact that
$ \forall \, x \in [ a, b ]^d \colon
\cE( x ) \in [ u, v ] $
hence
implies that
\begin{equation}
\lVert \vartheta \rVert_\infty
\leq
\max\{ 1, L, \lvert a \rvert, \lvert b \rvert, 2 \lvert u \rvert, 2 \lvert v \rvert \}
\leq c
.
\end{equation}
This and
\cref{eq:prop:approximation_error}
demonstrate that
\begin{equation}
\begin{split}
& \inf\nolimits_{ \theta \in [ -c, c ]^\bfd }
\sup\nolimits_{ x \in [ a, b ]^d }
    \lvert \clippedNN{\theta}{\bfl}{u}{v}( x ) - \cE( x ) \rvert^2
\\ &
\leq
\sup\nolimits_{ x \in [ a, b ]^d }
    \lvert \clippedNN{\vartheta}{\bfl}{u}{v}( x ) - \cE( x ) \rvert^2
\\ &
\leq
\biggl[
\frac{ 3 d L ( b - a ) }{ A^{ \nicefrac{1}{d} } }
\biggr]^2
=
\frac{ 9 d^2 L^2 ( b - a )^2 }{ A^{ \nicefrac{2}{d} } }
.
\end{split}
\end{equation}
Combining
this
with
\cref{eq:prop:main_estimate}
establishes for all
$ p \in ( 0, \infty ) $
that
\begin{equation}
\label{eq:overall_estimate}
\begin{split}
&
\Bigl(
\E\Bigl[
\Bigl( \medint{[ a, b ]^d}
        \lvert \clippedNN{\Theta_{\bfk}}{\bfl}{u}{v}( x ) - \cE( x ) \rvert^2
    \, \P_{ X_1 }( \dd x )
\Bigr)^{ \! p \, }
\Bigr]
\Bigr)^{ \! \nicefrac{1}{p} }
\\ &
\leq
\frac{ 9 d^2 L^2 ( b - a )^2 }{ A^{ \nicefrac{2}{d} } }
+
\frac{
4 ( v - u )
\bfL
( \lVert \bfl \rVert_{ \infty } + 1 )^\bfL
c^{ \bfL + 1 }
\max\{ 1, p \}
}{
K^{ [ \bfL^{-1} ( \lVert \bfl \rVert_{ \infty } + 1 )^{-2} ] }
}
\\ & \quad
+
\frac{
    18 \max\{ 1, ( v - u )^2 \}
    \bfL ( \lVert \bfl \rVert_{ \infty } + 1 )^2
    \max\{
        p,
        \ln( 3 M B c )
    \}
}{ \sqrt{M} }
.
\end{split}
\end{equation}
Moreover,
note that
the facts that
$ \max\{ 1, L, \lvert a \rvert, \lvert b \rvert \} \leq c $
and
$ ( b - a )^2 \leq ( \lvert a \rvert + \lvert b \rvert )^2 \leq 2 ( a^2 + b^2 ) $
yield that
\begin{equation}
\label{eq:coarse_estimate}
9 L^2 ( b - a )^2
\leq
18 c^2 ( a^2 + b^2 )
\leq
18 c^2 ( c^2 + c^2 )
=
36 c^4
.
\end{equation}
In addition,
the fact that
$ B \geq c \geq 1 $,
the fact that
$ M \geq 1 $,
and~\cref{item:lem:ln_estimate:2}
in
\cref{lem:ln_estimate}
ensure that
$ \ln( 3 M B c )
\leq
\tfrac{ 23 B }{18}
\ln( e M ) $.
This,
\cref{eq:coarse_estimate},
the fact that
$ ( v - u )
\leq 2 \max\{ \lvert u \rvert, \lvert v \rvert \}
= \max\{ 2 \lvert u \rvert, 2 \lvert v \rvert \}
\leq c
\leq B $,
and
the fact that
$ B \geq 1 $
prove for all
$ p \in ( 0, \infty ) $
that
\begin{equation}
\begin{split}
&
\frac{ 9 d^2 L^2 ( b - a )^2 }{ A^{ \nicefrac{2}{d} } }
+
\frac{
4 ( v - u )
\bfL
( \lVert \bfl \rVert_{ \infty } + 1 )^\bfL
c^{ \bfL + 1 }
\max\{ 1, p \}
}{
K^{ [ \bfL^{-1} ( \lVert \bfl \rVert_{ \infty } + 1 )^{-2} ] }
}
\\ &
+
\frac{
    18 \max\{ 1, ( v - u )^2 \}
    \bfL ( \lVert \bfl \rVert_{ \infty } + 1 )^2
    \max\{
        p,
        \ln( 3 M B c )
    \}
}{ \sqrt{M} }
\\ &
\leq
\frac{ 36 d^2 c^4 }{ A^{ \nicefrac{2}{d} } }
+
\frac{
4
\bfL
( \lVert \bfl \rVert_{ \infty } + 1 )^\bfL
c^{ \bfL + 2 }
\max\{ 1, p \}
}{
K^{ [ \bfL^{-1} ( \lVert \bfl \rVert_{ \infty } + 1 )^{-2} ] }
}
+
\frac{
    23 B^3
    \bfL ( \lVert \bfl \rVert_{ \infty } + 1 )^2
    \max\{ p, \ln( e M ) \}
}{ \sqrt{M} }
.
\end{split}
\end{equation}
Combining this with~\cref{eq:overall_estimate}
shows~\cref{eq:thm:main}.
\end{cproof}

\begin{corollary}
\label{cor:min_NN_architecture}
Let
$ d, \bfd, \bfL, M, K, N \in \N $,
$ L, a, u \in \R $,
$ b \in ( a, \infty ) $,
$ v \in ( u, \infty ) $,
$ c \in
[ \max\{ 1, L, \lvert a \rvert, \lvert b \rvert, 2 \lvert u \rvert, 2 \lvert v \rvert \}, \infty ) $,
$ B \in [ c, \infty ) $,
$ \bfl = ( \bfl_0, \bfl_1, \ldots, \bfl_\bfL ) \in \N^{ \bfL + 1 } $,
$ \bfN \subseteq \{ 0, 1, \ldots,
\allowbreak
N \} $,
assume
$ 0 \in \bfN $,
$ \bfl_0 = d $,
$ \bfl_\bfL = 1 $,
and
$ \bfd \geq \sum_{i=1}^{\bfL} \bfl_i( \bfl_{ i - 1 } + 1 ) $,
let
$ ( \Omega, \cF, \P ) $
be a probability space,
let
$ X_j \colon \Omega \to [ a, b ]^d $,
$ j \in \N $,
and
$ Y_j \colon \Omega \to [ u, v ] $,
$ j \in \N $,
be functions,
assume that
$ ( X_j, Y_j ) $,
$ j \in \{ 1, 2, \ldots, M \} $,
are i.i.d.\ random variables,
let
$ \cE \colon [ a, b ]^d \to [ u, v ] $
satisfy
$ \P $-a.s.\
that
$ \cE( X_1 )
= \E[ Y_1 \vert X_1 ] $,
assume for all
$ x, y \in [ a, b ]^d $
that
$ \lvert \cE( x ) - \cE( y ) \rvert \leq L \lVert x - y \rVert_{ 1 } $,
let
$ \Theta_{ k, n } \colon \Omega \to \R^{\bfd} $,
$ k, n \in \N_0 $,
and
$ \bfk \colon \Omega \to ( \N_0 )^2 $
be random variables,
assume
$ \bigl( \bigcup_{ k = 1 }^{ \infty }
\Theta_{ k, 0 }( \Omega ) \bigr)
\subseteq [ -B, B ]^\bfd $,
assume that
$ \Theta_{ k, 0 } $,
$ k \in \{ 1, 2, \ldots, K \} $,
are i.i.d.,
assume that
$ \Theta_{ 1, 0 } $ is continuous uniformly distributed on $ [ -c, c ]^\bfd $,
and let
$ \cR \colon \R^{\bfd} \times \Omega \to [ 0, \infty ) $
satisfy for all
$ \theta \in \R^{\bfd} $,
$ \omega \in \Omega $
that
\begin{align}
\cR( \theta, \omega )
=
\frac{1}{M}
\biggl[
\smallsum_{j=1}^M
    \lvert \clippedNN{\theta}{\bfl}{u}{v}( X_j( \omega ) ) - Y_j( \omega ) \rvert^2
\biggr]
\qquad
\text{and}
\\
\bfk( \omega ) \in
\argmin\nolimits_{ ( k, n ) \in \{ 1, 2, \ldots, K \} \times \bfN, \, \lVert \Theta_{ k, n }( \omega ) \rVert_{ \infty } \leq B }
\cR( \Theta_{ k, n }( \omega ), \omega )
\end{align}
(cf.~\cref{def:clipped_NN,def:p-norm}).
Then
it holds for all
$ p \in ( 0, \infty ) $
that
\begin{align*}
&
\Bigl(
\E\Bigl[
\Bigl( \medint{[ a, b ]^d}
        \lvert \clippedNN{\Theta_{\bfk}}{\bfl}{u}{v}( x ) - \cE( x ) \rvert^2
    \, \P_{ X_1 }( \dd x )
\Bigr)^{ \! \nicefrac{p}{2} \, }
\Bigr]
\Bigr)^{ \! \nicefrac{1}{p} }
\\ &
\leq
\frac{
3 d L ( b - a )
}{
[ \min(
    \{ \bfL \}
    \cup
    \{ \bfl_i \colon i \in \N \cap [ 0, \bfL ) \}
) ]^{ \nicefrac{1}{d} }
}
+
\frac{
2
[ ( v - u )
\bfL
( \lVert \bfl \rVert_{ \infty } + 1 )^\bfL
c^{ \bfL + 1 }
\max\{ 1, \nicefrac{p}{2} \}
]^{ \nicefrac{1}{2} }
}{
K^{ [ ( 2 \bfL )^{-1} ( \lVert \bfl \rVert_{ \infty } + 1 )^{-2} ] }
}
\\ & \quad \yesnumber
+
\frac{
    3 \max\{ 1, v - u \}
    ( \lVert \bfl \rVert_{ \infty } + 1 )
    [ \bfL
    \max\{
        p,
        2 \ln( 3 M B c )
    \}
    ]^{ \nicefrac{1}{2} }
}{ M^{ \nicefrac{1}{4} } }
\\ &
\leq
\frac{
6 d c^2
}{
[ \min(
    \{ \bfL \}
    \cup
    \{ \bfl_i \colon i \in \N \cap [ 0, \bfL ) \}
) ]^{ \nicefrac{1}{d} }
}
+
\frac{
2
\bfL
( \lVert \bfl \rVert_{ \infty } + 1 )^\bfL
c^{ \bfL + 1 }
\max\{ 1, p \}
}{
K^{ [ ( 2 \bfL )^{-1} ( \lVert \bfl \rVert_{ \infty } + 1 )^{-2} ] }
}
\\ & \quad
+
\frac{
    5 B^2
    \bfL ( \lVert \bfl \rVert_{ \infty } + 1 )
    \max\{ p, \ln( e M ) \}
}{ M^{ \nicefrac{1}{4} } }
\end{align*}
(cf.~\cref{item:lem:measurability:3} in \cref{lem:measurability}).
\end{corollary}
\begin{cproof}{cor:min_NN_architecture}
Throughout this proof
let
$ A \in ( 0, \infty ) $
be given by
\begin{equation}
\label{eq:choice_A}
A
=
\min(
    \{ \bfL \}
    \cup
    \{ \bfl_i \colon i \in \N \cap [ 0, \bfL ) \}
)
.
\end{equation}
Note that
\cref{eq:choice_A}
ensures that
\begin{equation}
\label{eq:assumption_bfL}
\begin{split}
\bfL
& \geq
A
=
A - 1 + 1
\geq
( A - 1 ) \mathbbm{1}_{ [ 2, \infty ) }( A ) + 1
\\ &
\geq
\bigl( A - \tfrac{A}{2} \bigr) \mathbbm{1}_{ [ 2, \infty ) }( A ) + 1
=
\tfrac{ A \mathbbm{1}_{ [ 2, \infty ) }( A ) }{ 2 } + 1
\geq
\tfrac{ A \mathbbm{1}_{ \smash{ ( 6^d, \infty ) } }( A ) }{ 2d } + 1
.
\end{split}
\end{equation}
Moreover,
the assumption that
$ \bfl_\bfL = 1 $
and
\cref{eq:choice_A}
imply that
\begin{equation}
\label{eq:assumption_bfl_1}
\bfl_1
=
\bfl_1
\mathbbm{1}_{ \{ 1 \} }( \bfL )
+
\bfl_1
\mathbbm{1}_{ [ 2, \infty ) }( \bfL )
\geq
\mathbbm{1}_{ \{ 1 \} }( \bfL )
+
A
\mathbbm{1}_{ [ 2, \infty ) }( \bfL )
=
A
\geq
A \mathbbm{1}_{ \smash{ ( 6^d, \infty ) } }( A )
.
\end{equation}
Moreover,
again~\cref{eq:choice_A}
shows for all
$ i \in \{ 2, 3, \ldots \} \cap [ 0, \bfL ) $
that
\begin{equation}
\label{eq:assumption_bfl_i}
\begin{split}
\bfl_i
& \geq
A
\geq
A
\mathbbm{1}_{ \smash{ [ 2, \infty ) } }( A )
\geq
\mathbbm{1}_{ \smash{ [ 2, \infty ) } }( A )
\max\{ A - 1, 2 \}
=
\mathbbm{1}_{ \smash{ [ 2, \infty ) } }( A )
\max\{ A - 4 + 3, 2 \}
\\ &
\geq
\mathbbm{1}_{ \smash{ [ 2, \infty ) } }( A )
\max\{ A - 2 i + 3, 2 \}
\geq
\mathbbm{1}_{ \smash{ ( 6^d, \infty ) } }( A )
\max\{ \nicefrac{A}{d} - 2i + 3, 2 \}
.
\end{split}
\end{equation}
Combining
\cref{eq:assumption_bfL}--\cref{eq:assumption_bfl_i}
and
\cref{thm:main}
(with
$ p \leftarrow \nicefrac{p}{2} $
for
$ p \in ( 0, \infty ) $
in the notation of \cref{thm:main})
establishes for all
$ p \in ( 0, \infty ) $
that
\begin{align*}
&
\Bigl(
\E\Bigl[
\Bigl( \medint{[ a, b ]^d}
        \lvert \clippedNN{\Theta_{\bfk}}{\bfl}{u}{v}( x ) - \cE( x ) \rvert^2
    \, \P_{ X_1 }( \dd x )
\Bigr)^{ \! \nicefrac{p}{2} \, }
\Bigr]
\Bigr)^{ \! \nicefrac{2}{p} }
\\ &
\leq
\frac{ 9 d^2 L^2 ( b - a )^2 }{ A^{ \nicefrac{2}{d} } }
+
\frac{
4 ( v - u )
\bfL
( \lVert \bfl \rVert_{ \infty } + 1 )^\bfL
c^{ \bfL + 1 }
\max\{ 1, \nicefrac{p}{2} \}
}{
K^{ [ \bfL^{-1} ( \lVert \bfl \rVert_{ \infty } + 1 )^{-2} ] }
}
\\ & \quad \yesnumber
+
\frac{
    18 \max\{ 1, ( v - u )^2 \}
    \bfL ( \lVert \bfl \rVert_{ \infty } + 1 )^2
    \max\{
        \nicefrac{p}{2},
        \ln( 3 M B c )
    \}
}{ \sqrt{M} }
\\ &
\leq
\frac{ 36 d^2 c^4 }{ A^{ \nicefrac{2}{d} } }
+
\frac{
4
\bfL
( \lVert \bfl \rVert_{ \infty } + 1 )^\bfL
c^{ \bfL + 2 }
\max\{ 1, \nicefrac{p}{2} \}
}{
K^{ [ \bfL^{-1} ( \lVert \bfl \rVert_{ \infty } + 1 )^{-2} ] }
}
+
\frac{
    23 B^3
    \bfL ( \lVert \bfl \rVert_{ \infty } + 1 )^2
    \max\{ \nicefrac{p}{2}, \ln( e M ) \}
}{ \sqrt{M} }
.
\end{align*}
This,
\cref{eq:choice_A},
and
the facts that
$ \bfL \geq 1 $,
$ c \geq 1 $,
$ B \geq 1 $,
and
$ \ln( e M ) \geq 1 $
demonstrate for all
$ p \in ( 0, \infty ) $
that
\begin{align*}
&
\Bigl(
\E\Bigl[
\Bigl( \medint{[ a, b ]^d}
        \lvert \clippedNN{\Theta_{\bfk}}{\bfl}{u}{v}( x ) - \cE( x ) \rvert^2
    \, \P_{ X_1 }( \dd x )
\Bigr)^{ \! \nicefrac{p}{2} \, }
\Bigr]
\Bigr)^{ \! \nicefrac{1}{p} }
\\ &
\leq
\frac{
3 d L ( b - a )
}{
[ \min(
    \{ \bfL \}
    \cup
    \{ \bfl_i \colon i \in \N \cap [ 0, \bfL ) \}
) ]^{ \nicefrac{1}{d} }
}
+
\frac{
2
[ ( v - u )
\bfL
( \lVert \bfl \rVert_{ \infty } + 1 )^\bfL
c^{ \bfL + 1 }
\max\{ 1, \nicefrac{p}{2} \}
]^{ \nicefrac{1}{2} }
}{
K^{ [ ( 2 \bfL )^{-1} ( \lVert \bfl \rVert_{ \infty } + 1 )^{-2} ] }
}
\\ & \quad
+
\frac{
    3 \max\{ 1, v - u \}
    ( \lVert \bfl \rVert_{ \infty } + 1 )
    [ \bfL
    \max\{
        p,
        2 \ln( 3 M B c )
    \}
    ]^{ \nicefrac{1}{2} }
}{ M^{ \nicefrac{1}{4} } }
\\ & \yesnumber
\leq
\frac{
6 d c^2
}{
[ \min(
    \{ \bfL \}
    \cup
    \{ \bfl_i \colon i \in \N \cap [ 0, \bfL ) \}
) ]^{ \nicefrac{1}{d} }
}
+
\frac{
2
[
\bfL
( \lVert \bfl \rVert_{ \infty } + 1 )^\bfL
c^{ \bfL + 2 }
\max\{ 1, \nicefrac{p}{2} \}
]^{ \nicefrac{1}{2} }
}{
K^{ [ ( 2 \bfL )^{-1} ( \lVert \bfl \rVert_{ \infty } + 1 )^{-2} ] }
}
\\ & \quad
+
\frac{
    5 B^3
    [ \bfL ( \lVert \bfl \rVert_{ \infty } + 1 )^2
    \max\{ \nicefrac{p}{2}, \ln( e M ) \}
    ]^{ \nicefrac{1}{2} }
}{ M^{ \nicefrac{1}{4} } }
\\ &
\leq
\frac{
6 d c^2
}{
[ \min(
    \{ \bfL \}
    \cup
    \{ \bfl_i \colon i \in \N \cap [ 0, \bfL ) \}
) ]^{ \nicefrac{1}{d} }
}
+
\frac{
2
\bfL
( \lVert \bfl \rVert_{ \infty } + 1 )^\bfL
c^{ \bfL + 1 }
\max\{ 1, p \}
}{
K^{ [ ( 2 \bfL )^{-1} ( \lVert \bfl \rVert_{ \infty } + 1 )^{-2} ] }
}
\\ & \quad
+
\frac{
    5 B^2
    \bfL ( \lVert \bfl \rVert_{ \infty } + 1 )
    \max\{ p, \ln( e M ) \}
}{ M^{ \nicefrac{1}{4} } }
.
\end{align*}
\end{cproof}

\subsection[Full strong error analysis with optimisation via stochastic gradient descent]{Full strong error analysis for the training of DNNs with optimisation via stochastic gradient descent with random initialisation}
\label{sec:SGD}

\begin{corollary}
\label{cor:SGD_full_error}
Let
$ d, \bfd, \bfL, M, K, N \in \N $,
$ L, a, u \in \R $,
$ b \in ( a, \infty ) $,
$ v \in ( u, \infty ) $,
$ c \in
[ \max\{ 1, L, \lvert a \rvert, \lvert b \rvert, 2 \lvert u \rvert, 2 \lvert v \rvert \}, \infty ) $,
$ B \in [ c, \infty ) $,
$ \bfl = ( \bfl_0, \bfl_1, \ldots, \bfl_\bfL ) \in \N^{ \bfL + 1 } $,
$ \bfN \subseteq \{ 0, 1, \ldots,
\allowbreak
N \} $,
$ ( \bfJ_n )_{ n \in \N } \subseteq \N $,
$ ( \gamma_n )_{ n \in \N } \subseteq \R $,
assume
$ 0 \in \bfN $,
$ \bfl_0 = d $,
$ \bfl_\bfL = 1 $,
and
$ \bfd \geq \sum_{i=1}^{\bfL} \bfl_i( \bfl_{ i - 1 } + 1 ) $,
let
$ ( \Omega, \cF, \P ) $
be a probability space,
let
$ X^{ k, n }_{ \smash{j} }
\colon \Omega \to [ a, b ]^d $,
$ k, n, j \in \N_0 $,
and
$ Y^{ k, n }_{ \smash{j} }
\colon \Omega \to [ u, v ] $,
$ k, n, j \in \N_0 $,
be functions,
assume that
$ ( X^{ 0, 0 }_{ \smash{j} }, Y^{ 0, 0 }_{ \smash{j} } ) $,
$ j \in \{ 1, 2, \ldots, M \} $,
are i.i.d.\ random variables,
let
$ \cE \colon [ a, b ]^d \to [ u, v ] $
satisfy
$ \P $-a.s.\
that
$ \cE( X_{ \smash{1} }^{ 0, 0 } )
= \E[ Y_{ \smash{1} }^{ 0, 0 } \vert X_{ \smash{1} }^{ 0, 0 } ] $,
assume for all
$ x, y \in [ a, b ]^d $
that
$ \lvert \cE( x ) - \cE( y ) \rvert \leq L \lVert x - y \rVert_{ 1 } $,
let
$ \Theta_{ k, n } \colon \Omega \to \R^{\bfd} $,
$ k, n \in \N_0 $,
and
$ \bfk \colon \Omega \to ( \N_0 )^2 $
be random variables,
assume
$ \bigl( \bigcup_{ k = 1 }^{ \infty }
\Theta_{ k, 0 }( \Omega ) \bigr)
\subseteq [ -B, B ]^\bfd $,
assume that
$ \Theta_{ k, 0 } $,
$ k \in \{ 1, 2, \ldots, K \} $,
are i.i.d.,
assume that
$ \Theta_{ 1, 0 } $ is continuous uniformly distributed on $ [ -c, c ]^\bfd $,
let
$ \cR^{ k, n }_{ \smash{J} } \colon \R^{\bfd} \times \Omega \to [ 0, \infty ) $,
$ k, n, J \in \N_0 $,
and
$ \cG^{ k, n } \colon \R^{\bfd} \times \Omega \to \R^{\bfd} $,
$ k, n \in \N $,
satisfy for all
$ k, n \in \N $,
$ \omega \in \Omega $,
$ \theta \in
\{ \vartheta \in \R^{\bfd} \colon
( \cR^{ k, n }_{ \smash{ \bfJ_n } } ( \cdot, \omega ) \colon \R^{\bfd} \to [ 0, \infty )
\text{ is differentiable at } \vartheta ) \} $
that
$ \cG^{ k, n }( \theta, \omega )
=
( \nabla_\theta \cR^{ k, n }_{ \smash{ \bfJ_n } } )( \theta, \omega ) $,
assume for all
$ k, n \in \N $
that
$ \Theta_{ k, n } = \Theta_{ k, n -1 } - \gamma_n \cG^{ k, n }( \Theta_{ k, n -1 } ) $,
and
assume for all
$ k, n \in \N_0 $,
$ J \in \N $,
$ \theta \in \R^{\bfd} $,
$ \omega \in \Omega $
that
\begin{align}
\cR^{ k, n }_{ \smash{J} }( \theta, \omega )
=
\frac{1}{J}
\biggl[
\smallsum_{j=1}^{J}
    \lvert \clippedNN{\theta}{\bfl}{u}{v}( X^{ k, n }_{ \smash{j} }( \omega ) ) - Y^{ k, n }_{ \smash{j} }( \omega ) \rvert^2
\biggr]
\qquad
\text{and}
\\
\bfk( \omega ) \in
\argmin\nolimits_{ ( l, m ) \in \{ 1, 2, \ldots, K \} \times \bfN, \, \lVert \Theta_{ l, m }( \omega ) \rVert_{ \infty } \leq B }
\cR^{ 0, 0 }_{ \smash{ M } }( \Theta_{ l, m }( \omega ), \omega )
\end{align}
(cf.~\cref{def:clipped_NN,def:p-norm}).
Then
it holds for all
$ p \in ( 0, \infty ) $
that
\begin{align*}
\label{eq:cor:SGD_full_error}
&
\Bigl(
\E\Bigl[
\Bigl( \medint{[ a, b ]^d}
        \lvert \clippedNN{\Theta_{\bfk}}{\bfl}{u}{v}( x ) - \cE( x ) \rvert^2
    \, \P_{ X^{ 0, 0 }_{ \smash{1}\vphantom{x} } }( \dd x )
\Bigr)^{ \! \nicefrac{p}{2} \, }
\Bigr]
\Bigr)^{ \! \nicefrac{1}{p} }
\\ &
\leq
\frac{
3 d L ( b - a )
}{
[ \min(
    \{ \bfL \}
    \cup
    \{ \bfl_i \colon i \in \N \cap [ 0, \bfL ) \}
) ]^{ \nicefrac{1}{d} }
}
+
\frac{
2
[ ( v - u )
\bfL
( \lVert \bfl \rVert_{ \infty } + 1 )^\bfL
c^{ \bfL + 1 }
\max\{ 1, \nicefrac{p}{2} \}
]^{ \nicefrac{1}{2} }
}{
K^{ [ ( 2 \bfL )^{-1} ( \lVert \bfl \rVert_{ \infty } + 1 )^{-2} ] }
}
\\ & \quad \yesnumber
+
\frac{
    3 \max\{ 1, v - u \}
    ( \lVert \bfl \rVert_{ \infty } + 1 )
    [ \bfL
    \max\{
        p,
        2 \ln( 3 M B c )
    \}
    ]^{ \nicefrac{1}{2} }
}{ M^{ \nicefrac{1}{4} } }
\\ &
\leq
\frac{
6 d c^2
}{
[ \min(
    \{ \bfL \}
    \cup
    \{ \bfl_i \colon i \in \N \cap [ 0, \bfL ) \}
) ]^{ \nicefrac{1}{d} }
}
+
\frac{
2
\bfL
( \lVert \bfl \rVert_{ \infty } + 1 )^\bfL
c^{ \bfL + 1 }
\max\{ 1, p \}
}{
K^{ [ ( 2 \bfL )^{-1} ( \lVert \bfl \rVert_{ \infty } + 1 )^{-2} ] }
}
\\ & \quad
+
\frac{
    5 B^2
    \bfL ( \lVert \bfl \rVert_{ \infty } + 1 )
    \max\{ p, \ln( e M ) \}
}{ M^{ \nicefrac{1}{4} } }
\end{align*}
(cf.~\cref{item:lem:measurability:3} in \cref{lem:measurability}).
\end{corollary}
\begin{cproof}{cor:SGD_full_error}
Observe that
\cref{cor:min_NN_architecture}
(with
$ ( X_j )_{ j \in \N } \leftarrow ( X^{ 0, 0 }_{ \smash{j} } )_{ j \in \N } $,
$ ( Y_j )_{ j \in \N } \leftarrow ( Y^{ 0, 0 }_{ \smash{j} } )_{ j \in \N } $,
$ \cR \leftarrow \cR^{ 0, 0 }_{ \smash{ M } } $
in the notation of \cref{cor:min_NN_architecture})
shows~\cref{eq:cor:SGD_full_error}.
\end{cproof}

\begin{corollary}
\label{cor:SGD_L1}
Let
$ d, \bfd, \bfL, M, K, N \in \N $,
$ L, a, u \in \R $,
$ b \in ( a, \infty ) $,
$ v \in ( u, \infty ) $,
$ c \in
[ \max\{ 1, L, \lvert a \rvert, \lvert b \rvert, 2 \lvert u \rvert, 2 \lvert v \rvert \}, \infty ) $,
$ B \in [ c, \infty ) $,
$ \bfl = ( \bfl_0, \bfl_1, \ldots, \bfl_\bfL ) \in \N^{ \bfL + 1 } $,
$ \bfN \subseteq \{ 0, 1, \ldots,
\allowbreak
N \} $,
$ ( \bfJ_n )_{ n \in \N } \subseteq \N $,
$ ( \gamma_n )_{ n \in \N } \subseteq \R $,
assume
$ 0 \in \bfN $,
$ \bfl_0 = d $,
$ \bfl_\bfL = 1 $,
and
$ \bfd \geq \sum_{i=1}^{\bfL} \bfl_i( \bfl_{ i - 1 } + 1 ) $,
let
$ ( \Omega, \cF, \P ) $
be a probability space,
let
$ X^{ k, n }_{ \smash{j} }
\colon \Omega \to [ a, b ]^d $,
$ k, n, j \in \N_0 $,
and
$ Y^{ k, n }_{ \smash{j} }
\colon \Omega \to [ u, v ] $,
$ k, n, j \in \N_0 $,
be functions,
assume that
$ ( X^{ 0, 0 }_{ \smash{j} }, Y^{ 0, 0 }_{ \smash{j} } ) $,
$ j \in \{ 1, 2, \ldots, M \} $,
are i.i.d.\ random variables,
let
$ \cE \colon [ a, b ]^d \to [ u, v ] $
satisfy
$ \P $-a.s.\
that
$ \cE( X_{ \smash{1} }^{ 0, 0 } )
= \E[ Y_{ \smash{1} }^{ 0, 0 } \vert X_{ \smash{1} }^{ 0, 0 } ] $,
assume for all
$ x, y \in [ a, b ]^d $
that
$ \lvert \cE( x ) - \cE( y ) \rvert \leq L \lVert x - y \rVert_{ 1 } $,
let
$ \Theta_{ k, n } \colon \Omega \to \R^{\bfd} $,
$ k, n \in \N_0 $,
and
$ \bfk \colon \Omega \to ( \N_0 )^2 $
be random variables,
assume
$ \bigl( \bigcup_{ k = 1 }^{ \infty }
\Theta_{ k, 0 }( \Omega ) \bigr)
\subseteq [ -B, B ]^\bfd $,
assume that
$ \Theta_{ k, 0 } $,
$ k \in \{ 1, 2, \ldots, K \} $,
are i.i.d.,
assume that
$ \Theta_{ 1, 0 } $ is continuous uniformly distributed on $ [ -c, c ]^\bfd $,
let
$ \cR^{ k, n }_{ \smash{J} } \colon \R^{\bfd} \times \Omega \to [ 0, \infty ) $,
$ k, n, J \in \N_0 $,
and
$ \cG^{ k, n } \colon \R^{\bfd} \times \Omega \to \R^{\bfd} $,
$ k, n \in \N $,
satisfy for all
$ k, n \in \N $,
$ \omega \in \Omega $,
$ \theta \in
\{ \vartheta \in \R^{\bfd} \colon
( \cR^{ k, n }_{ \smash{ \bfJ_n } } ( \cdot, \omega ) \colon \R^{\bfd} \to [ 0, \infty )
\text{ is differentiable at } \vartheta ) \} $
that
$ \cG^{ k, n }( \theta, \omega )
=
( \nabla_\theta \cR^{ k, n }_{ \smash{ \bfJ_n } } )( \theta, \omega ) $,
assume for all
$ k, n \in \N $
that
$ \Theta_{ k, n } = \Theta_{ k, n -1 } - \gamma_n \cG^{ k, n }( \Theta_{ k, n -1 } ) $,
and
assume for all
$ k, n \in \N_0 $,
$ J \in \N $,
$ \theta \in \R^{\bfd} $,
$ \omega \in \Omega $
that
\begin{align}
\cR^{ k, n }_{ \smash{J} }( \theta, \omega )
=
\frac{1}{J}
\biggl[
\smallsum_{j=1}^{J}
    \lvert \clippedNN{\theta}{\bfl}{u}{v}( X^{ k, n }_{ \smash{j} }( \omega ) ) - Y^{ k, n }_{ \smash{j} }( \omega ) \rvert^2
\biggr]
\qquad
\text{and}
\\
\bfk( \omega ) \in
\argmin\nolimits_{ ( l, m ) \in \{ 1, 2, \ldots, K \} \times \bfN, \, \lVert \Theta_{ l, m }( \omega ) \rVert_{ \infty } \leq B }
\cR^{ 0, 0 }_{ \smash{ M } }( \Theta_{ l, m }( \omega ), \omega )
\end{align}
(cf.~\cref{def:clipped_NN,def:p-norm}).
Then
\begin{equation}
\begin{split}
&
\E\Bigl[
\medint{[ a, b ]^d}
        \lvert \clippedNN{\Theta_{\bfk}}{\bfl}{u}{v}( x ) - \cE( x ) \rvert
    \, \P_{ X^{ 0, 0 }_{ \smash{1}\vphantom{x} } }( \dd x )
\Bigr]
\leq
\frac{
2
[ ( v - u )
\bfL
( \lVert \bfl \rVert_{ \infty } + 1 )^\bfL
c^{ \bfL + 1 }
]^{ \nicefrac{1}{2} }
}{
K^{ [ ( 2 \bfL )^{-1} ( \lVert \bfl \rVert_{ \infty } + 1 )^{-2} ] }
}
\\ & \quad
+
\frac{
3 d L ( b - a )
}{
[ \min\{ \bfL, \bfl_1, \bfl_2, \ldots, \bfl_{ \bfL - 1 } \} ]^{ \nicefrac{1}{d} }
}
+
\frac{
    3 \max\{ 1, v - u \}
    ( \lVert \bfl \rVert_{ \infty } + 1 )
    [ 2 \bfL \ln( 3 M B c ) ]^{ \nicefrac{1}{2} }
}{ M^{ \nicefrac{1}{4} } }
\\ &
\leq
\frac{
6 d c^2
}{
[ \min\{ \bfL, \bfl_1, \bfl_2, \ldots, \bfl_{ \bfL - 1 } \} ]^{ \nicefrac{1}{d} }
}
+
\frac{
    5 B^2
    \bfL ( \lVert \bfl \rVert_{ \infty } + 1 )
    \ln( e M )
}{ M^{ \nicefrac{1}{4} } }
+
\frac{
2
\bfL
( \lVert \bfl \rVert_{ \infty } + 1 )^\bfL
c^{ \bfL + 1 }
}{
K^{ [ ( 2 \bfL )^{-1} ( \lVert \bfl \rVert_{ \infty } + 1 )^{-2} ] }
}
\end{split}
\end{equation}
(cf.~\cref{item:lem:measurability:3} in \cref{lem:measurability}).
\end{corollary}
\begin{cproof2}{cor:SGD_L1}
Note that
Jensen's inequality
implies that
\begin{equation}
\E\Bigl[
\medint{[ a, b ]^d}
        \lvert \clippedNN{\Theta_{\bfk}}{\bfl}{u}{v}( x ) - \cE( x ) \rvert
    \, \P_{ X^{ 0, 0 }_{ \smash{1}\vphantom{x} } }( \dd x )
\Bigr]
\leq
\E\Bigl[
\Bigl( \medint{[ a, b ]^d}
        \lvert \clippedNN{\Theta_{\bfk}}{\bfl}{u}{v}( x ) - \cE( x ) \rvert^2
    \, \P_{ X^{ 0, 0 }_{ \smash{1}\vphantom{x} } }( \dd x )
\Bigr)^{ \! \nicefrac{1}{2} \, }
\Bigr]
.
\end{equation}
This and
\cref{cor:SGD_full_error}
(with
$ p \leftarrow 1 $
in the notation of \cref{cor:SGD_full_error})
complete
\end{cproof2}

\begin{corollary}
\label{cor:SGD_simplfied}
Let
$ d, \bfd, \bfL, M, K, N \in \N $,
$ L \in \R $,
$ c \in
[ \max\{ 2, L \}, \infty ) $,
$ B \in [ c, \infty ) $,
$ \bfl = ( \bfl_0, \bfl_1, \ldots, \bfl_\bfL ) \in \N^{ \bfL + 1 } $,
$ \bfN \subseteq \{ 0, 1, \ldots, N \} $,
$ ( \bfJ_n )_{ n \in \N } \subseteq \N $,
$ ( \gamma_n )_{ n \in \N } \subseteq \R $,
assume
$ 0 \in \bfN $,
$ \bfl_0 = d $,
$ \bfl_\bfL = 1 $,
and
$ \bfd \geq \sum_{i=1}^{\bfL} \bfl_i( \bfl_{ i - 1 } + 1 ) $,
let
$ ( \Omega, \cF, \P ) $
be a probability space,
let
$ X^{ k, n }_{ \smash{j} }
\colon \Omega \to [ 0, 1 ]^d $,
$ k, n, j \in \N_0 $,
and
$ Y^{ k, n }_{ \smash{j} }
\colon \Omega \to [ 0, 1 ] $,
$ k, n, j \in \N_0 $,
be functions,
assume that
$ ( X^{ 0, 0 }_{ \smash{j} }, Y^{ 0, 0 }_{ \smash{j} } ) $,
$ j \in \{ 1, 2, \ldots, M \} $,
are i.i.d.\ random variables,
let
$ \cE \colon [ 0, 1 ]^d \to [ 0, 1 ] $
satisfy
$ \P $-a.s.\
that
$ \cE( X_{ \smash{1} }^{ 0, 0 } )
= \E[ Y_{ \smash{1} }^{ 0, 0 } \vert X_{ \smash{1} }^{ 0, 0 } ] $,
assume for all
$ x, y \in [ 0, 1 ]^d $
that
$ \lvert \cE( x ) - \cE( y ) \rvert \leq L \lVert x - y \rVert_{ 1 } $,
let
$ \Theta_{ k, n } \colon \Omega \to \R^{\bfd} $,
$ k, n \in \N_0 $,
and
$ \bfk \colon \Omega \to ( \N_0 )^2 $
be random variables,
assume
$ \bigl( \bigcup_{ k = 1 }^{ \infty }
\Theta_{ k, 0 }( \Omega ) \bigr)
\subseteq [ -B, B ]^\bfd $,
assume that
$ \Theta_{ k, 0 } $,
$ k \in \{ 1, 2, \ldots, K \} $,
are i.i.d.,
assume that
$ \Theta_{ 1, 0 } $ is continuous uniformly distributed on $ [ -c, c ]^\bfd $,
let
$ \cR^{ k, n }_{ \smash{J} } \colon \R^{\bfd} \times \Omega \to [ 0, \infty ) $,
$ k, n, J \in \N_0 $,
and
$ \cG^{ k, n } \colon \R^{\bfd} \times \Omega \to \R^{\bfd} $,
$ k, n \in \N $,
satisfy for all
$ k, n \in \N $,
$ \omega \in \Omega $,
$ \theta \in
\{ \vartheta \in \R^{\bfd} \colon
( \cR^{ k, n }_{ \smash{ \bfJ_n } } ( \cdot, \omega ) \colon \R^{\bfd} \to [ 0, \infty )
\text{ is differentiable at } \vartheta ) \} $
that
$ \cG^{ k, n }( \theta, \omega )
=
( \nabla_\theta \cR^{ k, n }_{ \smash{ \bfJ_n } } )( \theta, \omega ) $,
assume for all
$ k, n \in \N $
that
$ \Theta_{ k, n } = \Theta_{ k, n -1 } - \gamma_n \cG^{ k, n }( \Theta_{ k, n -1 } ) $,
and
assume for all
$ k, n \in \N_0 $,
$ J \in \N $,
$ \theta \in \R^{\bfd} $,
$ \omega \in \Omega $
that
\begin{align}
\cR^{ k, n }_{ \smash{J} }( \theta, \omega )
=
\frac{1}{J}
\biggl[
\smallsum_{j=1}^{J}
    \lvert \clippedNN{\theta}{\bfl}{u}{v}( X^{ k, n }_{ \smash{j} }( \omega ) ) - Y^{ k, n }_{ \smash{j} }( \omega ) \rvert^2
\biggr]
\qquad
\text{and}
\\
\bfk( \omega ) \in
\argmin\nolimits_{ ( l, m ) \in \{ 1, 2, \ldots, K \} \times \bfN, \, \lVert \Theta_{ l, m }( \omega ) \rVert_{ \infty } \leq B }
\cR^{ 0, 0 }_{ \smash{ M } }( \Theta_{ l, m }( \omega ), \omega )
\end{align}
(cf.~\cref{def:clipped_NN,def:p-norm}).
Then
\begin{align*}
&
\E\Bigl[
\medint{[ 0, 1 ]^d}
        \lvert \clippedNN{\Theta_{\bfk}}{\bfl}{u}{v}( x ) - \cE( x ) \rvert
    \, \P_{ X^{ 0, 0 }_{ \smash{1}\vphantom{x} } }( \dd x )
\Bigr]
\\ &
\leq
\frac{
3 d L
}{
[ \min\{ \bfL, \bfl_1, \bfl_2, \ldots, \bfl_{ \bfL - 1 } \} ]^{ \nicefrac{1}{d} }
}
+
\frac{
    3
    ( \lVert \bfl \rVert_{ \infty } + 1 )
    [ 2 \bfL \ln( 3 M B c ) ]^{ \nicefrac{1}{2} }
}{ M^{ \nicefrac{1}{4} } }
+
\frac{
2
[ \bfL
( \lVert \bfl \rVert_{ \infty } + 1 )^\bfL
c^{ \bfL + 1 }
]^{ \nicefrac{1}{2} }
}{
K^{ [ ( 2 \bfL )^{-1} ( \lVert \bfl \rVert_{ \infty } + 1 )^{-2} ] }
}
\\ & \yesnumber
\leq
\frac{
d c^3
}{
[ \min\{ \bfL, \bfl_1, \bfl_2, \ldots, \bfl_{ \bfL - 1 } \} ]^{ \nicefrac{1}{d} }
}
+
\frac{
    B^3
    \bfL ( \lVert \bfl \rVert_{ \infty } + 1 )
    \ln( e M )
}{ M^{ \nicefrac{1}{4} } }
+
\frac{
\bfL
( \lVert \bfl \rVert_{ \infty } + 1 )^\bfL
c^{ \bfL + 1 }
}{
K^{ [ ( 2 \bfL )^{-1} ( \lVert \bfl \rVert_{ \infty } + 1 )^{-2} ] }
}
\end{align*}
(cf.~\cref{item:lem:measurability:3} in \cref{lem:measurability}).
\end{corollary}
\begin{cproof}{cor:SGD_simplfied}
Observe that
\cref{cor:SGD_L1}
(with
$ a \leftarrow 0 $,
$ u \leftarrow 0 $,
$ b \leftarrow 1 $,
$ v \leftarrow 1 $
in the notation of \cref{cor:SGD_L1}),
the facts that
$ B \geq c \geq \max\{ 2, L \} $
and
$ M \geq 1 $,
and~\cref{item:lem:ln_estimate:2}
in
\cref{lem:ln_estimate}
show that
\begin{align*}
&
\E\Bigl[
\medint{[ 0, 1 ]^d}
        \lvert \clippedNN{\Theta_{\bfk}}{\bfl}{u}{v}( x ) - \cE( x ) \rvert
    \, \P_{ X^{ 0, 0 }_{ \smash{1}\vphantom{x} } }( \dd x )
\Bigr]
\\ &
\leq
\frac{
3 d L
}{
[ \min\{ \bfL, \bfl_1, \bfl_2, \ldots, \bfl_{ \bfL - 1 } \} ]^{ \nicefrac{1}{d} }
}
+
\frac{
    3
    ( \lVert \bfl \rVert_{ \infty } + 1 )
    [ 2 \bfL \ln( 3 M B c ) ]^{ \nicefrac{1}{2} }
}{ M^{ \nicefrac{1}{4} } }
+
\frac{
2
[ \bfL
( \lVert \bfl \rVert_{ \infty } + 1 )^\bfL
c^{ \bfL + 1 }
]^{ \nicefrac{1}{2} }
}{
K^{ [ ( 2 \bfL )^{-1} ( \lVert \bfl \rVert_{ \infty } + 1 )^{-2} ] }
}
\\ &
\leq
\frac{
d c^3
}{
[ \min\{ \bfL, \bfl_1, \bfl_2, \ldots, \bfl_{ \bfL - 1 } \} ]^{ \nicefrac{1}{d} }
}
+
\frac{
    ( \lVert \bfl \rVert_{ \infty } + 1 )
    [ 23 B \bfL \ln( e M ) ]^{ \nicefrac{1}{2} }
}{ M^{ \nicefrac{1}{4} } }
+
\frac{
[ \bfL
( \lVert \bfl \rVert_{ \infty } + 1 )^\bfL
c^{ 2 \bfL + 2 }
]^{ \nicefrac{1}{2} }
}{
K^{ [ ( 2 \bfL )^{-1} ( \lVert \bfl \rVert_{ \infty } + 1 )^{-2} ] }
}
\\ & \yesnumber
\leq
\frac{
d c^3
}{
[ \min\{ \bfL, \bfl_1, \bfl_2, \ldots, \bfl_{ \bfL - 1 } \} ]^{ \nicefrac{1}{d} }
}
+
\frac{
    B^3
    \bfL ( \lVert \bfl \rVert_{ \infty } + 1 )
    \ln( e M )
}{ M^{ \nicefrac{1}{4} } }
+
\frac{
\bfL
( \lVert \bfl \rVert_{ \infty } + 1 )^\bfL
c^{ \bfL + 1 }
}{
K^{ [ ( 2 \bfL )^{-1} ( \lVert \bfl \rVert_{ \infty } + 1 )^{-2} ] }
}
.
\end{align*}
\end{cproof}

\section*{Acknowledgements}

This work has been
funded by the Deutsche Forschungsgemeinschaft (DFG, German Research Foundation) under Germany's Excellence Strategy EXC 2044-390685587, Mathematics M\"unster: Dy\-nam\-ics-Geometry-Structure,
by the Swiss National Science Foundation (SNSF) under the project
``Deep artificial neural network approximations for stochastic partial differential equations: Algorithms and convergence proofs''
(project number 184220),
and through
the ETH Research Grant \mbox{ETH-47 15-2}
``Mild stochastic calculus and numerical approximations for nonlinear stochastic evolution equations with L\'evy noise''.

\printbibliography

\end{document}